\documentclass[a4paper, final,leqno,10pt]{amsart}
\usepackage{amssymb}
\usepackage{epsfig}
\usepackage{enumitem}
\usepackage{framed}
\usepackage{graphicx}
\usepackage{subfig}
\usepackage{float}
\usepackage{pdflscape}
\usepackage[titletoc]{appendix}
\usepackage{mathtools} 
\usepackage{csquotes}
\RequirePackage{fix-cm}

\newcommand{\bA}{{\bf A}}
\newcommand{\bb}{{\bf b}}
\newcommand{\bB}{{\bf B}}

\newcommand{\be}{{\bf e}}
\newcommand{\bfn}{{\bf f}}
\newcommand{\bH}{{\bf H}}
\newcommand{\bI}{{\bf I}}

\newcommand{\bM}{{\bf M}}
\newcommand{\bn}{{\bfa n}}
\newcommand{\bo}{{\bf 0}}

\newcommand{\bT}{{\bf T}}

\newcommand{\bU}{{\bf U}}
\newcommand{\bv}{{\bf v}}
\newcommand{\bV}{{\bf V}}
\newcommand{\bw}{{\bf w}}

\newcommand{\bx}{{\bfa x}}
\newcommand{\by}{{\bf y}}

\newcommand{\Rn}{\mathbb{R}}

\newcommand{\bbeta}{\boldsymbol{\beta}}

\newcommand{\balp}{\boldsymbol{\alpha}}
\newcommand{\pder}[2]{\frac{\partial#1}{\partial#2}}

\newcommand{\pderi}[1]{\frac{\partial}{\partial#1}}
\newcommand{\spdr}[2]{\frac{\partial^2#1}{\partial#2^2}}
\newcommand{\bfa}[1]{\mbox{\boldmath $ #1 $}}
\usepackage[sort&compress]{natbib}

\begin{document}

 \title[Mode Isolation of RDSs]
{A computational approach for mode isolation for reaction-diffusion systems on arbitrary geometries}\author[L. Murphy, C. Venkataraman and A. Madzvamuse]{Laura Murphy\and Chandrasekhar Venkataraman \and Anotida Madzvamuse}
  \address[L.~Murphy]{University of Sussex, School of Mathematical and Physical Sciences, Department of Mathematics, Brighton BN1 9QH, UK.}
  \email[L.~Murphy]{L.R.Murphy@sussex.ac.uk}
  \address[C.~Venkataraman]{Mathematical Institute, North Haugh, University of St Andrews, Fife, UK. KY16 9SS.}
  \email[C.~Venkataraman]{cv28@st-andrews.ac.uk}
  \address[A.~Madzvamuse]{University of Sussex, School of Mathematical and Physical Sciences, Department of Mathematics, Brighton BN1 9QH, UK.}
  \email[A.~Madzvamuse]{A.Madzvamuse@sussex.ac.uk}

\date{}

\begin{abstract}In this article we present a computational framework for isolating spatial patterns arising in the steady states of reaction-diffusion systems. Such systems have been used to model many different phenomena in areas such as developmental and cancer biology, cell motility and material science. Often one is interested in identifying parameters which will lead to a particular pattern. To attempt to answer this, we compute eigenpairs of the Laplacian on a variety of domains and use linear stability analysis to determine parameter values for the system that will lead to spatially inhomogeneous steady states whose patterns correspond to particular eigenfunctions. This method has previously been used on domains and surfaces where the eigenvalues and eigenfunctions are found analytically in closed form. Our contribution to this methodology is that we numerically compute eigenpairs on arbitrary domains and surfaces. Here we present various examples and demonstrate that mode isolation is straightforward especially for low eigenvalues. Additionally we see that if two or more eigenvalues are in a permissible range then the inhomogeneous steady state can be a linear combination of  the respective eigenfunctions. Finally we show an example which suggests that pattern formation is robust on similar surfaces in cases that the surface either has or does not have a boundary.

\end{abstract}
\maketitle
\section{Introduction}
In his seminal work, \cite{turing1952} presented an elegant mathematical theory of reaction-diffusion type for pattern formation in developmental biology. He showed that, via a symmetry breaking, a homogeneous state which is linearly stable in the absence of diffusion may be driven unstable in the presence of diffusion to give rise to the emergence of a spatially inhomogeneous pattern. This process is now well known as {\it diffusion-driven instability} or {\it Turing instability}. Since then, reaction-diffusion systems have been proposed and applied to model many phenomena including cancer invasion and angiogenesis in cancer biology \citep{chap,chaplain,gatenby}, pattern formation in developmental biology \citep{hunding,maini}, wound healing in biomedicine \citep{dale,sherratt}, cell motility \citep{mogilner,moek,uduak} and material science \citep{bozzini,krinsky} among many others.  Despite their numerous applications, Turing's theory of pattern formation has been widely criticised mainly due to the lack of robustness of the model system to changes in the parameters as well as the lack of experimental evidence of the existence of so-called morphogens with varying diffusivities. Only recently has the existence of chemical morphogens been experimentally validated in hair follicle pattern formation by \cite{sick2006}. 

To-date mode selection and parameter identification for reaction-diffusion systems have been mainly carried out on regular planar domains and surfaces where the eigenvalue problem can be analytically solved to yield analytical forms of the wave numbers as well as their corresponding eigenfunctions \citep{ano,madzvamuse2003,uduak}. In this work, we will depart from this framework and extend computationally mode selection and parameter identification to include arbitrary domains and stationary surfaces. First, we will solve the eigenvalue problem numerically using finite elements on planar domains or surface finite elements on smooth surfaces, respectively, to obtain the eigenmodes and their corresponding eigenfunctions. Here, we employ the Krylov-Schur algorithm \citep{stewart} for solving the resulting algebraic system arising from the finite element discretisation. Second, we then pick an eigenmode to which we apply the necessary and sufficient conditions for Turing diffusion-driven instability in order to isolate reaction-kinetic model parameter values within a reaction-diffusion system. This process can be loosely thought of as an inverse problem for model parameter identification. Once the parameter values are isolated, the full reaction-diffusion system is then solved with these isolated parameter values to obtain an inhomogeneous spatially varying solution which is then compared to the numerically computed eigenfunction on the domain or surface.  Alternatively, one could pose the following problem to which this methodology will provide insightful information which is otherwise out of reach with the current methodology: {\it Given a biological pattern on a domain or surface and a plausible reaction-diffusion system, what are the model parameter values within this reaction-diffusion system that will give rise to the observed pattern?} This article provides a theoretical and computational framework to answer such a question.

It must be observed that the eigenvalue problem and the reaction-diffusion system are both solved by a similar numerical method, the finite element method in multi-dimensions \citep{johnson}. The finite element method is well known for its capability to deal with complex irregular geometries \citep{barreira,elliot,venkataraman}. Alternative numerical methods such as finite differences \citep{beckett}, spectral methods \citep{chap,ruuth} and finite volume methods among others could be used but with considerable efforts in dealing with geometrical complexities. As mentioned above one interpretation of our approach is that it provides a means of estimating parameter values such that the pattern predicted by linear stability analysis is close to a desired pattern. It must be noted that in many cases the steady state pattern may not be an eigenfunction (or a linear combination of the eigenfunctions) of the Laplacian on the given domain. This is since the nonlinear terms play a role in the resultant steady state pattern \citep{murray2003}. In such a setting our approach may provide parameters which serve as a suitable initial guess for a more advanced parameter identification algorithm \citep{croft,garvie}

The remainder of this article is structured as follows. In Section \ref{sec:model} we introduce the mathematical model which we study in this work. We summarise the necessary and sufficient conditions for Turing {\it diffusion-driven instability} in Section \ref{sec:conditions}. We then detail how mode selection and parameter identification are carried out. In Sections \ref{sec:isolation} and \ref{sec:genisolation} we outline the new theoretical and computational framework for mode selection and parameter identification. The use of the finite element method is described in Section \ref{sec:ferds}. We then give specific examples in 2- and 3-dimensions for regular (by which we mean domains on which analytic expressions for the eigenfunctions are available) as well as general domains and surfaces. We discuss the implications of our framework in the context of current methodologies and conclude that given a biological pattern and a reaction-diffusion system, our approach provides a useful tool for estimating parameter values which may give rise to the observed pattern.

\section{Mathematical model framework}\label{sec:model}
In order to illustrate with clarity  the novelty of our approach, we first introduce the standard theoretical framework for reaction-diffusion systems in multi-dimensions \citep{murray2003}.  Let $\Omega  \subset {\mathbb{R}}^m$ $(m=1,2,3)$   be a simply connected bounded stationary volume for all time $t \in I=[0,t_F]$, $t_F > 0$ and ${\partial \Omega}$ be the surface
 boundary enclosing $\Omega$. Also let $ {\bfa  u} = \left( u \left( {\bfa x},t \right), v \left( {\bfa x},t \right) \right)^T$ be a vector of two chemical concentrations at position ${\bfa x} \in \Omega \subset  {\mathbb{R}}^m$ and time $t\in I$. The evolution equations for reaction-diffusion systems in the absence of cross-diffusion can be obtained from the application of the law of mass conservation and the extended Fick's first law \citep{murray2003,turing1952} to yield the dimensional system 
\begin{equation} \label{eq:rden}
 \begin{cases}
 \begin{cases}
  u_t = D_u \Delta  u  + f(u,v),\\  
  
  v_t = D_v\Delta v + g (u,v), 
  \end{cases}
  \quad x \in \Omega, \; t> 0, \\
  \\
{\bfa n} \cdot \nabla u= {\bfa n} \cdot \nabla v=0,\; x\,\text{on}\,\partial\Omega, \; t \ge 0, \\ \\
u(x,0)= u_0(x),\; \text{and} \; v(x,0)=v_0(x),\; x\,\text{on}\,\Omega, \; t=0,
 \end{cases}
\end{equation}
where $\Delta$ denotes the usual cartesian Laplace operator, $D_u>0$ and $D_v>0$ are diffusion coefficients. Here, {\bfa n} is the unit outward normal to $\partial\Omega$.  Initial conditions are prescribed through non-negative bounded functions $u_0 (x)$ and $v_0 (x)$.  In the above, $f (u,v)$ and $g (u,v)$ represent nonlinear reactions. 

In the case of surfaces, the Laplace operator is replaced by the Laplace-Beltrami operator $\Delta_\Gamma$, where $\Gamma$ is the (smooth) surface. This can be described as follows (For more details we refer the interested reader to see \cite{dziukelliott}). If $f:\Gamma\to\Rn$ is differentiable at $x\in\Gamma$ we can define the {\it tangential gradient} of $f$ at $x\in\Gamma$ by
\begin{equation}
 \nabla_\Gamma f =\nabla\bar{f}-\nabla\bar{f}\cdot\bn\bn.
\end{equation}
Here $\bar{f}$ is a smooth extension of $f:\Gamma\to\Rn$ to an $(n+1)$-dimensional neighbourhood $U$ of the surface $\Gamma$, so that $\bar{f}|_\Gamma=f$. $\nabla$ is the gradient in $\Rn^{n+1}$ and $\bn$ is the unit normal. The {\it Laplace-Beltrami operator} applied to a twice differentiable function $f\in C^2(\Gamma)$ is given by
\begin{equation}
 \Delta_\Gamma f = \nabla_\Gamma\cdot\nabla_\Gamma f.
\end{equation}
It must be observed that if the surface does not have a boundary, no boundary conditions are needed. If the surface has a boundary, we assume homogeneous Neumann boundary conditions. 

Since the reaction terms are nonlinear, analytical solutions cannot normally be obtained. Therefore we investigate solution behaviour using linear stability theory and numerical methods. Linear stability analysis is one way of determining the behaviour of a nonlinear system near a given stationary point, normally a uniform steady state, of the given system. The idea is to find under what conditions on the nonlinear reaction kinetics is the uniform steady state linearly asymptotically stable in the absence of diffusion. When diffusion is introduced, the uniform steady state is driven unstable in what is now known as the process of {\it diffusion-driven instability} with the system converging to a spatially inhomogeneous steady state, thereby giving rise to patterning \citep{murray2003,turing1952}. The mathematical treatment of the derivation of the necessary conditions for {\it diffusion-driven} instability requires solving the well known eigenvalue problem, with $W$ a solution of
\begin{subequations}\label{eq:equn}
 \begin{gather}
  \Delta W+k^2W=0, \quad \bx\in \Omega,\\
 ({\bfa n}\cdot\nabla)W=0, \quad {\bfa x} \in \partial \Omega,  
 \end{gather}
\end{subequations}
where the solution pairs ($k$ (eigenvalues), $W_k(\bx)$ (eigenfunctions) obtained either analytically on certain spatial domains or numerically for the general case) of this vector equation can be compared to the spatially inhomogeneous steady state solutions of \eqref{eq:rden}, with good agreement expected near primary bifurcation points.

This approach is generally called mode isolation. The most famous exploration of this problem is the celebrated article "Can one hear the shape of the drum?" by Mark Kac \citeyearpar{kac}. The question being asked is if one knows all the eigenvalues of the eigenvalue problem is it possible to determine the domain? It was later proven by Gordon, Webb and Wolpert \citeyearpar{gww} that the answer is no and they gave examples of distinct regions with identical eigenvalues.

Other work concerned with mode isolation and linear stability theory for reaction-diffusion systems can be found in \cite{chap} and \cite{ano}, here the validation has been mainly restricted to special domains and volumes where the eigenvalue problem can be solved analytically. In this work we will depart from this framework, instead we will compute approximations of the eigenpairs on arbitrary, simply connected domains, volumes and surfaces. We then use these eigenvalues to calculate, by use of the Turing-parameter space restrictions, appropriate model parameter values. This approach can be thought to be analogous to an inverse parameter identification approach whereby, given the eigenvalues and eigenfunctions solving the eigenvalue problem (\ref{eq:equn}), find model parameter values that would give rise to an inhomogeneous spatially varying solution similar to that exhibited by the eigenfunction. To confirm numerical predictions, we use the computed model parameter values to solve the full nonlinear reaction-diffusion systems and compare approximated eigenfunctions on these arbitrary domains, volumes and surfaces to the spatially inhomogeneous solutions obtained numerically.

To proceed, next we show the two-component form which we will work with and state the conditions for {\it diffusion-driven} instability. These will help us to isolate particular modes.

\section{Conditions for diffusion driven instability for reaction-diffusion systems}\label{sec:conditions}
All two component reaction-diffusion systems of the form \eqref{eq:rden} can be non-dimensionalised and scaled to take the form
\begin{subequations}\label{eq:non}
 \begin{align}
  u_t = \gamma f(u,v)+\Delta u, \quad v_t = \gamma g(u,v)+d\Delta v, \quad \bx\in \Omega\subset \Rn^n, \; t\in [0,\infty], \label{eq:parta} \\
  (\bn\cdot \nabla)\begin{pmatrix}u \\ v\end{pmatrix}=0 \quad \bx \in \partial\Omega  \; t\in [0,\infty],\\
  u(\bx, 0), v(\bx, 0) \text{ given},
 \end{align}
\end{subequations}
where $u=u(\bx,t), v=v(\bx,t)$, $d$ is the ratio of diffusion coefficients, $f(u,v)$ and $g(u,v)$ describe the reaction kinetics. For simplicity, we assume that $f$ and $g$ are continuously differentiable, $\gamma$ can be described as the relative strength of the reaction terms or alternatively as the domain size.
We have zero flux boundary conditions (homogeneous Neumann) because we want only internal sources of instability, ie. self-organisation of the system.
 A uniform steady state $(u_s,v_s)$ is a fixed point where $(u,v)=(u_s,v_s)$, constant in time and space, satisfies \eqref{eq:non}, i.e. $(u_t,v_t)|_{u=u_s, v=v_s}=\bo$.
We can find the steady state by solving $f(u_s,v_s)=g(u_s,v_s)=0$. \\
The conditions for instability due to diffusion are well known (see, for example \cite{murray2003}). Firstly, in the absence of diffusion, the steady state $(u_s,v_s)$ is linearly stable if and only if the partial derivatives of $f$ and $g$ at $(u_s,v_s)$ satisfy
\begin{equation}
f_u+g_v<0 \text{ and }  f_ug_v-f_vg_u>0.
\end{equation}
Linear stability analysis considering small perturbations from the equilibrium $\bw(\bx,t)=(\hat{u}(\bx,t),\hat{v}(\bx,t))$ leads us to the system
\begin{equation}
 \bw_t=\gamma\begin{pmatrix}f_u & f_v \\ g_u & g_v\end{pmatrix}\bw+
\begin{pmatrix}1 & 0 \\ 0 & d  \end{pmatrix} \Delta\bw,
\end{equation}
which can be solved by method of separation of variables to yield 
\begin{equation}
 \bw(\bx,t)=\sum_kc_ke^{\lambda t}W_k(\bx),
\end{equation}
where $W_k(\bx)$ solve the eigenvalue problem 
\begin{subequations}
 \begin{gather}
  \Delta W+k^2W=0 \quad \bx\in \Omega  \label{equn},\\
 (\bn\cdot\nabla)W=0 \quad \bx\in \partial \Omega \label{neumann}.
 \end{gather}
\end{subequations}
These are modes that will decay with time unless the wavenumber $k^2$ satisfies
\begin{equation}\label{eq:disprel}
c(k^2)=d(k^2)^2-\gamma(df_u+g_v)k^2+\gamma^2(f_ug_v-f_vg_u)<0,
\end{equation}
this means that instability will occur if 
\begin{equation}
  df_u+g_v>0, \quad\quad (df_u+g_v)^2-4d(f_ug_v-f_vg_u)>0
\end{equation}
and $k^2$ lies in the range $k_-^2<k^2<k_+^2$ where
\begin{equation}\label{eq:kroots}
  k^2_{\pm}=\gamma\frac{(df_u+g_v)\pm\sqrt{(df_u+g_v)^2-4d(f_ug_v-f_vg_u)}}{2d}.
\end{equation}
We exploit this range to isolate particular patterns/modes. The unstable modes will correspond to the eigenfunctions of the Laplacian (or Laplace-Beltrami) on the chosen domain or surface with the selected boundary conditions and $k^2$ the associated eigenvalues. The effect of varying $d$ and $\gamma$ on \eqref{eq:disprel} is shown in Figure \ref{fig:ck2}.\\
In summary the necessary conditions for diffusion driven instability are
\begin{subequations}
 \begin{align}
  f_u+g_v&<0, &f_ug_v-f_vg_u>0, \\
  df_u+g_v&>0, &(df_u+g_v)^2-4d(f_ug_v-f_vg_u)>0.
 \end{align}
\end{subequations}
Additionally, the sufficient conditions for  patterning  formation  are  that  one must be able to isolate distinct real wave numbers and that the domain must be large enough \citep{madzvamuse2010,madzvamuse2015,murray2003}.

\subsection{Examples of reaction kinetics}
For illustrative purposes, we consider three classical reaction kinetics as summarised below. The work presented in this article holds true for other similar reaction kinetics capable of generating Turing patterns.
\subsubsection{Schnakenberg or activator-depleted substrate Kinetics}
The Schnakenberg kinetics \citep{schnak} are a condensed version of the well documented Brusselator model describing a series of autocatalytic reactions also known as activator-depleted models \citep{giermein,prig}, and can characterised by
\begin{equation}
A\rightleftharpoons X \quad\quad
B+X\to Y+D \quad\quad
2X+Y\to 3X.
\end{equation}
Using the Law of Mass Action and the non-dimensionalisation of $f$ and $g$, within system \eqref{eq:non}, we obtain that
\begin{equation}\label{eq:schnakkin}
f(u,v)=a-u+u^2v \quad \text{ and }\quad g(u,v)=b-u^2v,
\end{equation}
where $a$ and $b$ are positive parameters.
\subsubsection{Gierer-Meinhart Kinetics}
One of the models proposed by \cite{giermein} describes an system whereby an "activator" activates the production of an "inhibitor" which inhibits the production of the activator. Again the non-dimensionalised form can be obtained  
\begin{equation}\label{eq:gm}
f(u,v)=a-bu+\frac{u^2}{v(1+ku^2)}, \quad \text{and} \quad g(u,v)=u^2-v,
\end{equation}
where $a$ and $b$ are positive parameters (representing constant production rate and linear degradation respectively) and $k$ can be thought of as the saturation concentration of $u$.

\subsubsection{Thomas Kinetics}
The Thomas model \citep{thomas} is an immobilized-enzyme substrate-inhibition mechanism which can be written in  non-dimensional form as 
\begin{equation}\label{eq:tho}
f(u,v)=a-u-\frac{\rho uv}{1+u+Ku^2}, \quad
g(u,v)=\alpha b - \alpha v-\frac{\rho uv}{1+u+Ku^2},
\end{equation}
where $a$, $\rho$, $K$, $\alpha$, $\beta$ are all non-negative parameters. This can be interpreted as in \cite{murray1982} by saying that $u$ and $v$
\begin{itemize}
\item are generated by constant production $a$ and $\alpha b$ respectively,
\item decay linearly proportional to $u$ and $\alpha v$ respectively and
\item are used up in a substrate inhibition manner  $\frac{\rho uv}{1+u+Ku^2}$.
\end{itemize}

\section{Overview on mode isolation for reaction-diffusion systems}\label{sec:isolation}
The goal of mode isolation is to choose parameters, in our case ($d,\gamma$), so that a trajectory starting from a small random perturbation from the steady state will evolve into a spatial pattern generated by one that corresponds, or at least is close to, a chosen eigenfunction of the Laplacian on that domain.
Wavenumber isolation of reaction-diffusion systems is described by \cite{ano} in one dimension, squares and triangles. In \cite{uduak} wavenumbers of a visco-elastic model are isolated on the unit disk. We use similar ideas in the present work. The basic steps are as follows.
\begin{enumerate}
\item Determine a subset of eigenpairs of the Laplacian with suitable boundary conditions on the domain. For special domains this can be done analytically but in general must be done numerically.
\item Compute the dispersal relation \eqref{eq:disprel} for the chosen reaction kinetics (this is independent of the geometry) and the range of admissible wave numbers as a function of $d$ and $\gamma$.
\item Compute $d^*$ and $\gamma^*$ such that only one of the eigenvalues (wave numbers) computed in step 1 is in the range.
\item In order to compare with the patterned state, solve the reaction-diffusion system numerically with computed parameter values and compare with the numerically computed eigenfunctions.
\end{enumerate}
 It is possible to implement the above procedure simply because if a domain is bounded and the boundary is sufficiently regular, the Neumann Laplacian has a discrete spectrum of infinitely many non-negative eigenvalues with no finite accumulation point
\begin{equation}
 0<\lambda_1\leq\lambda_2\leq\cdots,\lambda_n\to\infty
\end{equation}
and this is due to the spectral theorem for compact self-adjoint operators (\citeauthor{benguria}, 2016; \citealp{kreyszig,taylor}).

The aim is to have an algorithm to find the parameter values $d$ and $\gamma$ for a given eigenpair $(k^2,W)$ such that only patterns analogous to $W$ will grow. For this, one needs that the corresponding $k$ is in the range defined in \eqref{eq:kroots}
\begin{equation}\label{eq:LR}
 \gamma L=k^2_-<k<k^2_+=\gamma R
\end{equation}
where
\begin{subequations}
\begin{align} 
 L=\frac{(df_u+g_v)-\sqrt{(df_u+g_v)^2-4d(f_ug_v-f_vg_u)}}{2d},\\ R=\frac{(df_u+g_v)+\sqrt{(df_u+g_v)^2-4d(f_ug_v-f_vg_u)}}{2d},
\end{align}
\end{subequations}
and that no other $k$ is in this range. 
\begin{figure}
\centering
\subfloat[$\gamma=15$]{\includegraphics[width=58mm]{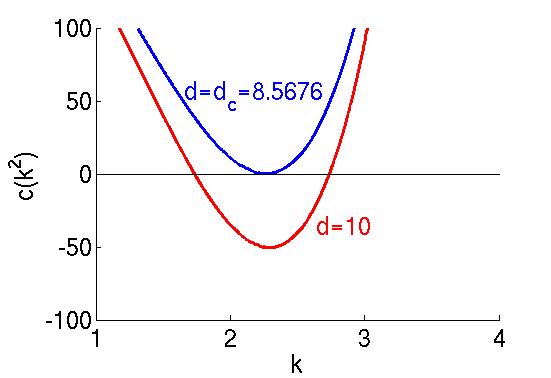}}
\subfloat[d = 10]{\includegraphics[width=58mm]{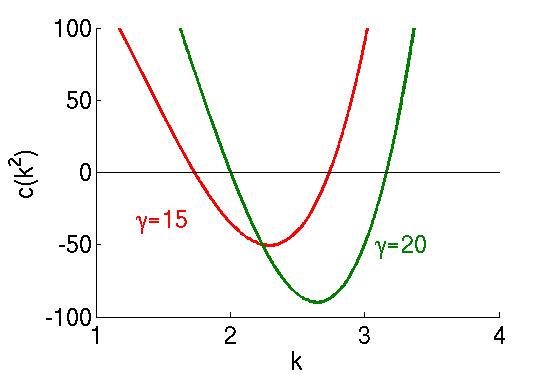}}
\caption{Here the dispersal relation \eqref{eq:disprel} is plotted (for Schnakenberg kinetics). For a fixed value of $\gamma$, when $d$ is below the critical value $d_c$, $c(k^2)$ has no roots so no modes can be isolated. As $d$ increases as does the difference between the two roots so there is more chance the value of $k$ we seek will be between $k^2_-$ and $k^2_+$. Similarly, for a fixed value of $d$, increasing $\gamma$ causes both $k^2_-$ and $k^2_+$ to increase (Colour version online)}\label{fig:ck2}
\end{figure}
In other words, the sign of the polynomial $c(k^2)$ for a given $k$ determines if the mode will grow. Figure \ref{fig:ck2} illustrates how the graph of $c(k^2)$ changes as $d$ and $\gamma$ are varied. We define the critical diffusion ratio $d_c$ as the root of
\begin{equation}
d_c^2f^2_u + 2(2f_vg_u - f_ug_v)d_c + g^2_v = 0.
\end{equation}
We find $(k^2,W)$ either analytically or numerically. Then we propose the following algorithm described in pseudo-code:

{\bf Input:} $d=d_c+\epsilon$, $\epsilon\approx d_c/5$, $\gamma>0$, $f, g$ and the $k_{l,n}$ that we wish to be uniquely isolated.
\begin{enumerate}
\item Compute $k^2_-$ and $k^2_+$ from \eqref{eq:LR}.
\item If $k^2_{l,n}<k^2_-$ increase $\gamma$ by 1 (this number is arbitrary but should be small). This moves the curve to higher values of $k$.
\item If $k^2_{l,n}<k^2_+$ decrease $\gamma$ by 1. This moves the curve to lower values of $k$.
\item If there exists another $k_{l,n}^*\not= k_{l,n}$ such that $k^2_-<k_{l,n}^{*2}<k^2_+$ then decrease $\epsilon$ by $d_c/100$. This shifts the curve upwards so the difference between $k^2_-$ and $k^2_+$ is smaller.
\item If $k_{l,n}$ is uniquely isolated END. If not go to 3.
\end{enumerate}
{\bf Output:} The appropriate $d,\gamma$.\\
Note that we cannot have $d<d_c$ (because then $c(k^2)$ would have no roots) nor $\gamma<0$ (because $k^2>0$).

\section{Finite element method for reaction diffusion systems}\label{sec:ferds}

In order to validate that our mode isolation algorithm does indeed isolate the desired unstable mode, we will simulate the reaction-diffusion systems under consideration with the computed parameter values. To do this we employ a finite element method for the space discretisation and an implicit-explicit time-stepping scheme for the temporal approximation \citep{lakkis2013,madzvamuse2006,ruuth}. 
 
In order to compute a finite element approximation, we write the weak formulation of \eqref{eq:non} as follows:
Find $u,v\in L^2(0,T;H^1(\Omega))$ such that for all $\phi\in H^1(\Omega)$ we have 
 \begin{align}\label{eqn:wf}
\begin{cases} \int_\Omega u_t\phi + \int_\Omega\nabla u \cdot \nabla\phi&=\gamma\int_\Omega f(u,v)\phi,\\
\int_\Omega v_t \phi+d\int_\Omega\nabla v \cdot \nabla\phi&=\gamma\int_\Omega g(u,v)\phi,
\end{cases}\quad \bx\in\Omega,\;t>0.
\end{align}
In this work we shall assume the well posedness of the weak formulation above. We note that for suitable parameter values existence and uniqueness of a classical solution, and hence a weak solution, to \eqref{eq:non} may be shown for example by the method of invariant regions proposed and analysed by Sm\"{o}ller \citeyearpar{smoller}.

\subsection{Spatial discretisation}

We define the computational domain $\Omega_h$ by requiring that $\Omega_h$ is a polyhedral approximation to $\Omega$.  We define $T_h$ to be a triangulation of $\Omega_h$ made up of non-degenerate elements $\kappa_i$, i.e., $T_h=\bigcup_i\{\kappa_i\}$. We define the finite element space $V_h := \{ v_h\in C^0(\Omega):v_h|_\kappa \text{ is linear}\}$. The semidiscrete (space discrete) finite element approximation to (\ref{eqn:wf}) seeks a pair $(U,V)\in V_h^2$ such that
 \begin{align}\label{eqn:fe_scheme}
  \begin{cases} \int_{\Omega_h}  U_t\phi+ \int_{\Omega_h}\nabla U \cdot \nabla\phi = \gamma\int_{\Omega_h} I_h\left[f(U,V)\right]\phi,\\
\int_{\Omega_h} V_t\phi+ d\int_{\Omega_h}\nabla V \cdot \nabla\phi = \gamma\int_{\Omega_h} I_h\left[g(U,V)\right]\phi,
  \end{cases} \forall \phi \in V_h,
 \end{align}
  where we use the Lagrange interpolant of the initial data into $V_h$ as initial conditions for the scheme.
In order to illustrate a concrete example of the scheme, we focus on the reaction-diffusion system with Schnakenberg kinetics (\ref{eq:schnakkin}).
The finite element approximation (\ref{eqn:fe_scheme}) with the Schnakenberg kinetics can be written in matrix-vector form as follows
\begin{subequations}\label{eqn:mv_form_schnak}
 \begin{align}
 \bM \balp_t+\bA \balp=\gamma\left[a\bH-\bM \balp+ \bM (\balp)^2 \bbeta\right],\\
 \bM  \bbeta_t+d\bA \bbeta=\gamma\left[b\bH- \bM (\balp)^2 \bbeta\right],
 \end{align}
\end{subequations}
where $\balp$ and $\bbeta$ are the coefficient vectors of the finite element functions $U$ and $V$ respectively  and
\[
M_{i,j}=\int_{\Omega_h} \phi_i \phi_j, \quad A_{i,j}=\int_{\Omega_h} \nabla \phi_i\cdot\nabla\phi_j \quad \text{and} \quad H_j=\int_{\Omega_h}\phi_j.
\]

\subsection{Temporal discretisation}
\begin{figure}
\centering
\includegraphics[width=115mm]{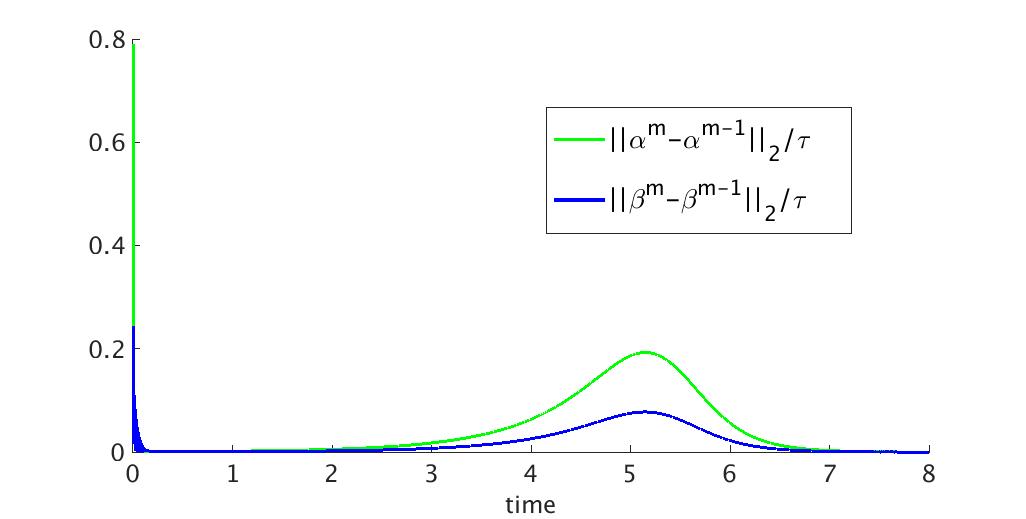}
\caption{Plot of the $L_2$ norm  of the discrete time-derivative over time for the example shown in Figure 8(b). There is an initial decay due to diffusion followed by a growth because of the exponentially growing modes which eventually decays, due to the dominant nonlinear terms (Colour version online)}
\label{l2si}
\end{figure}

 For the temporal discretisation we employ an IMEX method \citep{lakkis2013,madzvamuse2006,ruuth} in which the diffusive term is treated implicitly and the reaction terms are treated explicitly, for simplicity we employ a uniform timestep $\tau$. Introducing the shorthand for a time discrete sequence of functions, $f^n=f(t_n)$, 
the fully discrete scheme we employ reads, for $n=0,1,\dots$, given $(U^n, V^n)\in V_h^2$ find $(U^{n+1}, V^{n+1})\in V_h^2$
such that, $\forall \phi \in V_h$,
 \begin{align}\label{eqn:discrete_fe_scheme}
  \begin{cases} \int_{\Omega_h}  \frac{1}{\tau}\left(U^{n+1}-U^n\right)\phi+ \int_{\Omega_h}\nabla U^{n+1} \cdot \nabla\phi = \gamma\int_{\Omega_h} I_h\left[f(U^n,V^n)\right]\phi,\\
 \int_{\Omega_h}  \frac{1}{\tau}\left(V^{n+1}-V^n\right)\phi+ d\int_{\Omega_h}\nabla V^{n+1} \cdot \nabla\phi = \gamma\int_{\Omega_h} I_h\left[g(U^n,V^n)\right]\phi,
  \end{cases}
 \end{align}
 where we use Lagrange interpolant of the initial data into $V_h$ as initial conditions for the scheme.
This leads us to the following matrix vector form
\begin{subequations}\label{eq:matrixsystem}
 \begin{align}
 \left(\frac{1}{\tau}\bM+\bA\right)  \balp^{m+1}=\gamma\left[a\bH-\bM \balp^m+\bM (\balp^m)^2\bbeta^m \right]+\frac{1}{\tau}\bM \balp^m,\\
 \left(\frac{1}{\tau}\bM+d\bA\right) \bbeta^{m+1} =\gamma\left[b\bH- \bM (\balp^m)^2\bbeta^m\right]+\frac{1}{\tau}\bM \bbeta^m.
 \end{align}
\end{subequations}
 Since we are interested in convergence to a spatially inhomogeneous steady state, for the stopping criteria we use the $L_2$ norm of the approximate time derivative of the discrete solution, stopping the computation if this decreases below some tolerance (see Figure \ref{l2si}). 

\subsection{Numerical computations}\label{sec:parameters}
We take the parameter values as shown in Table \ref{table:parameters}, for the initial data we use small quasi-random perturbations around the uniform steady state values. The linear system \eqref{eq:matrixsystem} is solved using the conjugate gradient method \citep{dealii,golub,CG}.
\begin{table}[H]
 \caption{Parameters for reaction kinetic models and the corresponding uniform steady states. The uniform states for Schnakenberg kinetics were obtained analytically while for the Gierer-Meinhardt and Thomas reaction kinetics these were calculated computationally using the Newton-Raphson method \citep{arfken,ano}.}
 \label{table:parameters}
 \begin{tabular}{lcccccccc}
 \hline\noalign{\smallskip}
 Model & a & b & k & K & $\alpha$ & $\rho$ & $u_s$ & $v_s$ \\ 
\noalign{\smallskip}\hline\noalign{\smallskip}
 Schnakenberg & 0.9 & 0.1 & & & & & 1 & 0.9 \\ 
 Gierer-Meinhart & 0.1 & 1 & 0.5 & & & & 0.8395 & 0.7047 \\ 
 Thomas & 150 & 100 & & 0.05 & 1.5 & 13 & 37.74 & 25.16 \\ 
\noalign{\smallskip}\hline
 \end{tabular}
\end{table}

\subsection{Convergence to a steady state}
Figure \ref{l2si} plots the $L_2$ norm of the discrete time derivative of $U$ and $V$ against the elapsed time. To begin with the difference is large. This quickly decays due to diffusion then there is a rapid growth, because of the exponentially growing modes. The time derivative eventually starts to decay due to the effects of the nonlinear terms that act to bound the exponentially growing solution thereby giving rise to a spatially inhomogeneous steady state. 

\section{Isolating modes on general domains}\label{sec:genisolation}
On arbitrary domains, analytical solutions for the eigenvalue problem are not typically available  but approximate eigenpairs can be computed numerically. Numerically approximating these pairs is a significant challenge. In general, as we are only typically interested in a small number of eigenpairs, it is not necessary to find all solution pairs, however for our approach to mode isolation to remain applicable, it is important that we obtain consecutive pairs.\\
As previously stated, the eigenvalue problem we wish to solve is as follows,
\begin{equation}\label{eq:geneig}
 \begin{cases}
  \Delta W+k^2W=0, \quad \bx\in\Omega, \\
 (\bn\cdot\nabla)W=0, \quad \bx\in \partial\Omega.
 \end{cases}
\end{equation}
To approximate the solution we  employ the finite element method for the spatial discretisation outlined in Section \ref{sec:ferds}. We work with the weak formulation of the eigenvalue problem and look for an approximate eigenpairs $(W_h,k_h^2)\in V_h\times\mathbb{R}_+$ (where $V_h$ contains all continuous piecewise linear functions on a given mesh) such that\begin{equation}
\int_\Omega \nabla W_h\cdot\nabla\phi=k^2\int_\Omega W_h\cdot\phi, \qquad \forall \phi \in V_h.
\end{equation}
As in \eqref{eqn:mv_form_schnak} this may be written in matrix-vector form, we want to find $(\balp,k_h^2)\in\mathbb{R}^m\times\mathbb{R}_+$, where $m$ is the dimension of $V_h$ such that
\begin{equation}
\bA\balp = k^2 \bM\balp, \label{eq:ep}
\end{equation}
where $\bA$ and $\bM$ are  stiffness and mass matrices defined respectively, by
\begin{equation}
A_{i,j}=\int_{\Omega_h} \nabla \phi_i\cdot\nabla\phi_j \quad \text{and} \quad M_{i,j}=\int_{\Omega_h} \phi_i \phi_j.
\end{equation}
This is a generalised eigenvalue problem. 
We use the package {\bf deal.II} \citep{dealii} for its approximation using SLEPc and the Krylov-Schur algorithm. For completeness we give a  description of the algorithm employed in Appendix \ref{krsch}.

\begin{figure}
\centering
\subfloat[Unit sphere]{\includegraphics[width=28mm]{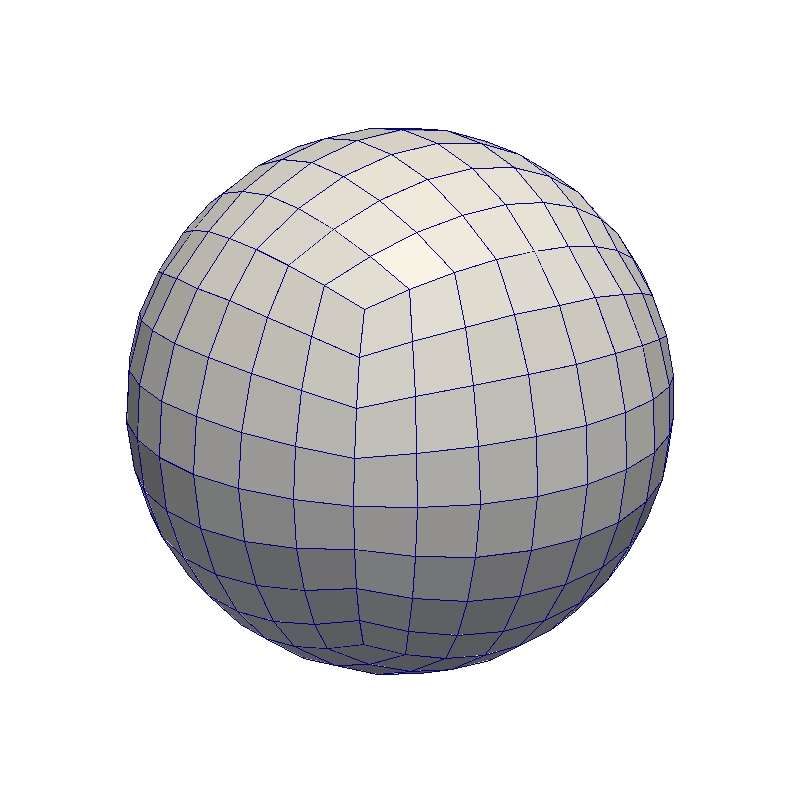}\label{fig:sphmesh}}\;
\subfloat[Unit sphere cut to show inside]{\includegraphics[width=28mm]{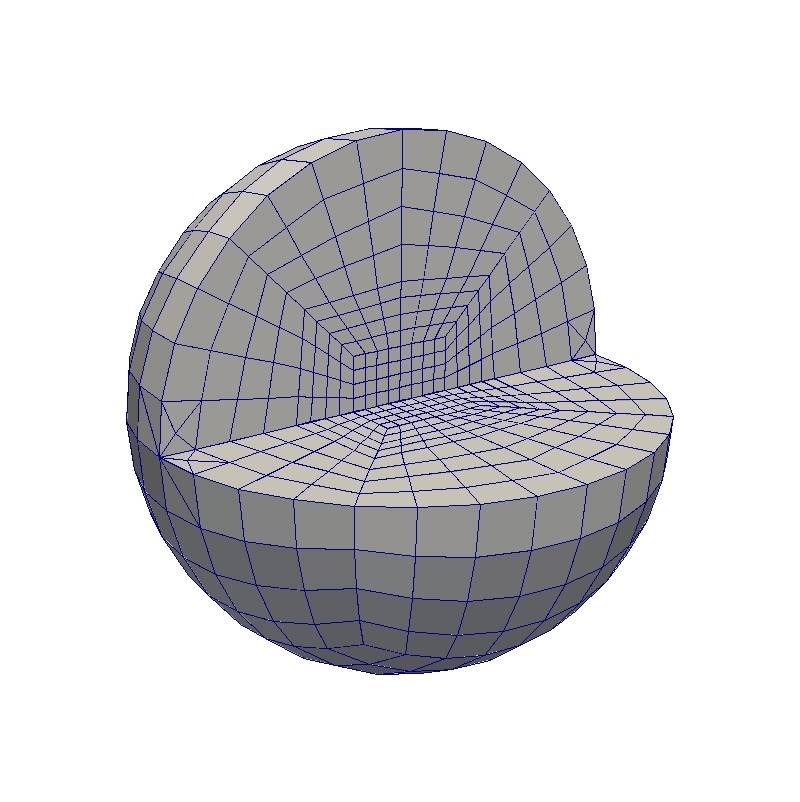}\label{fig:sphcut}}\;
\subfloat[Surface of unit sphere]{\includegraphics[width=28mm]{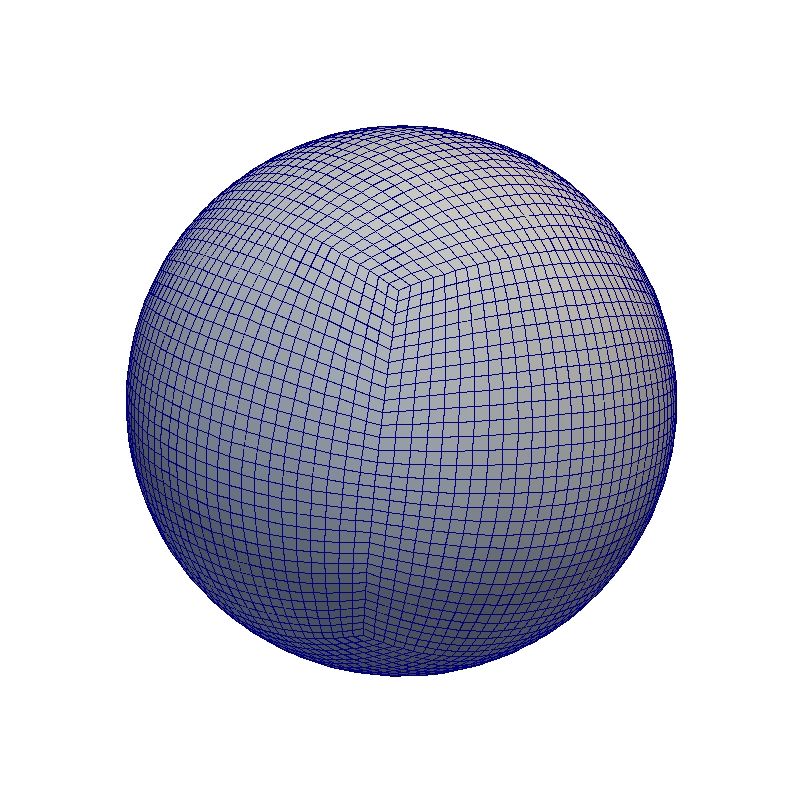}\label{fig:ssmesh}}\;
\subfloat[Surface of unit sphere cut to show inside]{\includegraphics[width=28mm]{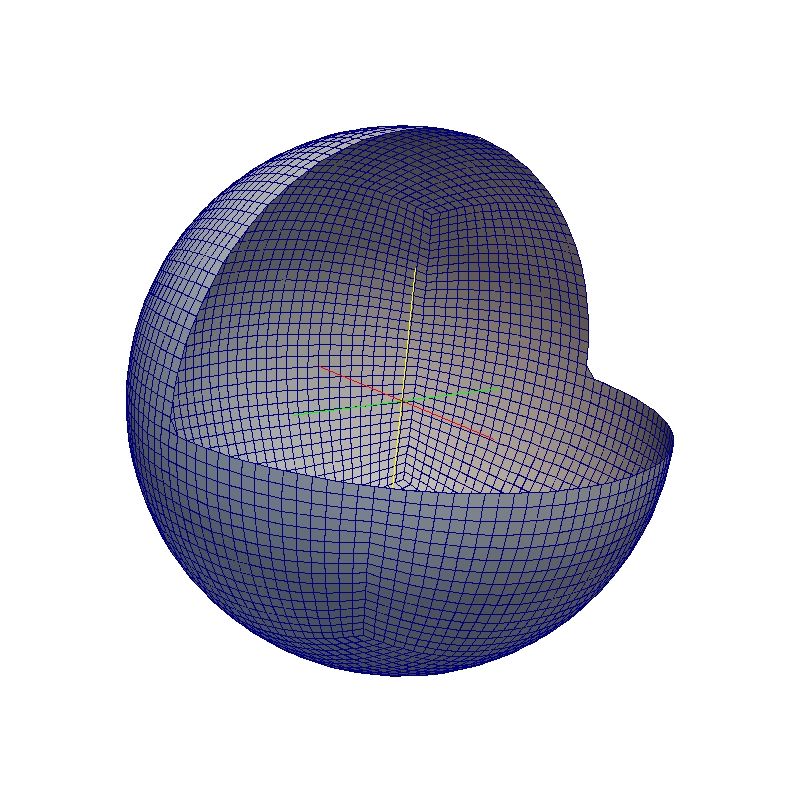}}\\
\subfloat[Ellipse]{\includegraphics[width=50mm]{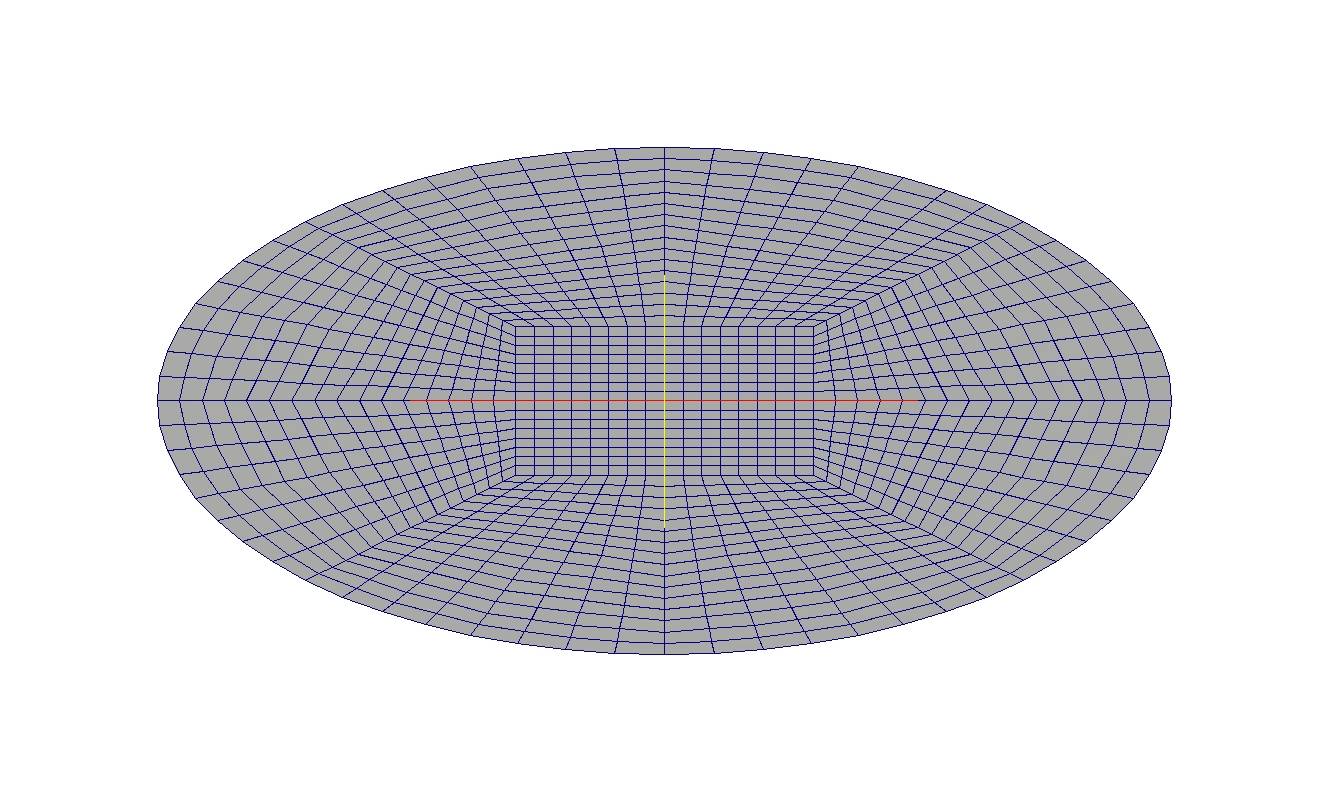}}\;
\subfloat[Dumbell mesh]{\includegraphics[width=30mm]{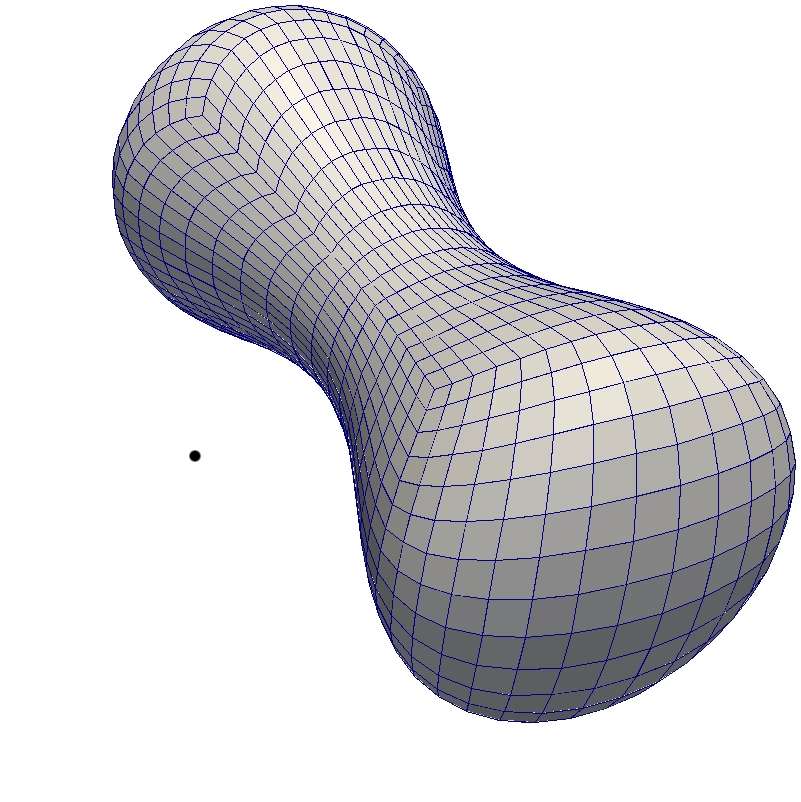}\label{fig:dumbbellgrid}}\;
\subfloat[Inner structure of dumbell mesh]{\includegraphics[width=30mm]{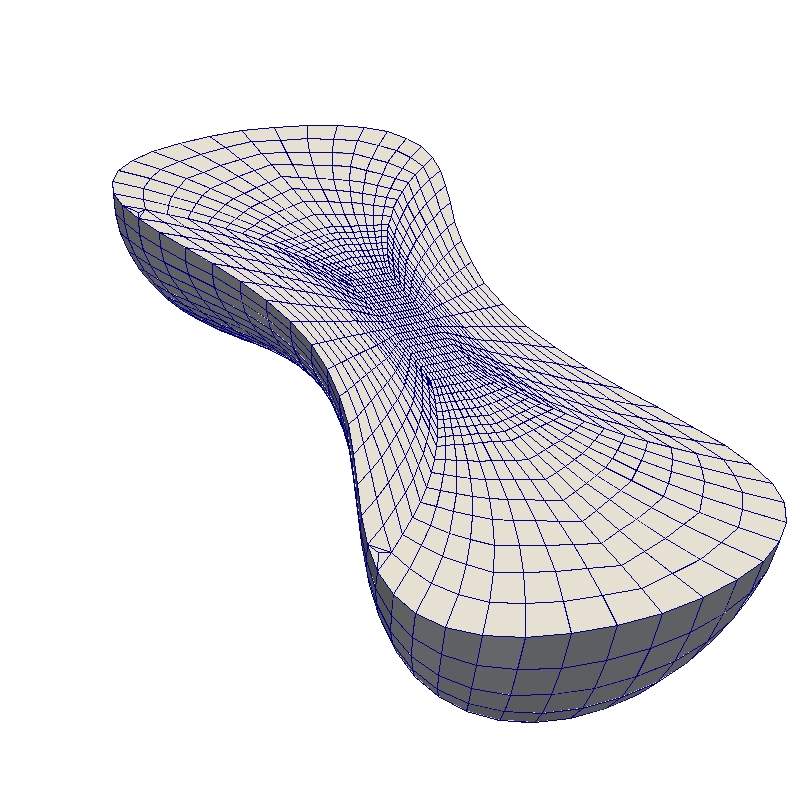}}\\
\subfloat["fish" mesh]{\includegraphics[width=50mm]{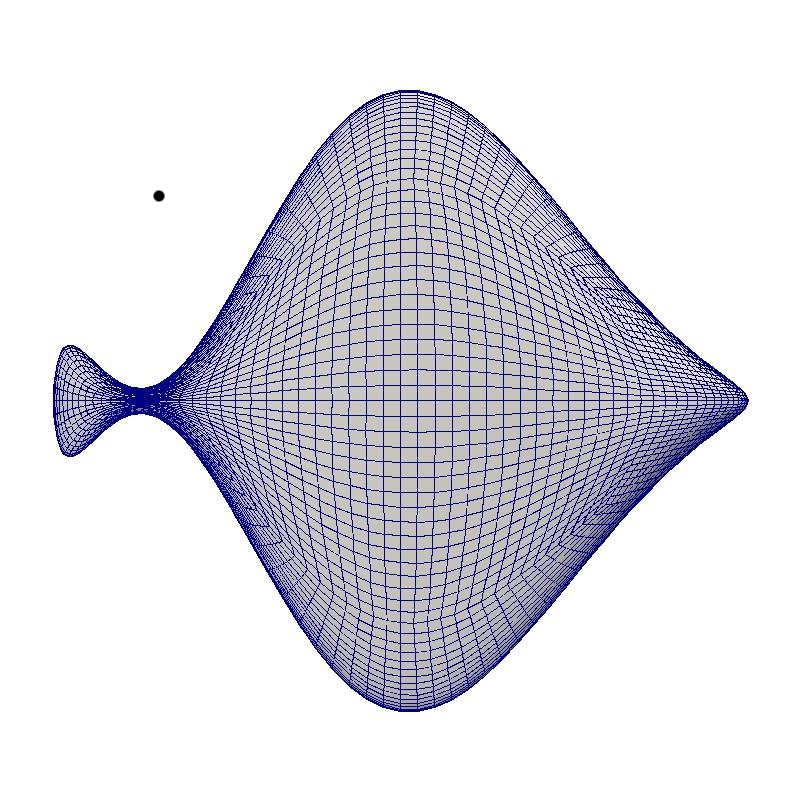}\label{fig:smmesh}}\;
\subfloat["eel" meshes (with and without boundary)]{\includegraphics[trim = 20cm 4cm 20cm 10cm, clip,width=50mm]{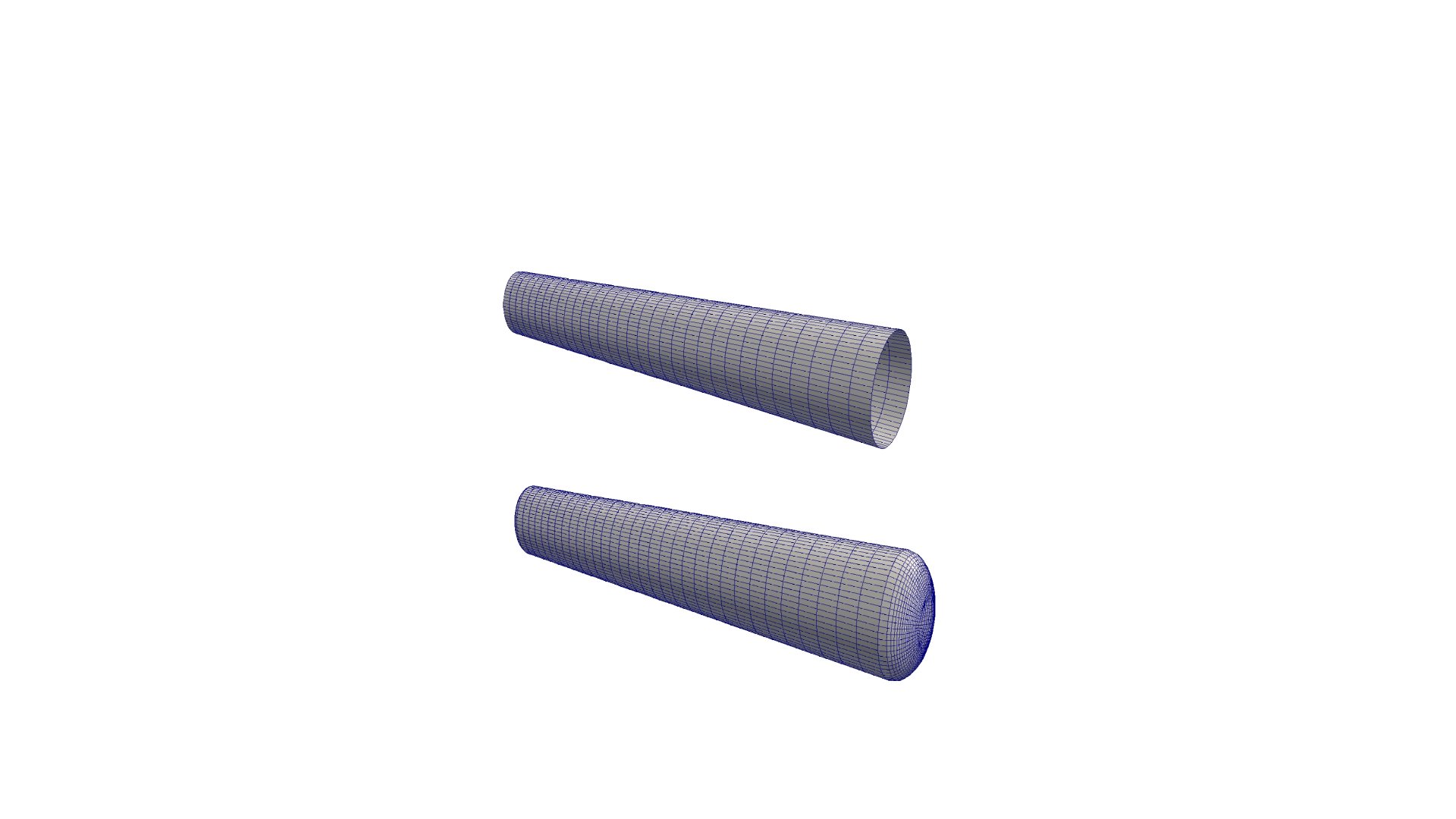}}
\caption{Examples of mesh generation for different volumes and surfaces: (a-c) Mesh generation on the unit sphere. (d) The ellipse which is a deformation of a circle mesh. (e-f) The dumbbell is a deformation of the bulk of a sphere. (g) The "fish" shape is a deformation of the surface of a sphere. (h) An "eel" is modelled by a cylinder with an open boundary and additionally as the same cylinder with added rounded ends}\label{fig:meshes}
\end{figure}
\section{Mesh generation}
All the mesh generation is carried out using the {\bf deal.II} library. We use hexahedral meshes for the volumes and quadrilaterals for the ellipse and surfaces. In Figure \ref{fig:meshes} we exhibit different meshes generated by this package on which we will carry out computations. We also consider smooth surfaces; these meshes are generated by creating a triangulation $\Omega_h$ of the bulk of the domain $\Omega$ then the surface triangulation is defined by collecting the faces of the elements of the bulk triangulation that lie on the surface ($\Gamma_h=\Omega_h|_d\Omega$), i.e., the surface mesh is the trace of the volume mesh (in the example of the cylinder with open ends we use only the elements on the curved surface). For this reason the equations are not being approximated on the actual surface but on an approximation of it. For more details on surface mesh generation the reader is referred to \cite{dealii} and the references therein.

\section{Comparisons of eigenfunctions and spatially inhomogeneous steady states}
\subsection{Example 1: Sphere}
We start by considering the unit sphere, a domain for which the eigenvalue problem can be solved analytically.
\subsubsection{Eigenvalues and eigenfunctions of the Laplacian in the bulk of the unit sphere}

In order to solve \eqref{eq:equn} on the sphere, we convert the eigenvalue problem into spherical coordinates. The eigenvalue problem in spherical coordinates is as follows \citep{arfken,morimoto},
\begin{equation*}
\Delta w + k^2 w =\frac{1}{r^2}\pderi{r}\left(r^2\pder{w}{r}\right)+\frac{1}{r^2\sin\theta} \pderi{\theta}\left(\sin\theta\pder{w}{\theta}\right)+\frac{1}{r^2\sin\theta}\spdr{w}{\phi} + k^2 w=0,
\end{equation*}
with homogeneous Neumann boundary conditions.
The solutions of the above eigenvalue problem are well known and are obtained using separation of variables \citep{arfken,morimoto}. Following \cite{arfken} (p. 424-428) we find an infinite number of solutions of the form
\begin{subequations}\label{eq:efsph}
\begin{align*}
& \quad\quad\quad\quad w_{l,n}^m(r,\theta,\phi)=A_{l,n}^mJ_{l+\frac{1}{2}}(j'_{l+\frac{1}{2},n}r)e^{im\phi}P_l^m(\cos\theta),\\
\text{where} & 
 \begin{cases}
  & l,m,n \text{ all integers such that } |m|\leq l \leq n, \\
  & A_{l,n}^m \text{ are constants},\\
  & J_\alpha(x)=\sum_{j=0}^\infty\frac{(-1)^j}{j!\Gamma(1+j+\alpha)}\left(\frac{x}{2}\right)^{2j+\alpha} \\ & \text{ with } \Gamma(n)=(n-1)! \text{  (i.e. a Bessel function of the first kind)},\\
  & P_l^m(x) \text{ are associated Legendre polynomials},\\
  & j'_{l+\frac{1}{2},n} \text{ are zeros of the differential of the spherical Bessel function.}
 \end{cases}
\end{align*}
\end{subequations}
\begin{figure}
\centering 
\subfloat[$w_{1,1}^1$]{\includegraphics[trim = 50mm 50mm 0mm 30mm, clip, width=55mm]{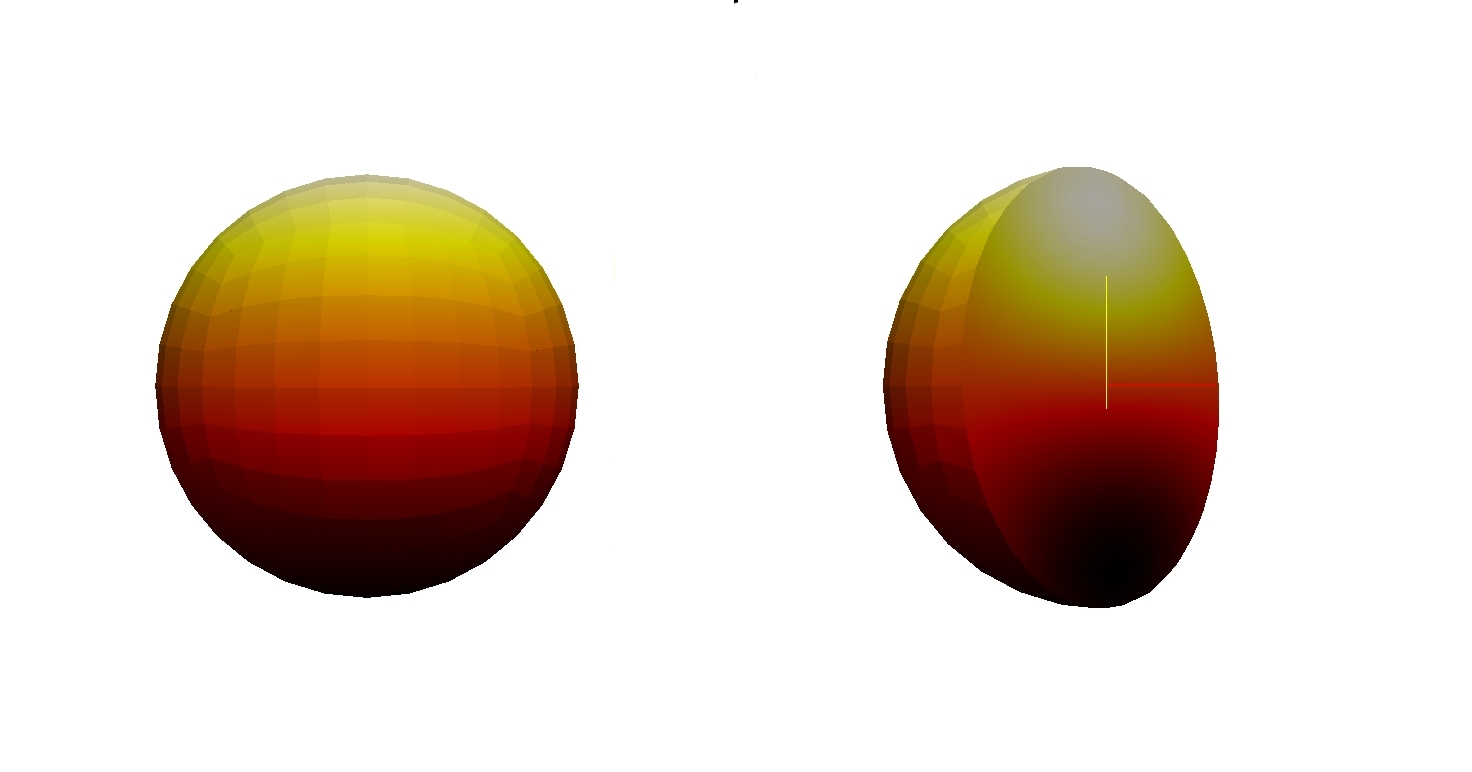}\label{fig:sphef1}}\;
\subfloat[$w_{2,1}^0$]{\includegraphics[trim = 0mm 50mm 50mm 30mm, clip, width=55mm]{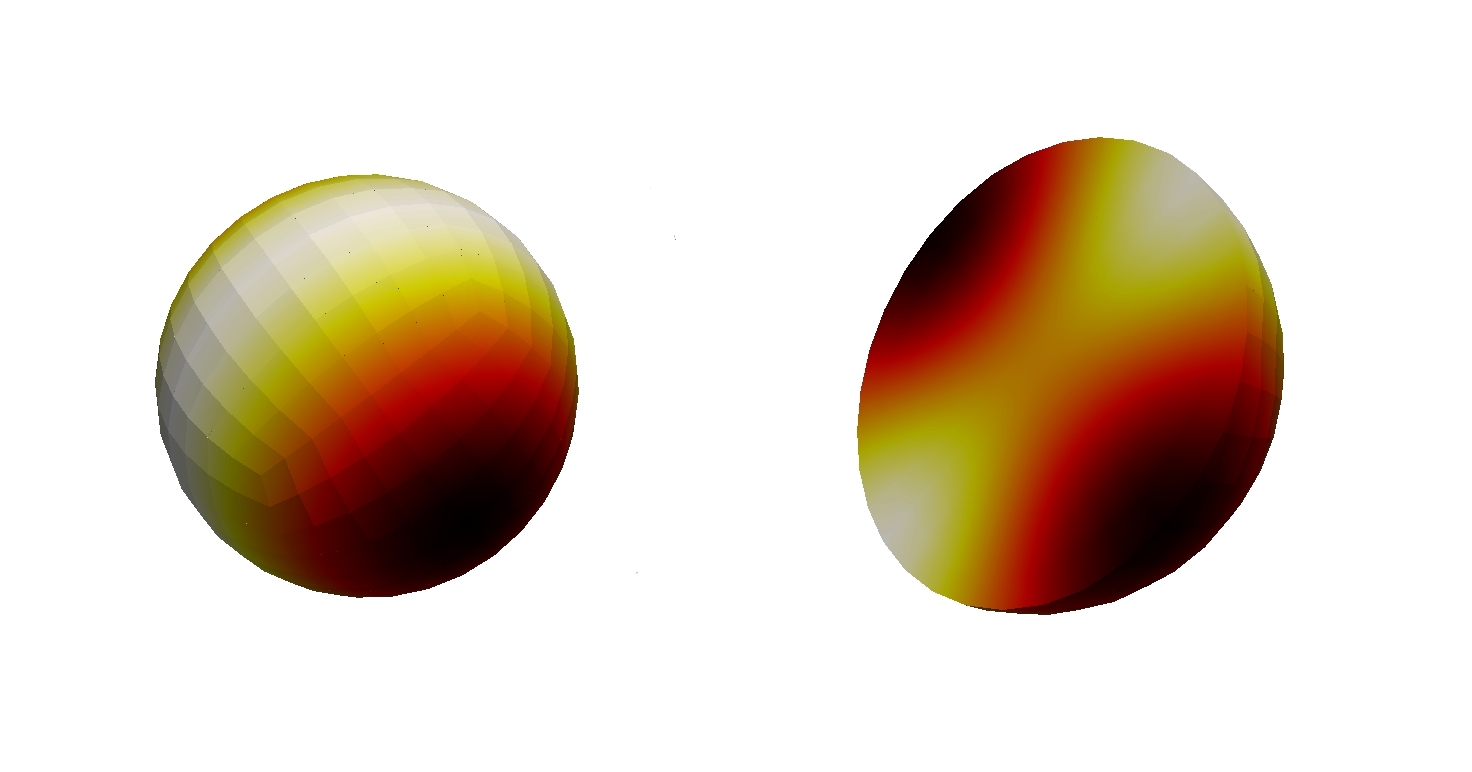}\label{fig:sphef2}}\\
\subfloat[$w_{3,1}^0$]{\includegraphics[trim = 50mm 50mm 0mm 30mm, clip, width=55mm]{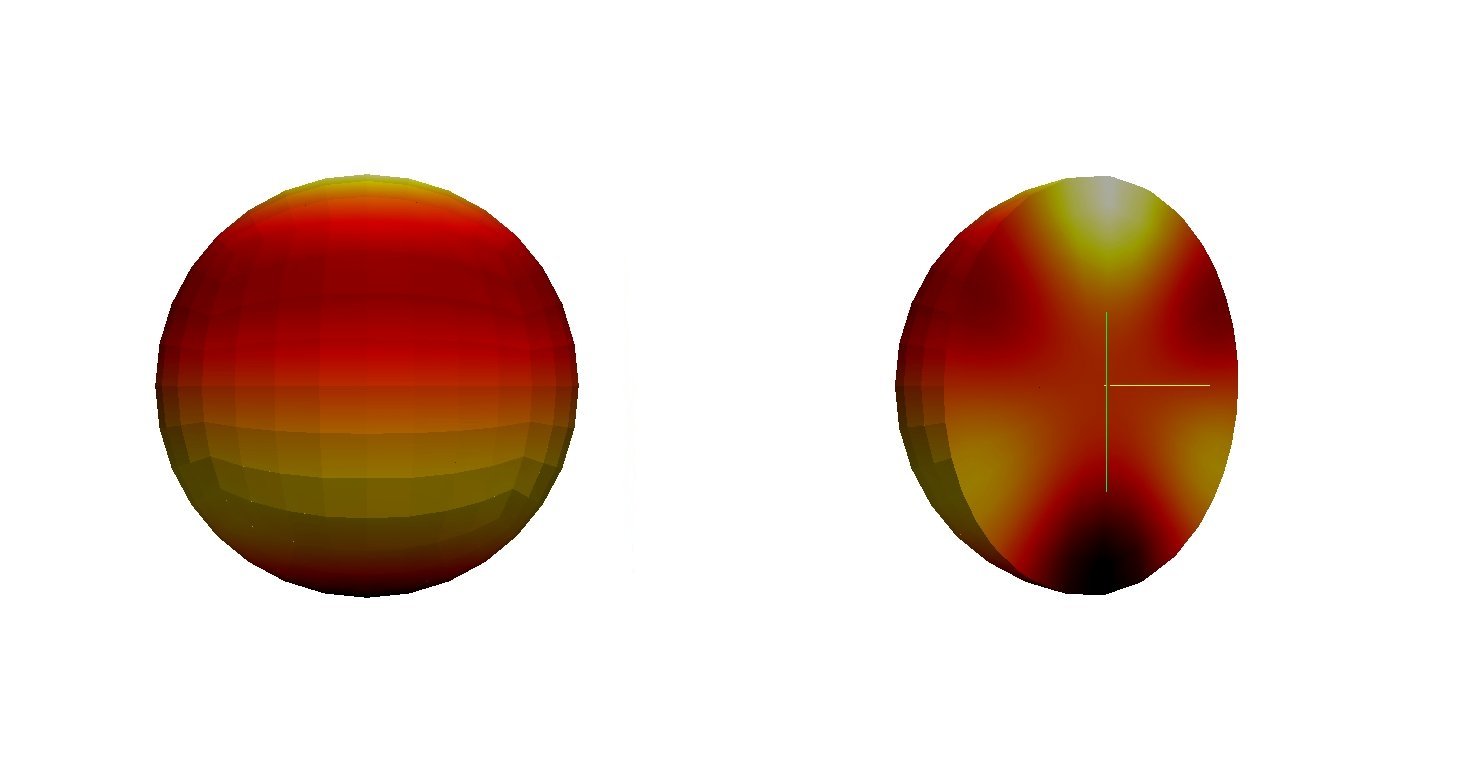}}\;
\subfloat[$w_{3,1}^{-2}$]{\includegraphics[trim = 0mm 45mm 50mm 30mm, clip, width=55mm]{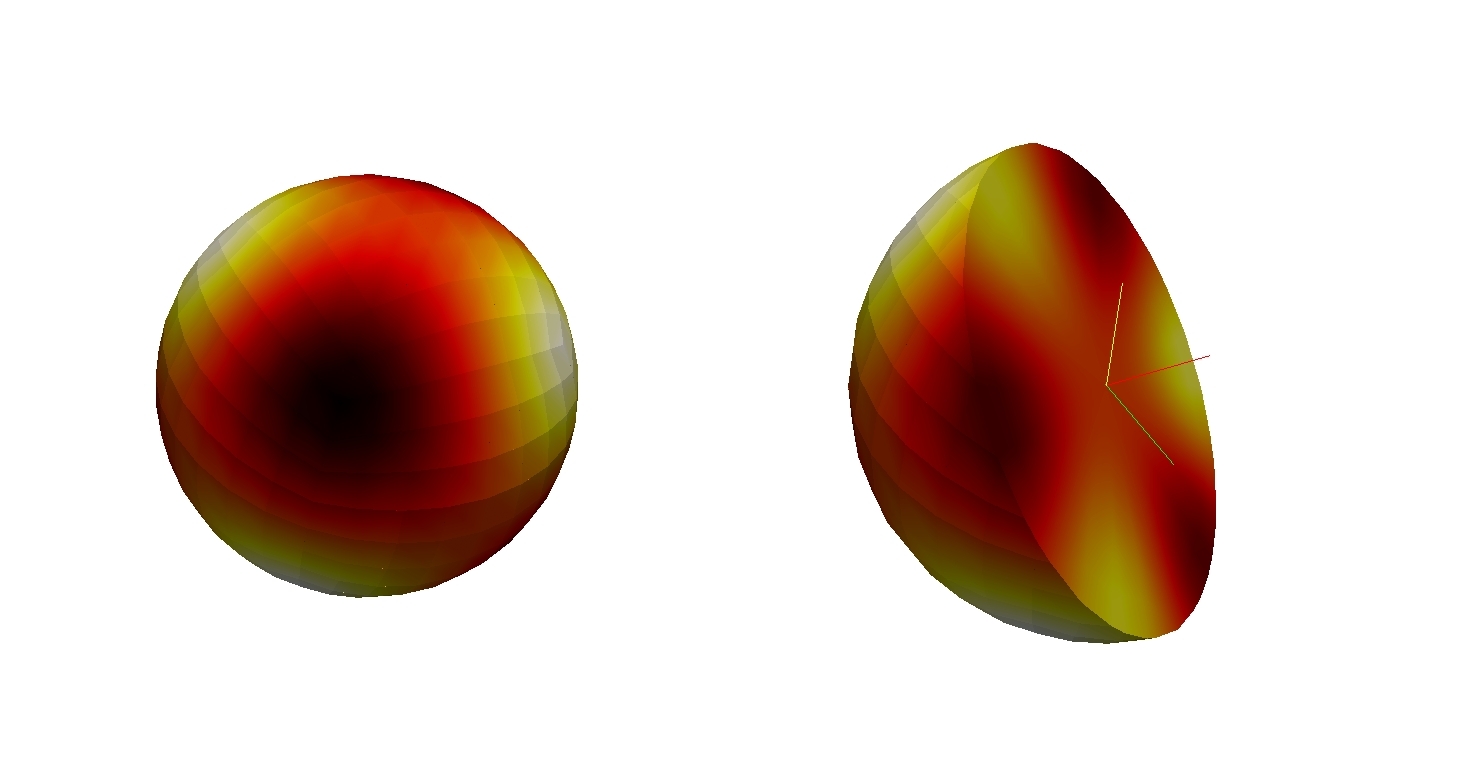}}\\
\subfloat[$w_{4,1}^{-3}$]{\includegraphics[trim = 50mm 45mm 0mm 30mm, clip, width=55mm]{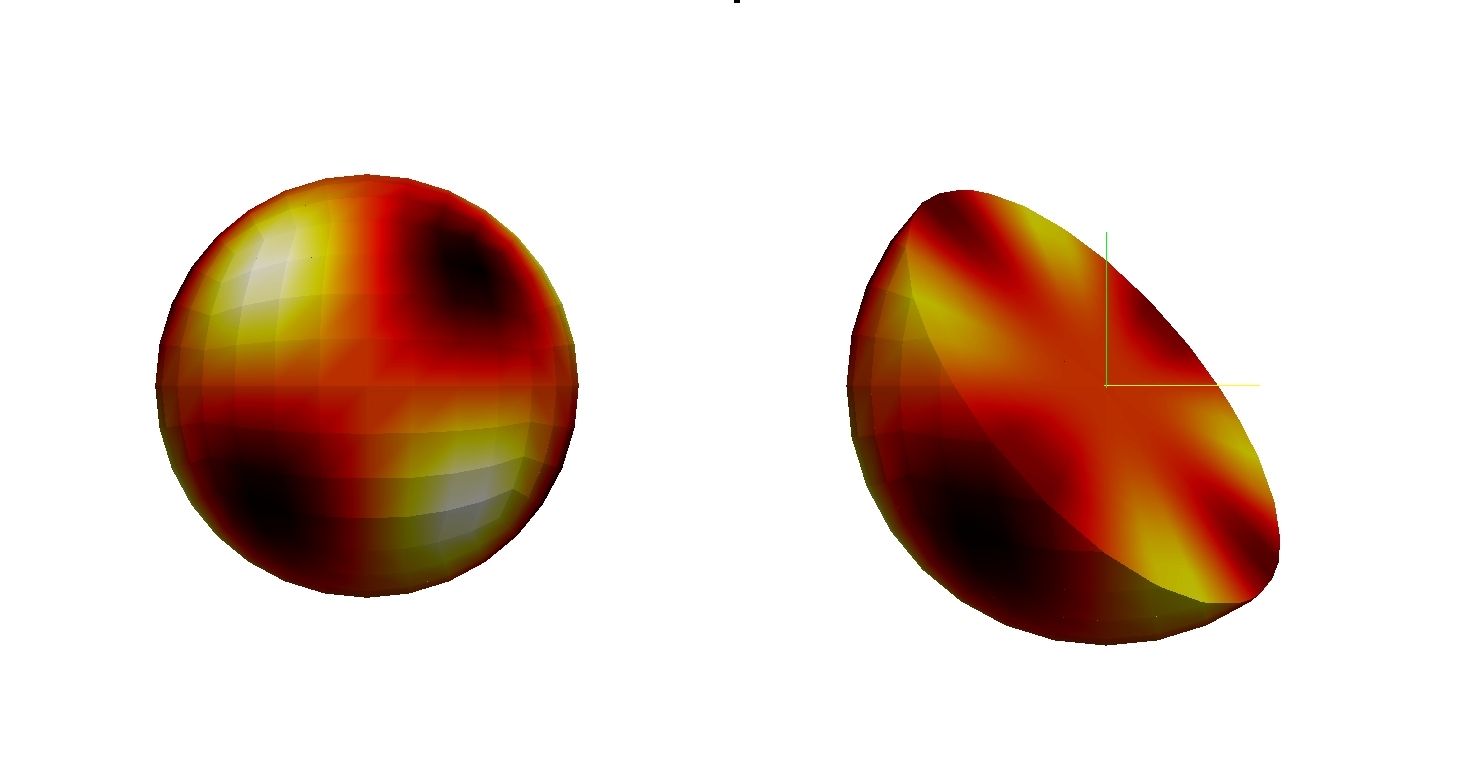}}
\caption{Analytical solutions to the eigenvalue problem on the unit sphere i.e. \eqref{eq:efsph} for selected values of $l$, $m$ and $n$. These are plotted using {\bf deal.II} (Colour version online)}\label{fig:sphereefs}
\end{figure}
We can find the eigenvalues $k_{l,n}^2=(j'_{l+\frac{1}{2},n})^2$ numerically (using the fact that $J'_{l+\frac{1}{2},n}=\frac{l}{k}J_{l+\frac{1}{2}}(k)-J_{l+\frac{3}{2}}(k)$). It follows that for each eigenvalue $\lambda_{l,n}=k_{l,n}^2$ there are $2l+1$ possible eigenfunctions. Figure \ref{fig:sphereefs} shows the eigenfunctions for some selected values of $l$, $m$ and $n$. For example $k_{1,1}=2.08158$ is the first zero of $J_{\frac{3}{2}}(x)$ and corresponds to the eigenfunctions
\[
w_{1,1}^m(r,\theta,\phi)=J_{\frac{3}{2}}(k_{1,1}r)e^{im\phi}P_1^m(\cos\theta), \text{ with } m=-1,0,1.
\]
The spherical Bessel function is given by $J_{\frac{3}{2}}(k_{1,1}r)=\frac{\sin(k_{1,1}r)}{(k_{1,1}r)^2}-\frac{\sin(k_{1,1}r)}{(k_{1,1}r)^2}$.
Meanwhile $Y_1^m=e^{im\phi}P_1^m(\cos\theta)$ are spherical harmonics whose real parts can be written in cartesian coordinates as $Y_1^{-1}=\sqrt{\frac{3}{4\pi}}\cdot\frac{y}{r}$, $Y_1^{0}=\sqrt{\frac{3}{4\pi}}\cdot\frac{z}{r}$ and $Y_1^{1}=\sqrt{\frac{3}{4\pi}}\cdot\frac{x}{r}$. Since the system we are solving is not sensitive to polarity we can consider these to be equivalent. Figure \ref{fig:sphereefs}\subref{fig:sphef1} shows a plot of the eigenfunction 
\[
 w_{1,1}^1=(\frac{\sin(k_{1,1}r)}{(k_{1,1}r)^2}-\frac{\sin(k_{1,1}r)}{(k_{1,1}r)^2})\cdot\frac{x}{r},
\]
where as usual $r^2=x^2+y^2+z^2$.
The second example, $k_{2,1}=3.34209$ corresponds to the eigenfunctions
\[
w_{2,1}^m(r,\theta,\phi)=J_{\frac{5}{2}}(k_{2,1}r)e^{im\phi}P_2^m(\cos\theta) \text{, with } -l\leq m\leq l.
\]
Choosing $m=0$, converting the above to cartesian coordinates and taking the real part gives
\begin{equation*}
\begin{split}
w_{2,1}^0 & (x,y,z)= \\ & \left(\left(\frac{3}{k_{2,1}^2r^2}-1\right)\frac{\sin(k_{2,1}r)}{k_{2,1}r}-\frac{3\cos(k_{2,1}r)}{k_{2,1}^2r^2}\right)
\left(\frac{1}{4}\sqrt{\frac{5}{\pi}}\cdot\frac{-x^2-y^2+2z^2}{r^2}\right).
\end{split}
\end{equation*}
The function $w_{2,1}^0$ is plotted in Figure \ref{fig:sphereefs}\subref{fig:sphef2}.

\subsubsection{Mode isolation on the sphere}
Using the method described in Section \ref{sec:isolation} and the values given in Table \ref{table:parameters} we can isolate the wavenumbers for the reaction-diffusion system with Schnakenberg kinetics and these are shown in Table \ref{table:schnak}. We can do the same for Thomas and Gierer-Meinhart (Table \ref{table:gmtom}). In these cases the interval $[k_-,k_+]$ is centered on $k_{l,n}$.
\begin{table}[H]
\caption{Given $d$ and $\gamma$ from the first two columns we obtain values for $k_-$ and $k_+$ and this means that particular given wavenumbers are isolated on the sphere for the reaction-diffusion system with Schnakenberg kinetics.}\label{table:schnak}
\begin{center}
 \begin{tabular}{| c | c | c | c | c |}
 \hline
$d$ & $\gamma$ & $k_-$ & $k_+$ & Wavenumbers excited   \\ \hline \hline
10 & 15 & 1.7321 & 2.7386 & $k_{1,1}=2.08158$ \\ \hline
10 & 40 & 2.8284 & 4.4721 & $k_{2,1}=3.34209$ \\ \hline
9 & 60 & 3.9319 & 5.0866 & $k_{0,2}=4.49341$, $k_{3,1}=4.51410$ \\ \hline
8.81 & 85 & 4.8575 & 5.8955 & $k_{4,1}=5.64670$ \\ \hline
\end{tabular}
\end{center}
\end{table}

\begin{table}[H]
\caption{The values of $d$ and $\gamma$ which isolate the given wavenumbers on the sphere for the Gierer-Meinhart and Thomas reaction kinetics.}
\label{table:gmtom}
\begin{center}
 \begin{tabular}{| r | r | c |}
 \hline
Gierer-Meinhart & Thomas & Wavenumbers excited \\ \hline \hline
d=74 $\gamma$=30 & d=30 $\gamma$=15 & $k_{1,1}$ \\ \hline
d=74 $\gamma$=80 & d=30 $\gamma$=40 & $k_{2,1}$ \\ \hline
d=74 $\gamma$=160 & d=28, $\gamma$=60 & $k_{0,2}$, $k_{3,1}$ \\ \hline
d=72 $\gamma$=200 & d=27.5 $\gamma$=90 & $k_{4,1}$ \\ \hline
\end{tabular}
\end{center}
\end{table}

\subsubsection{Simulations of the reaction-diffusion systems on the unit sphere}\label{sec:sphere}
\begin{figure}
\centering 
\subfloat[$\gamma=15$, d=10]{\includegraphics[trim = 10mm 30mm 10mm 30mm, clip,width=55mm]{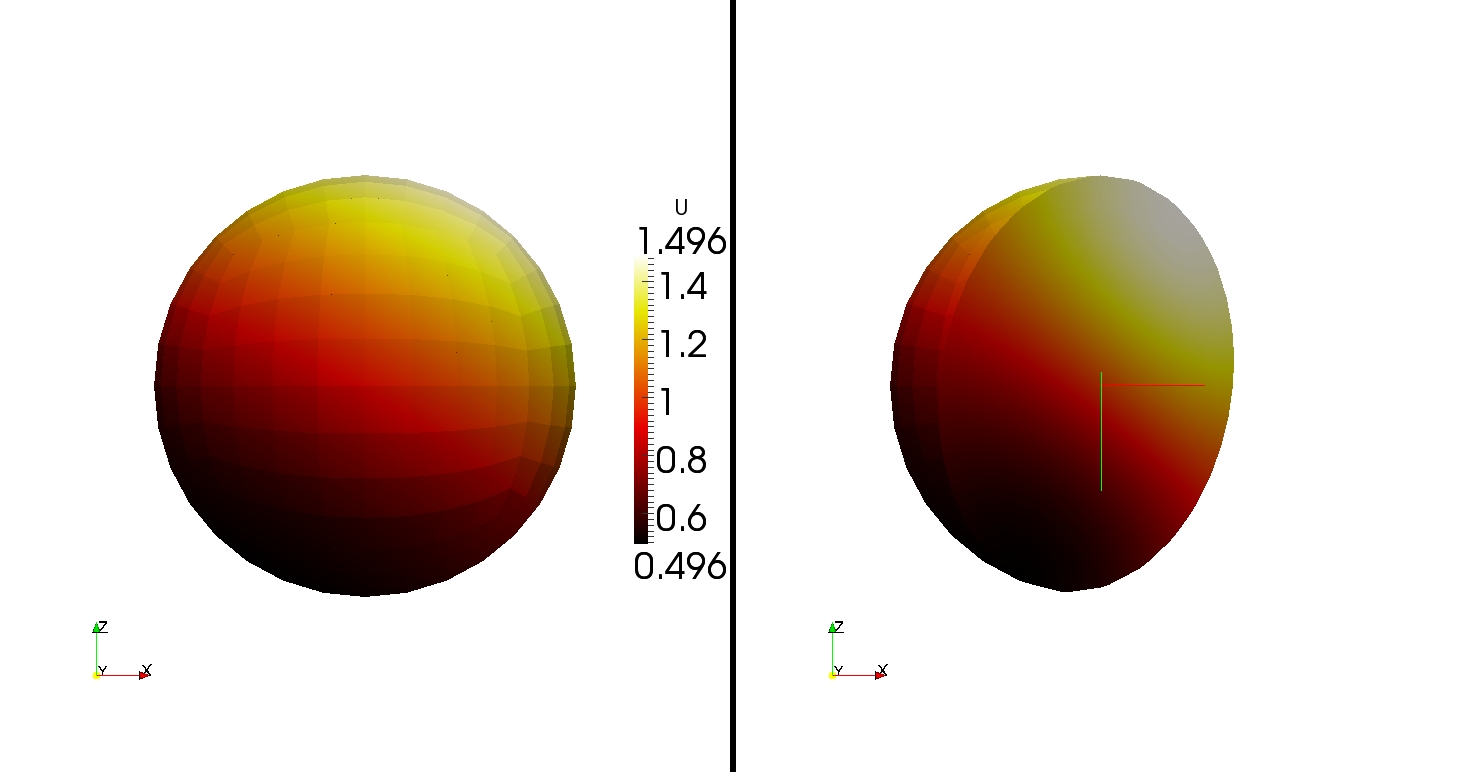}}\;\;
\subfloat[$\gamma=40$, d=10]{\includegraphics[trim = 10mm 30mm 10mm 30mm, clip,width=55mm]{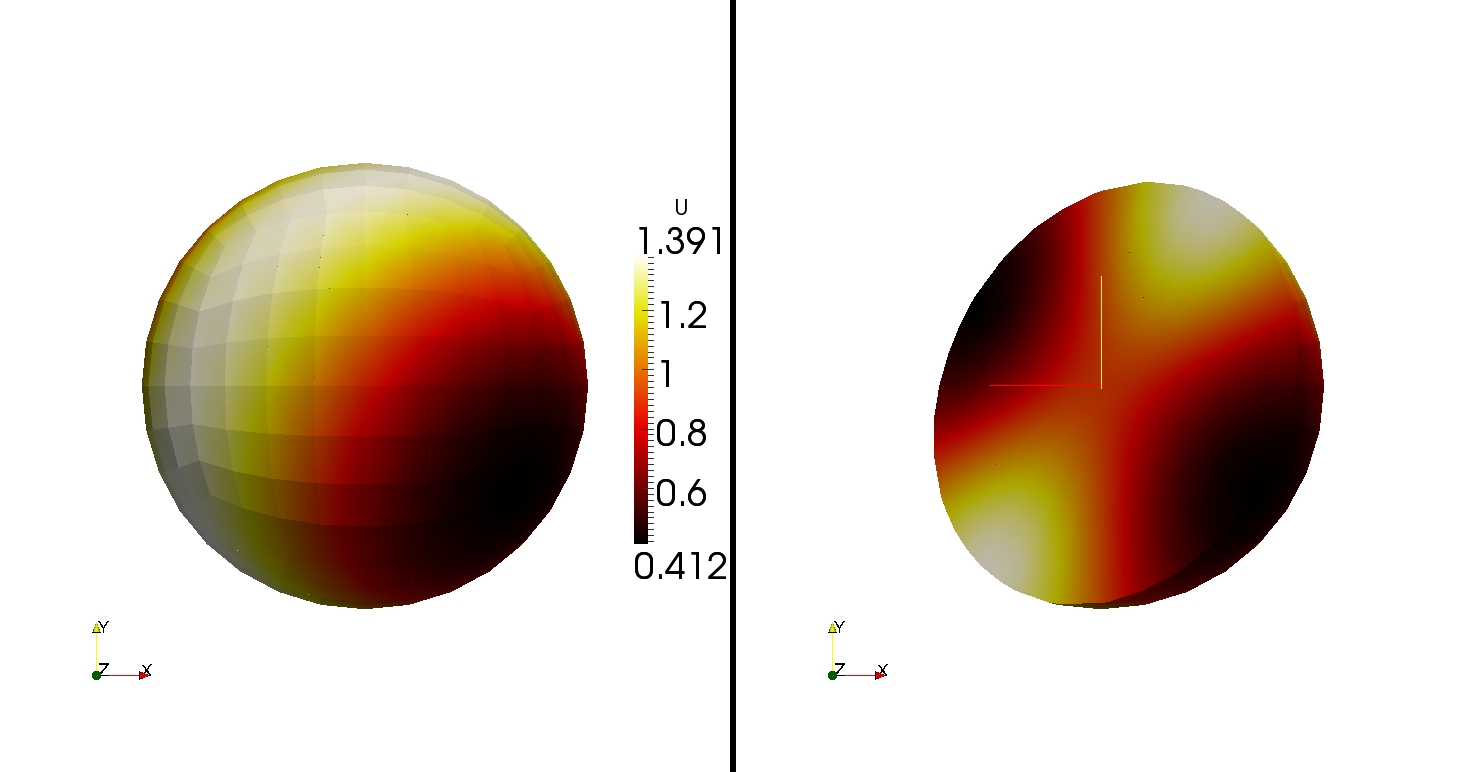}}\\
\subfloat[$\gamma=70$, d=9]{\includegraphics[trim = 10mm 30mm 10mm 30mm, clip,width=55mm]{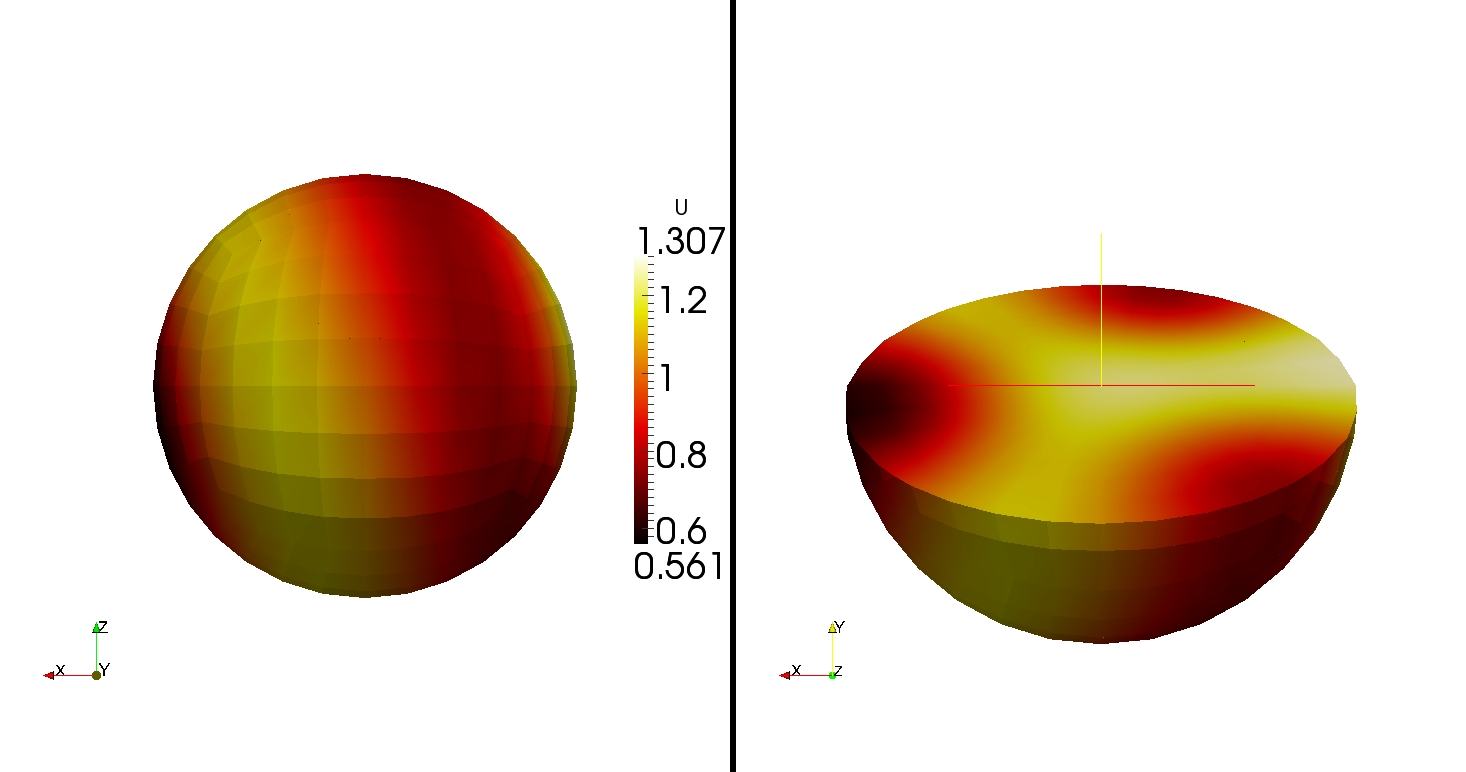}}\;\;
\subfloat[$\gamma=85$, d=8.81]{\includegraphics[trim = 10mm 30mm 10mm 30mm, clip,width=55mm]{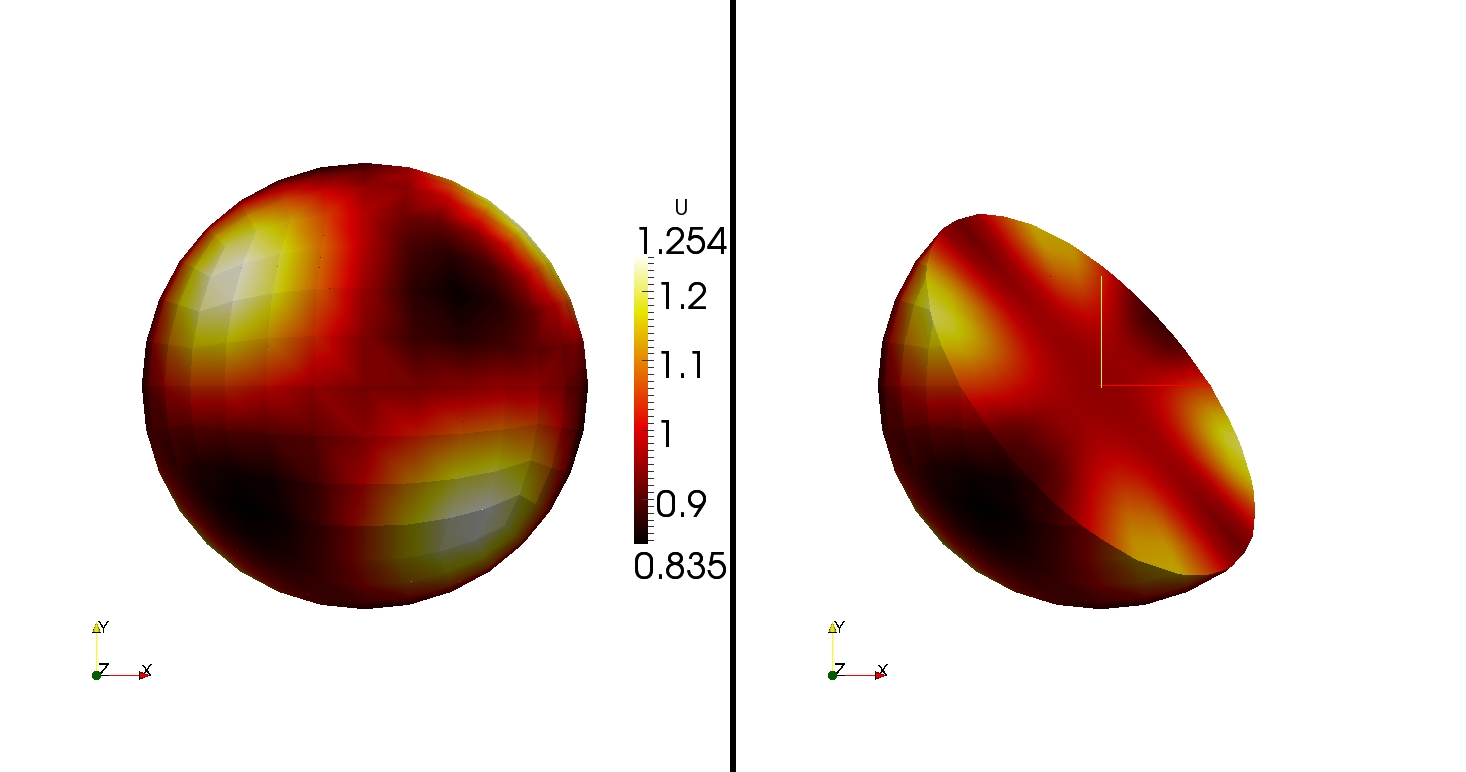}}
\caption{Converged solutions of system \eqref{eq:non} with Schnakenberg kinetics \eqref{eq:schnakkin}. These solutions represent the species $u$. The isolated modes are $w_{1,1}^1$, $w_{2,1}^0$, $w_{3,1}^0$ and $w_{4,1}^{-3}$ (Colour version online)}\label{fig:schnakspheres}
%
\subfloat[GM, $\gamma=80$, d=74]{\includegraphics[trim = 10mm 30mm 10mm 30mm, clip,width=55mm]{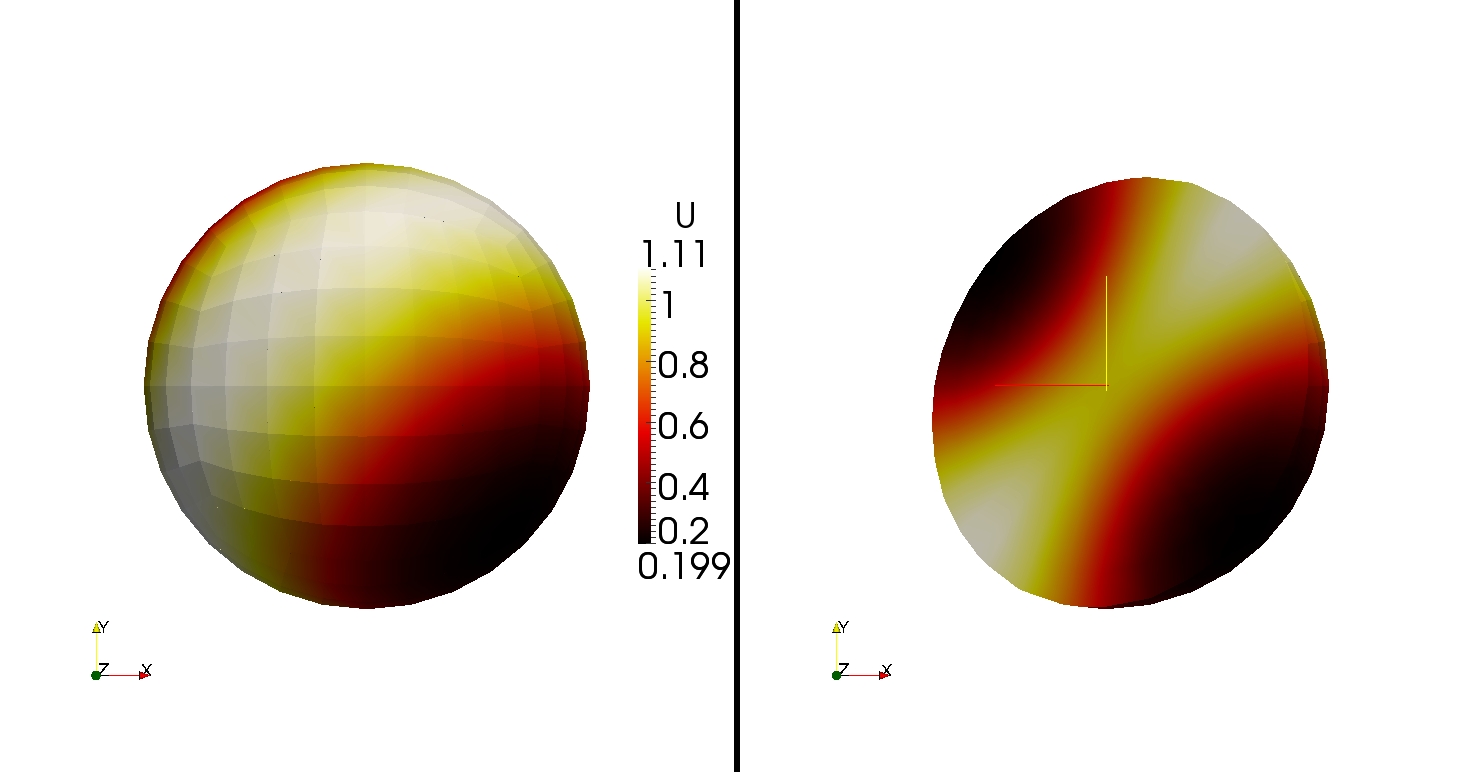}}\;\;
\subfloat[Thomas, $\gamma=40$, d=30]{\includegraphics[trim = 10mm 30mm 10mm 30mm, clip,width=55mm]{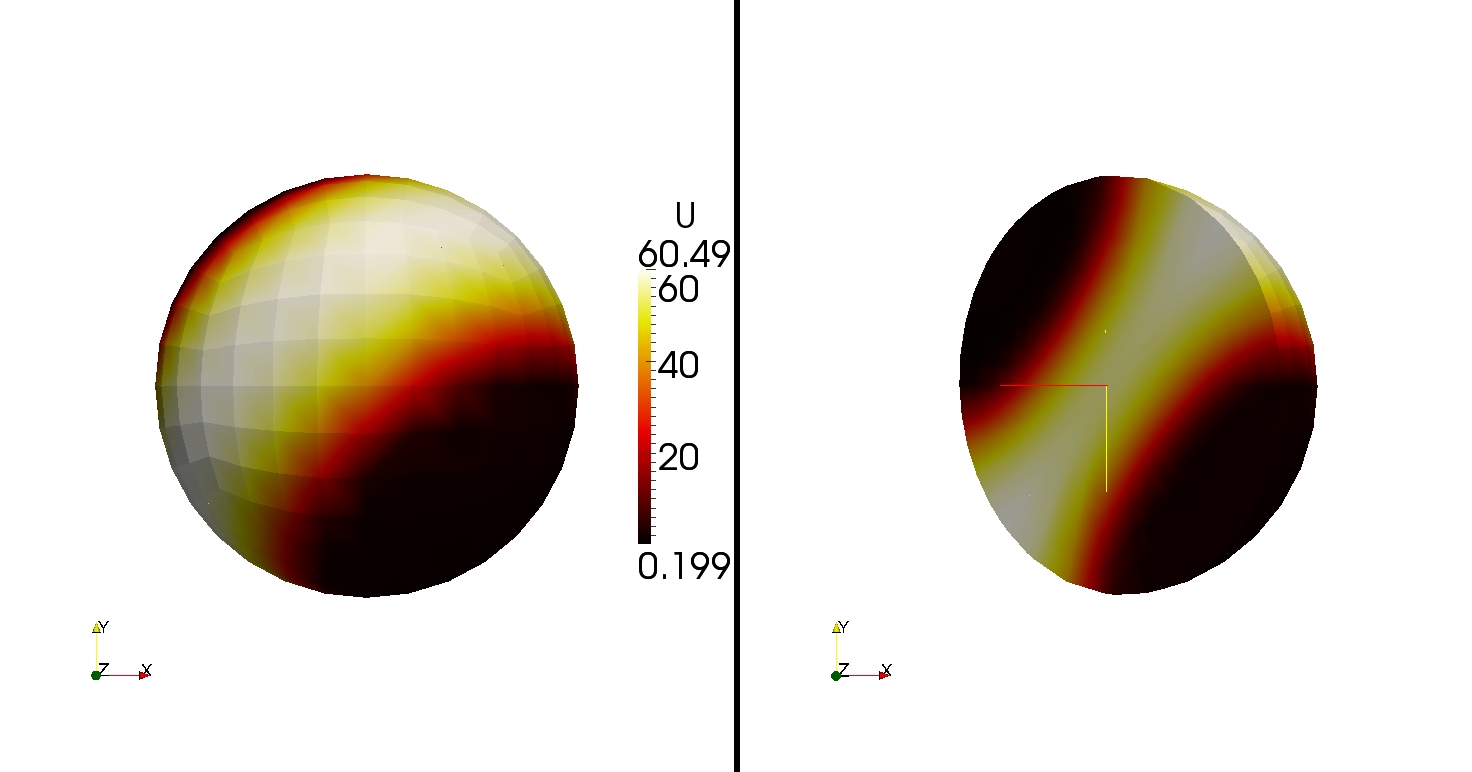}}\\
\subfloat[GM, $\gamma=160$, d=74]{\includegraphics[trim = 10mm 30mm 10mm 30mm, clip,width=55mm]{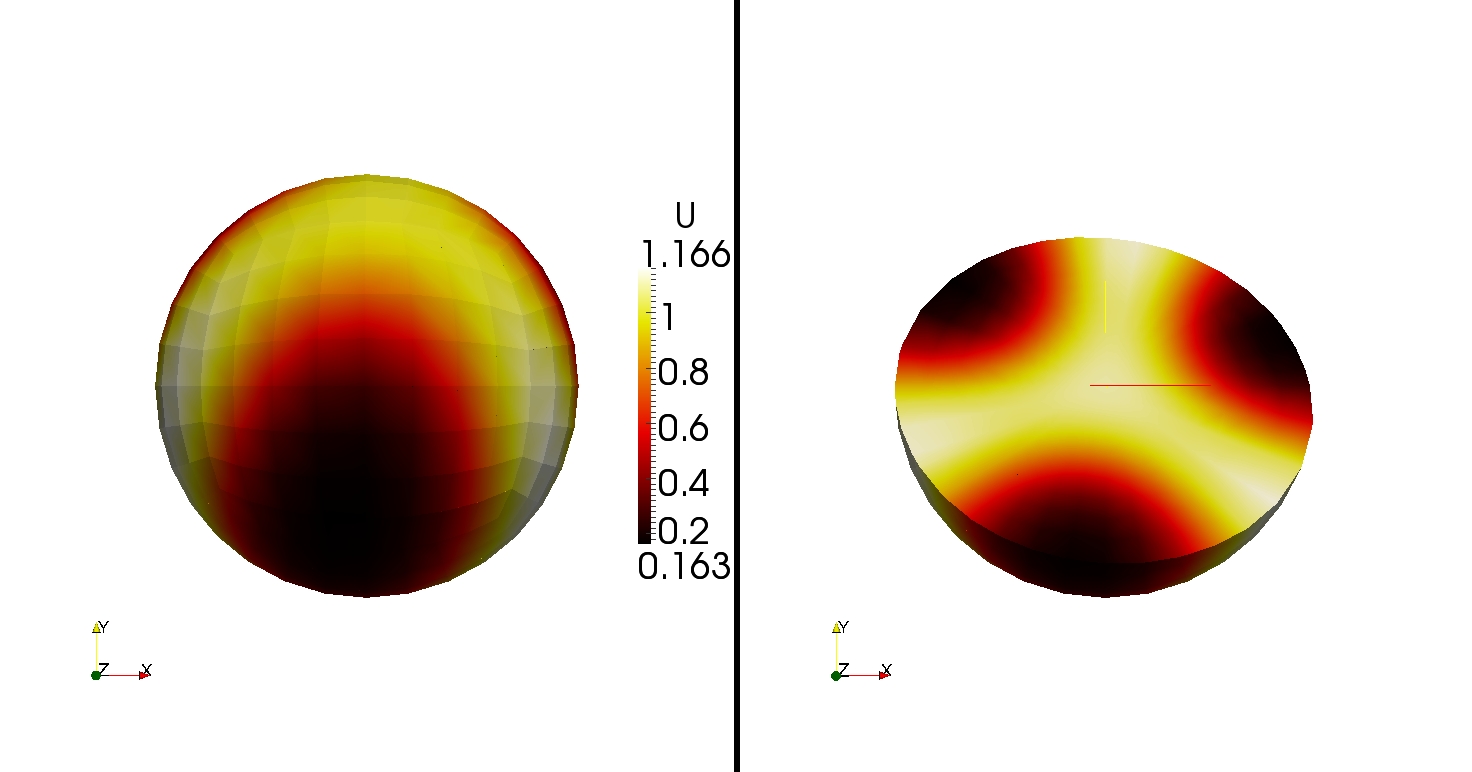}}\;\;
\subfloat[Thomas, $\gamma=70$, d=28]{\includegraphics[trim = 10mm 30mm 10mm 30mm, clip,width=55mm]{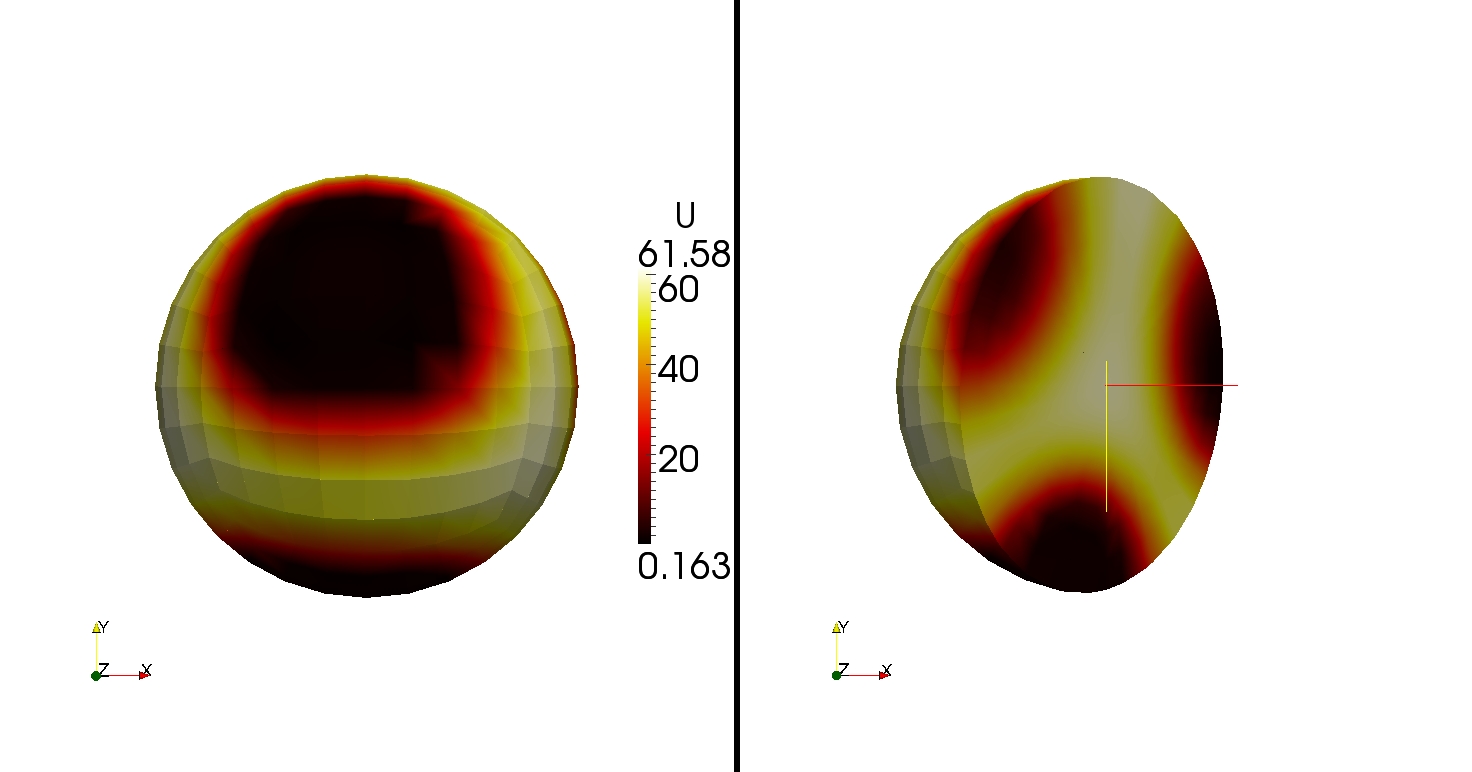}}\\
\subfloat[GM, $\gamma=200$, d=72]{\includegraphics[trim = 10mm 30mm 10mm 30mm, clip,width=55mm]{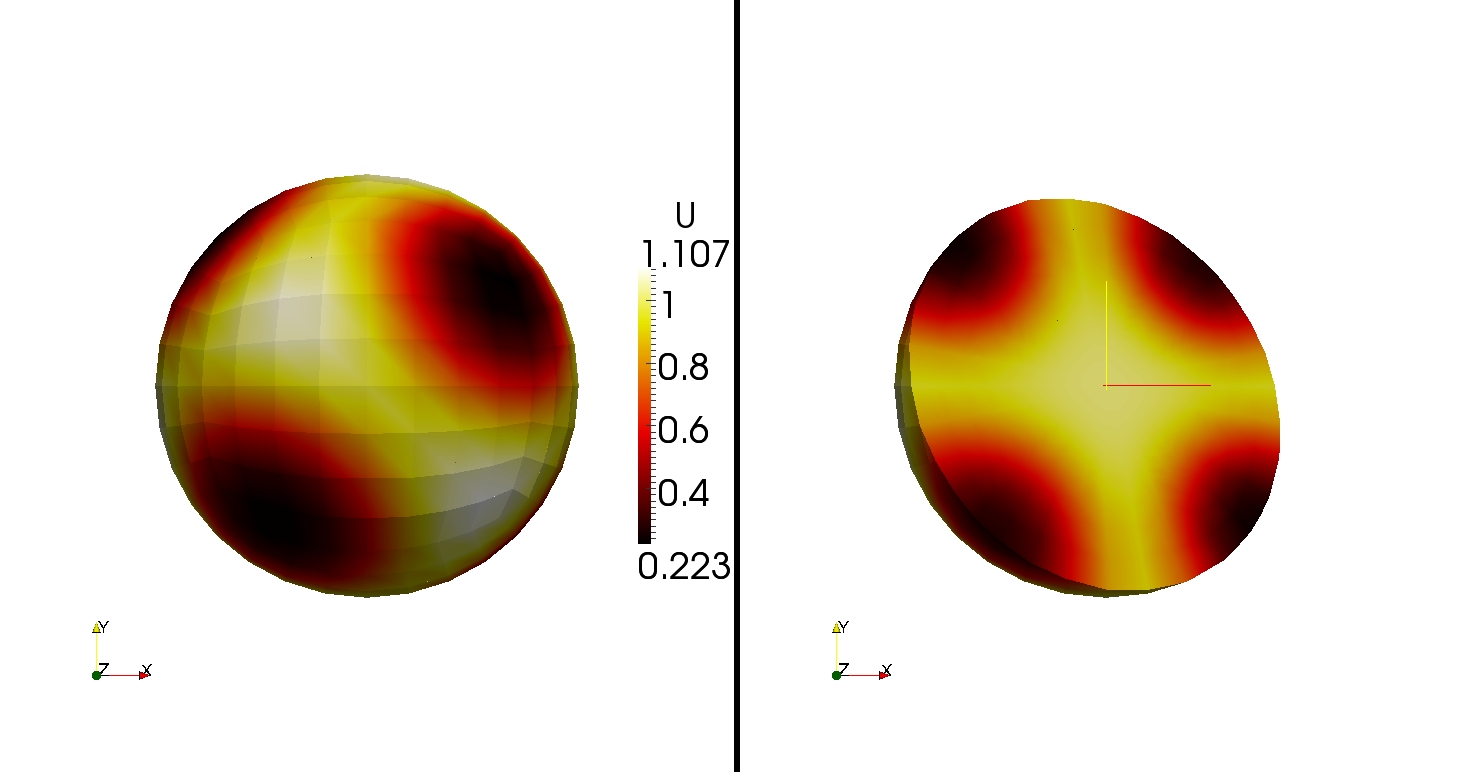}}\;\;
\subfloat[Thomas, $\gamma=90$, d=27.5]{\includegraphics[trim = 10mm 30mm 10mm 30mm, clip,width=55mm]{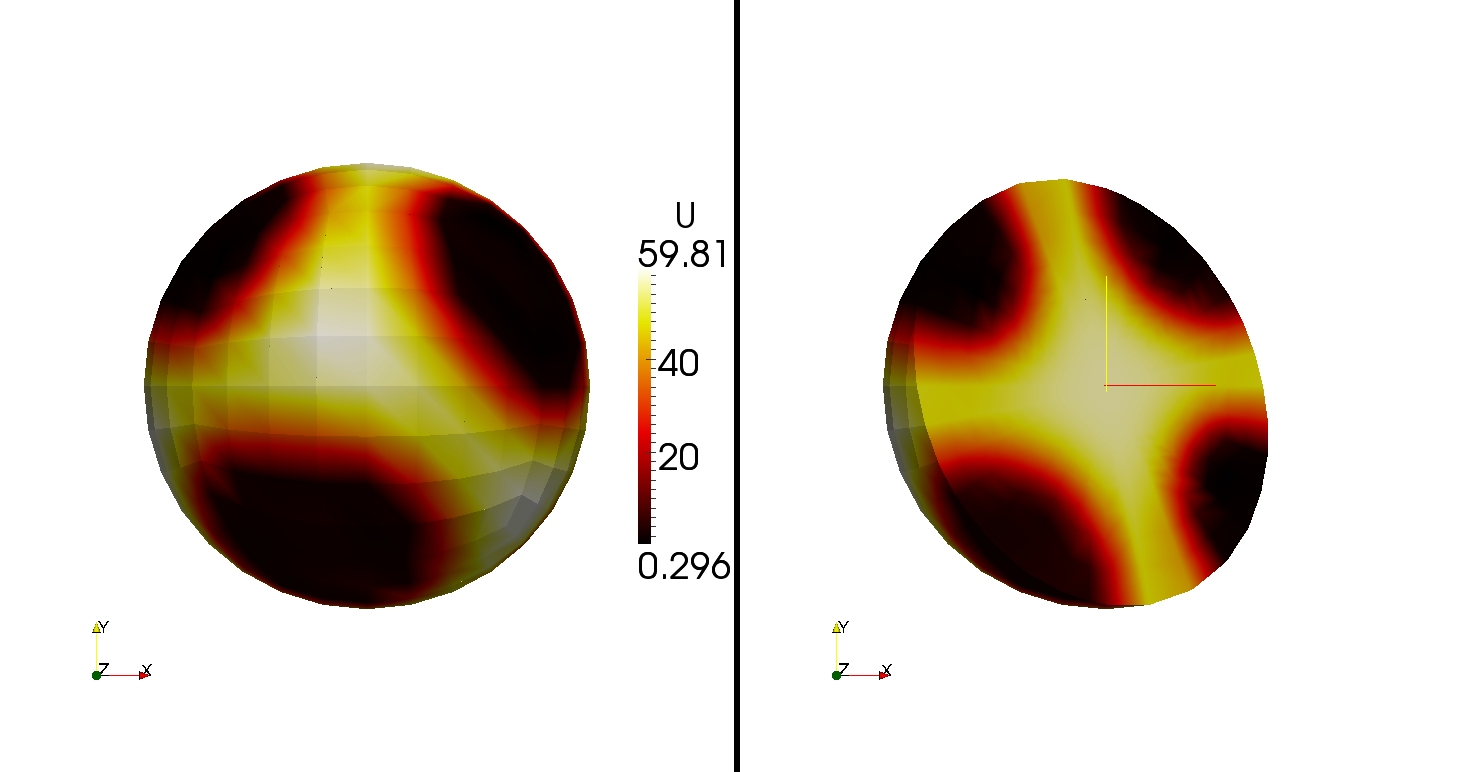}}
\caption{Converged solutions of system \eqref{eq:non} for the species $u$ with Gierer-Meinhart kinetics \eqref{eq:gm} on the left with isolated modes $w_{2,1}^0$, $w_{3,1}^3$ and $w_{4,1}^{-3}$ and Thomas \eqref{eq:tho} on the right with isolated modes $w_{2,1}^0$, $w_{3,1}^{-2}$ and $w_{4,1}^{-3}$ (Colour version online)}\label{fig:gmthomspheres}
\end{figure}
Solving using {\bf deal.II} we use the mesh shown in Figure \ref{fig:meshes}\subref{fig:sphmesh}. The timestep is taken to be $\tau=10^{-3}$. We take the initial conditions to be a small random perturbation from the previously computed homogeneous steady state. So for the reaction-diffusion system with Schnakenberg kinetics, at each point in the grid we set the initial conditions to be:
\begin{equation}
\alpha^0=0.995+0.01\epsilon, \quad \beta^0=0.895+0.01\epsilon,
\end{equation} 
where $\epsilon$ is a uniformly distributed random variable between $0$ and $1$. \\
For each eigenvalue there are a number of different eigenfunctions. Computing using the values obtained with mode isolation, the solution converges to either one of the eigenfunctions or a linear combination. These converged solutions are shown in Figures \ref{fig:schnakspheres} and \ref{fig:gmthomspheres}. It is possible to force the solution to converge to an eigenfunction (which it does not appear to with random initial perturbation) by making a suitable choice of initial condition, for example a perturbation of the desired eigenfunction, suitably scaled. Hence, in the case where multiple wave numbers are excited, pattern selection is heavily influenced by the choice of initial conditions which act as the basin of attraction, one of the major criticisms of Turing's theory for pattern formation \citep{BardandLauder1974}.

\begin{figure}
\subfloat[$\lambda_1=3.52$]{\includegraphics[trim = 30mm 0 30mm 0, clip, width=32mm]{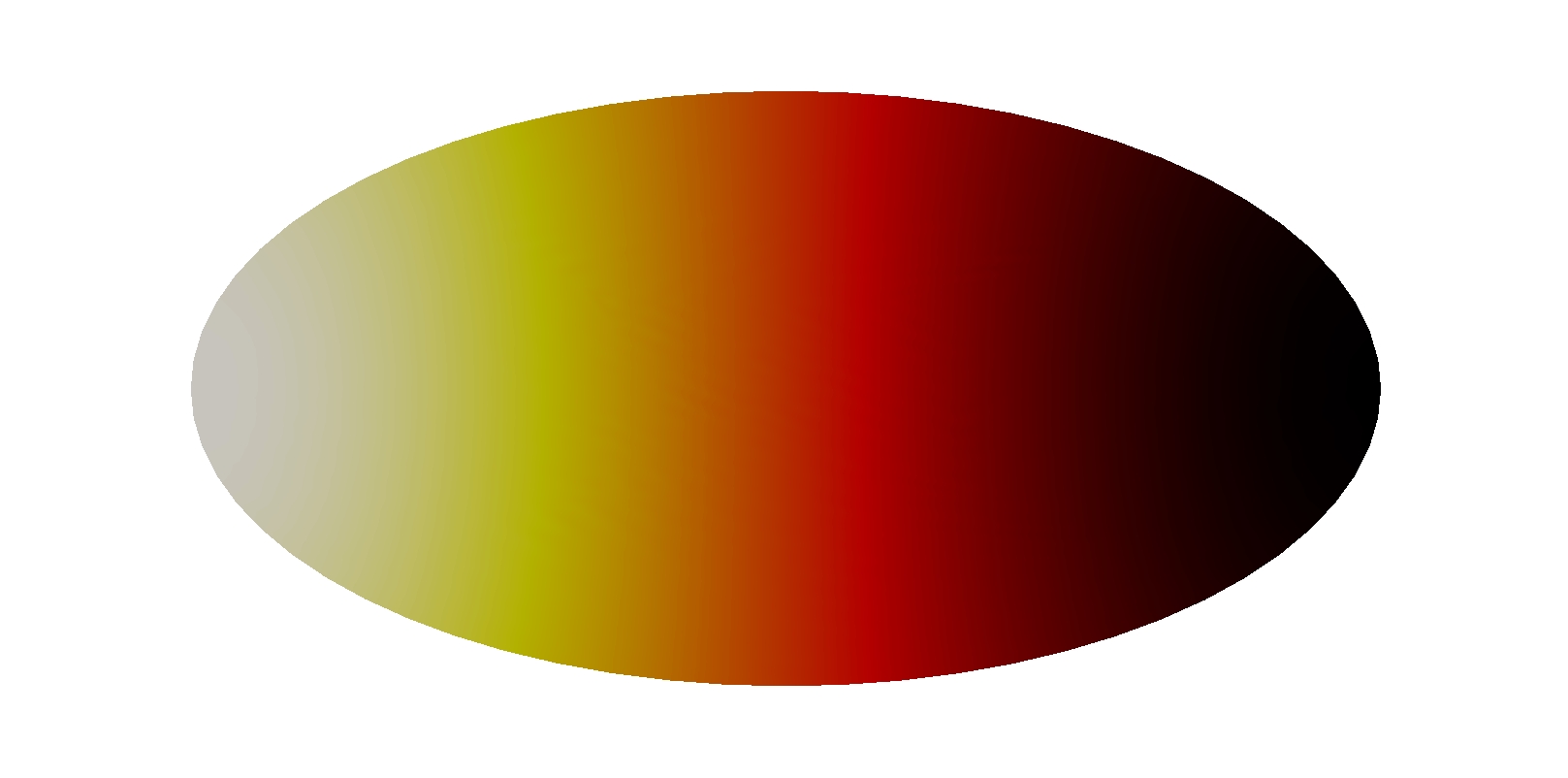}}\;
\subfloat[$\lambda_2=11.74$]{\includegraphics[trim = 30mm 0 30mm 0, clip, width=32mm]{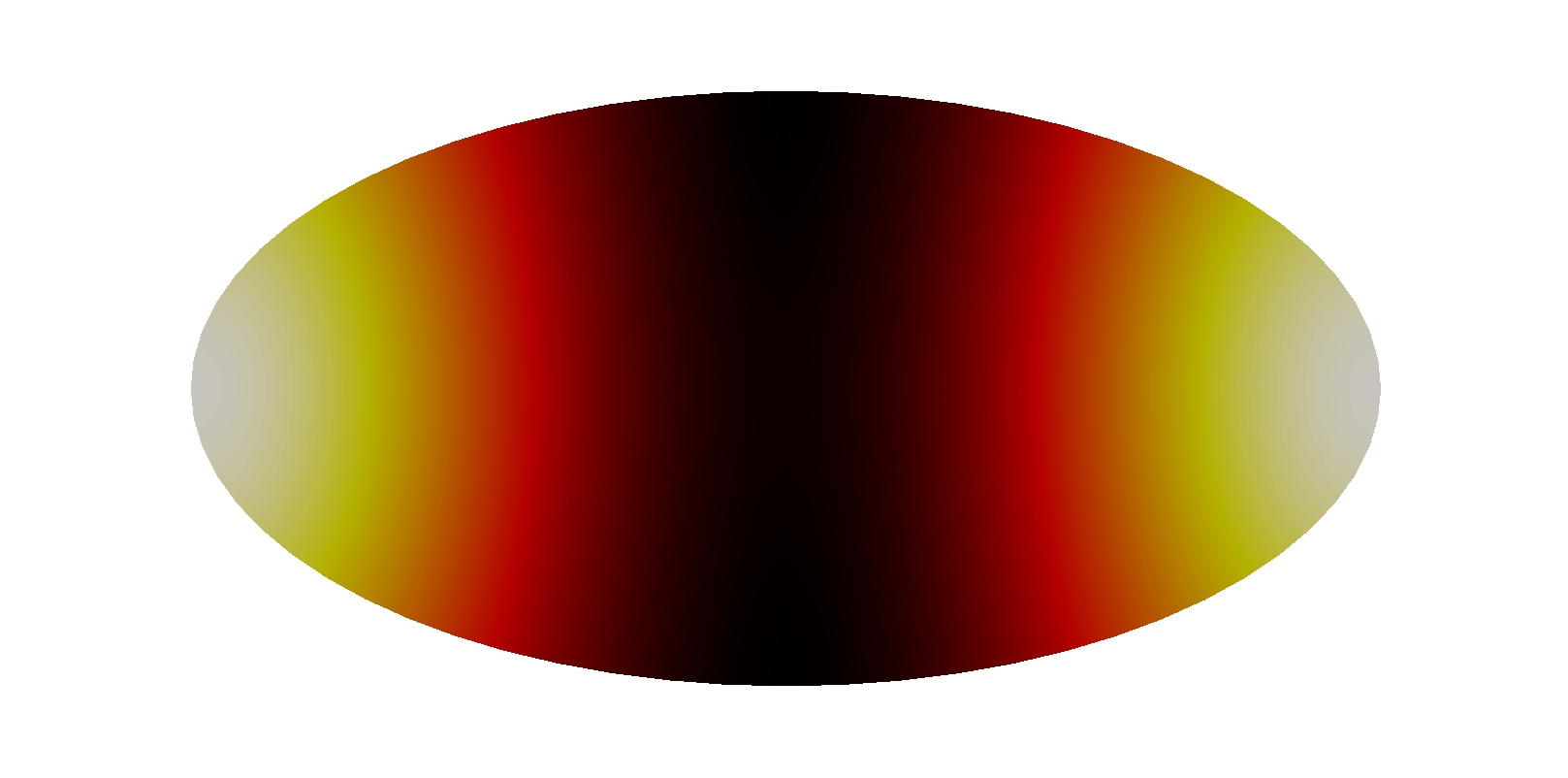}}\;
\subfloat[$\lambda_3=12.52$]{\includegraphics[trim = 70mm 25mm 0 0, clip, width=39mm]{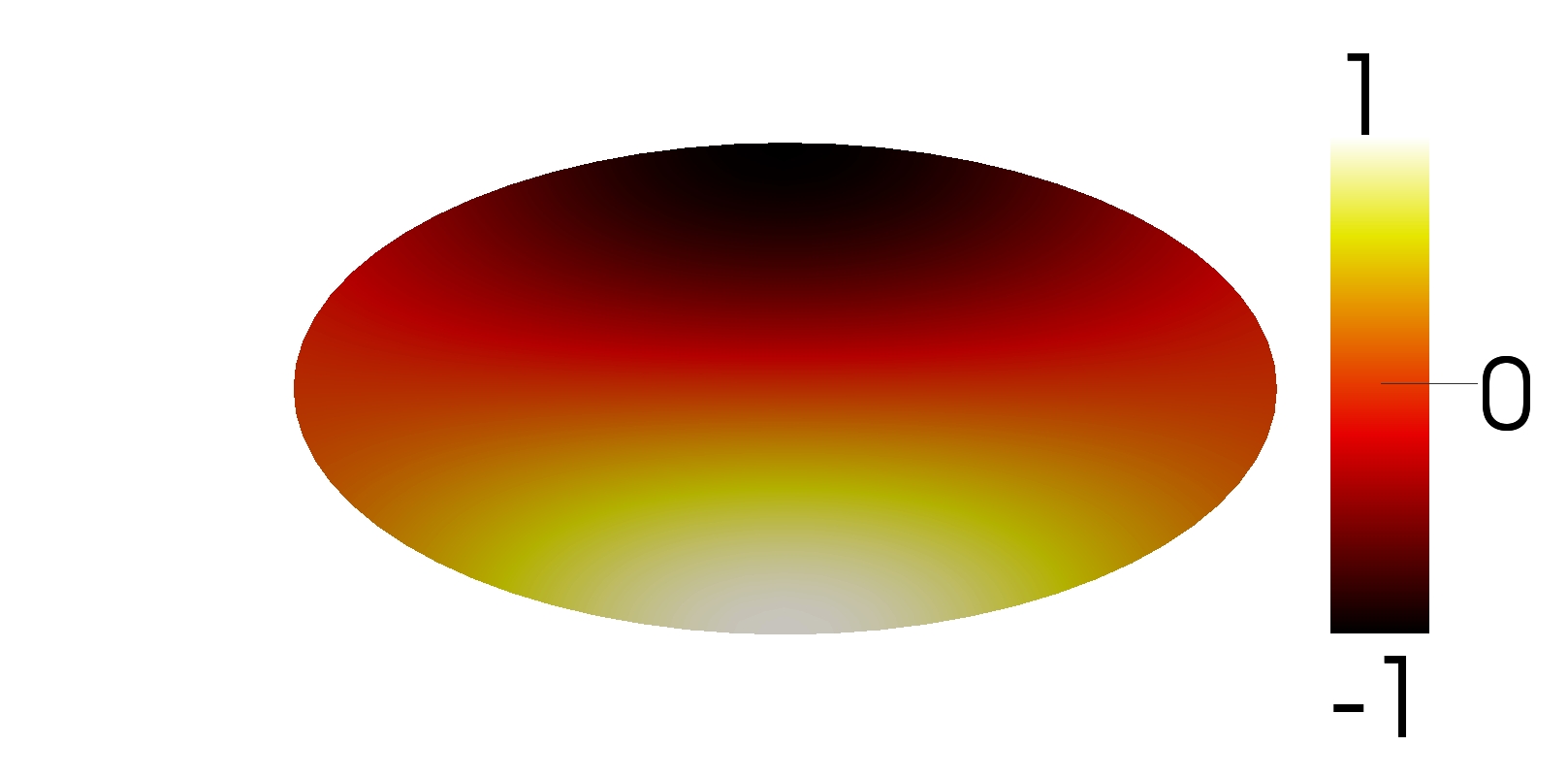}}\\
\subfloat[$\lambda_4=21.63$]{\includegraphics[trim = 30mm 0 30mm 0, clip, width=32mm]{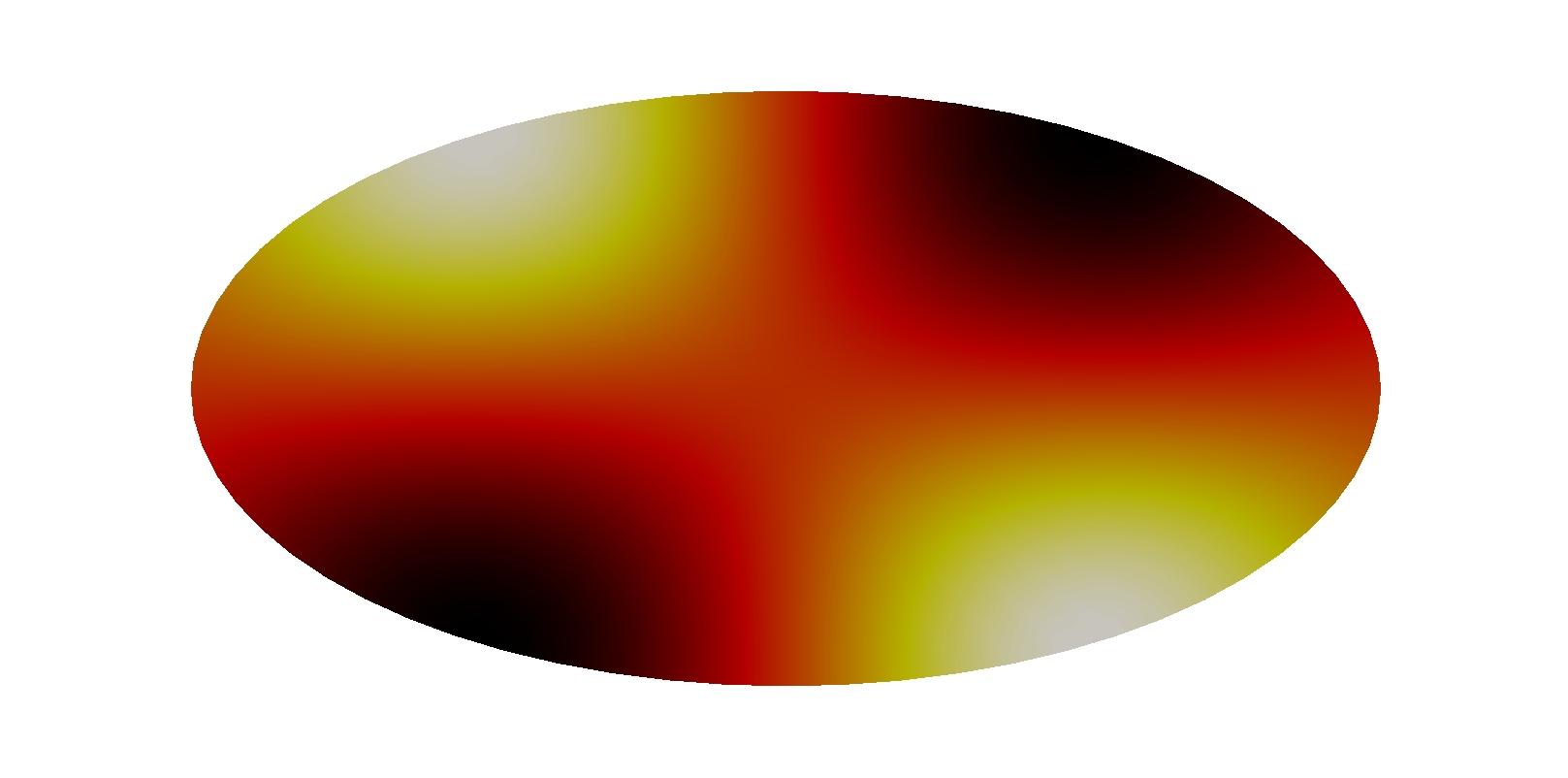}}\;
\subfloat[$\lambda_5=24.51$]{\includegraphics[trim = 30mm 0 30mm 0, clip, width=32mm]{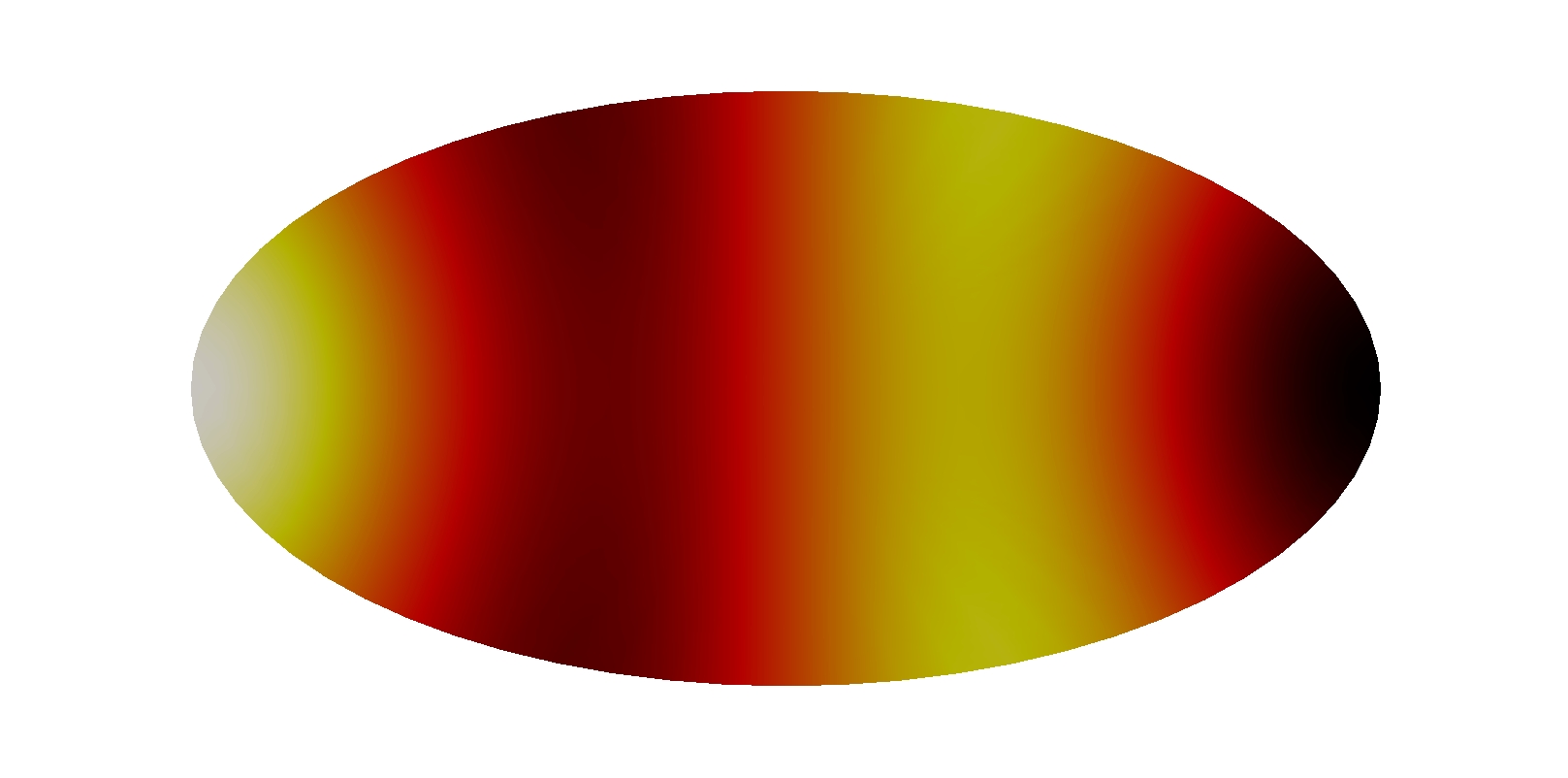}}\;
\subfloat[$\lambda_6=34.30$]{\includegraphics[trim = 30mm 0 30mm 0, clip, width=32mm]{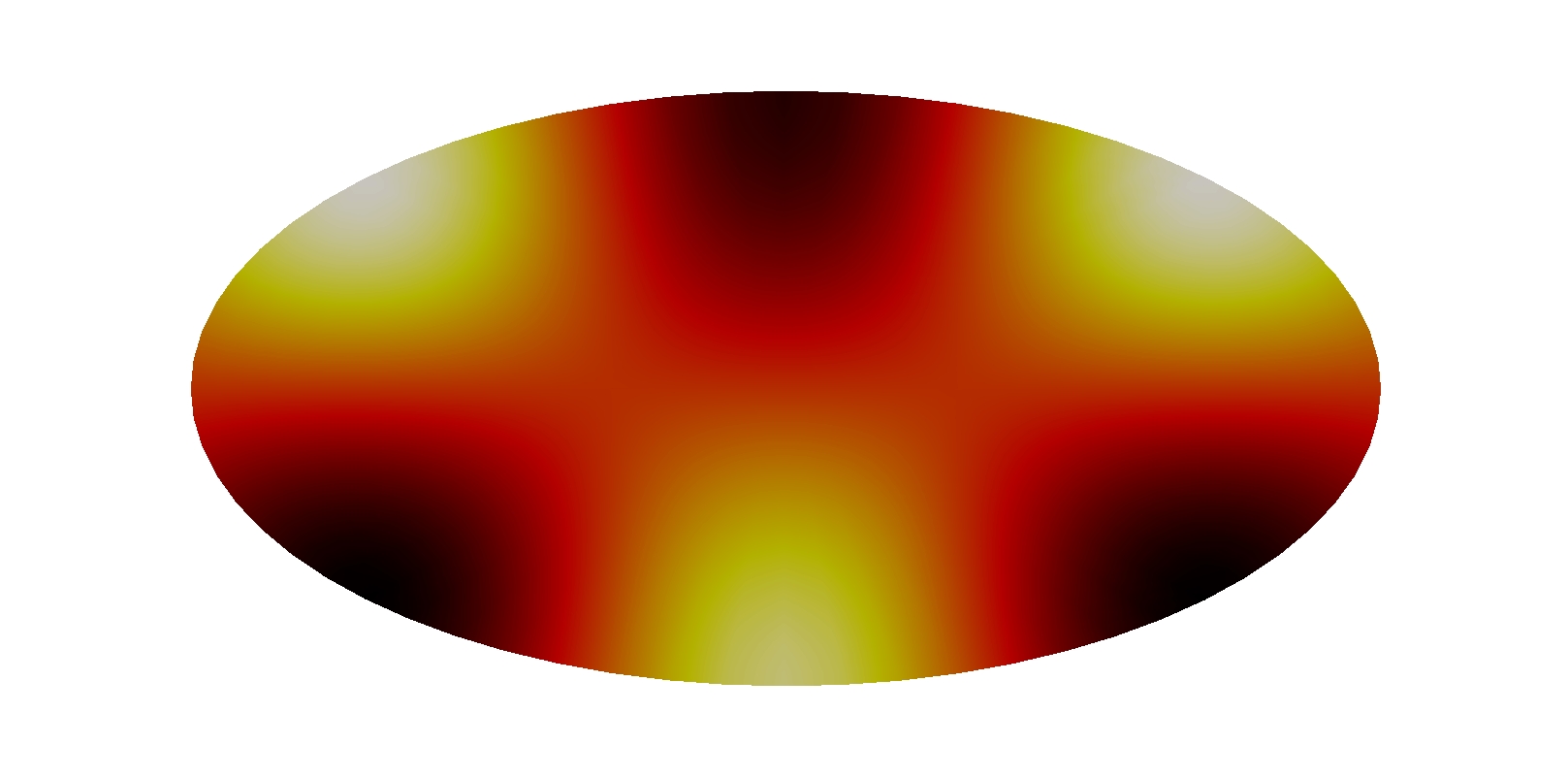}}\\
\subfloat[$\lambda_7=41.75$]{\includegraphics[trim = 30mm 0 30mm 0, clip, width=32mm]{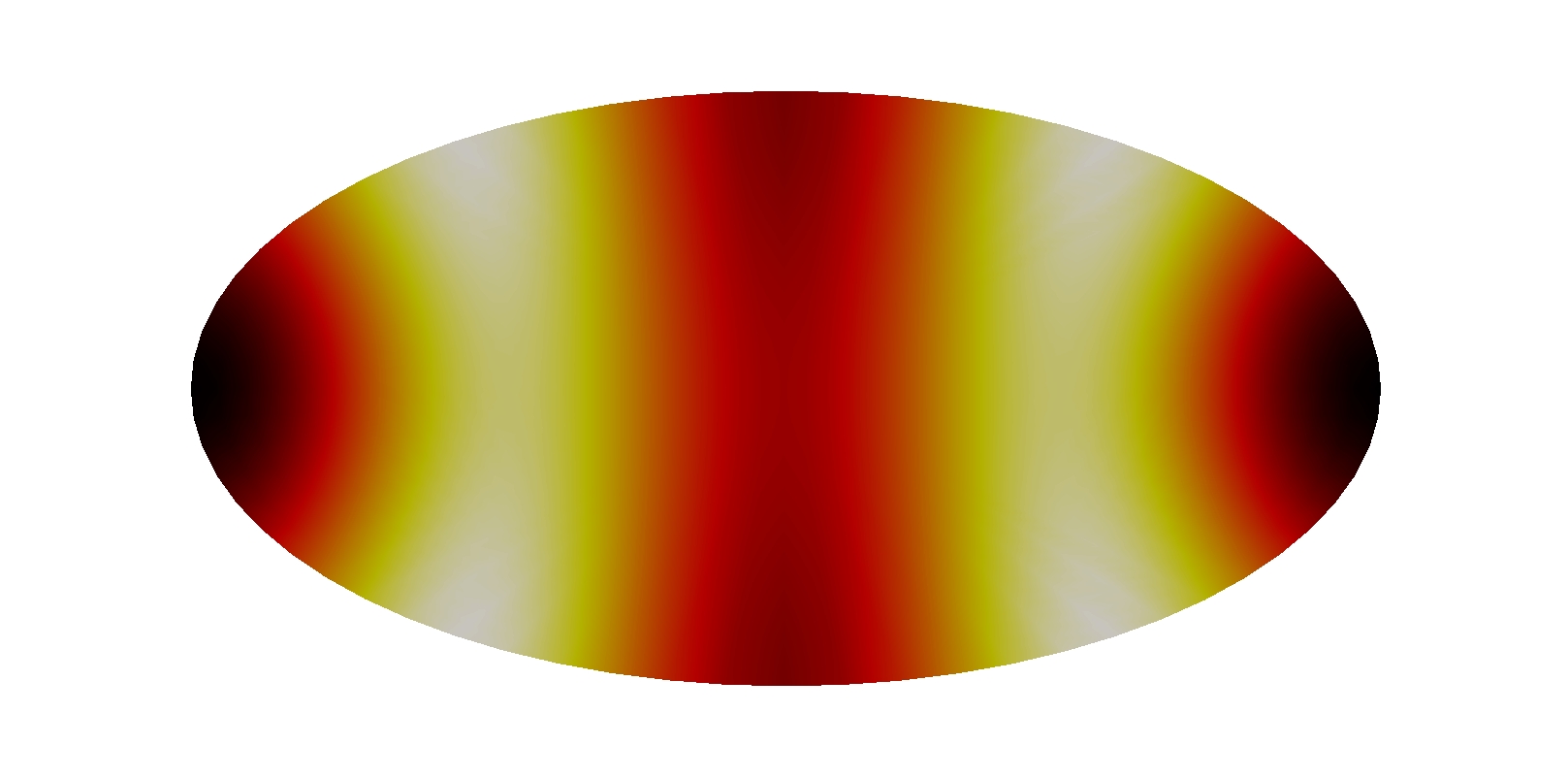}}\;
\subfloat[$\lambda_8=45.88$]{\includegraphics[trim = 30mm 0 30mm 0, clip, width=32mm]{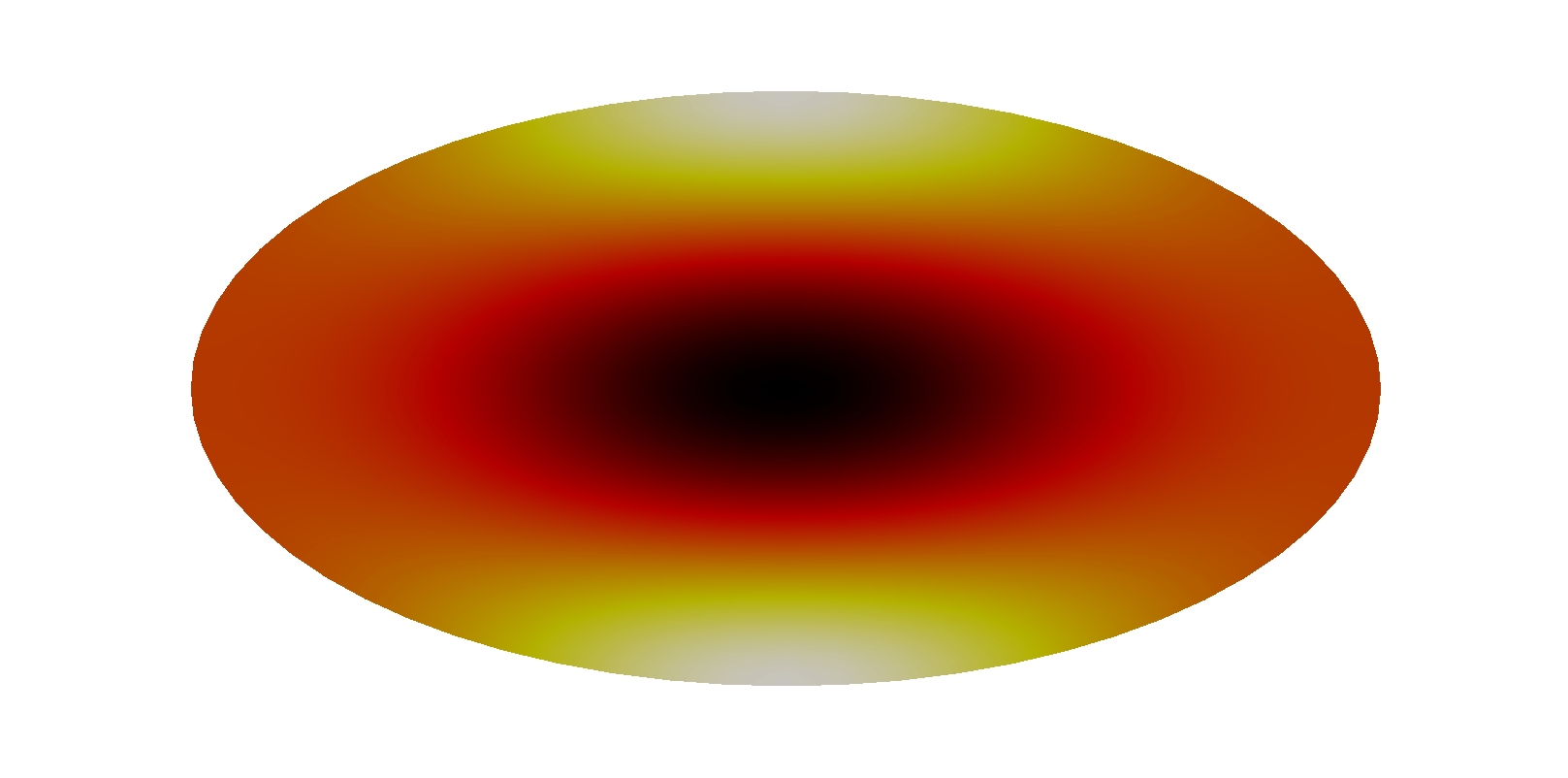}}\;
\subfloat[$\lambda_9=50.97$]{\includegraphics[trim = 30mm 0 30mm 0, width=34mm]{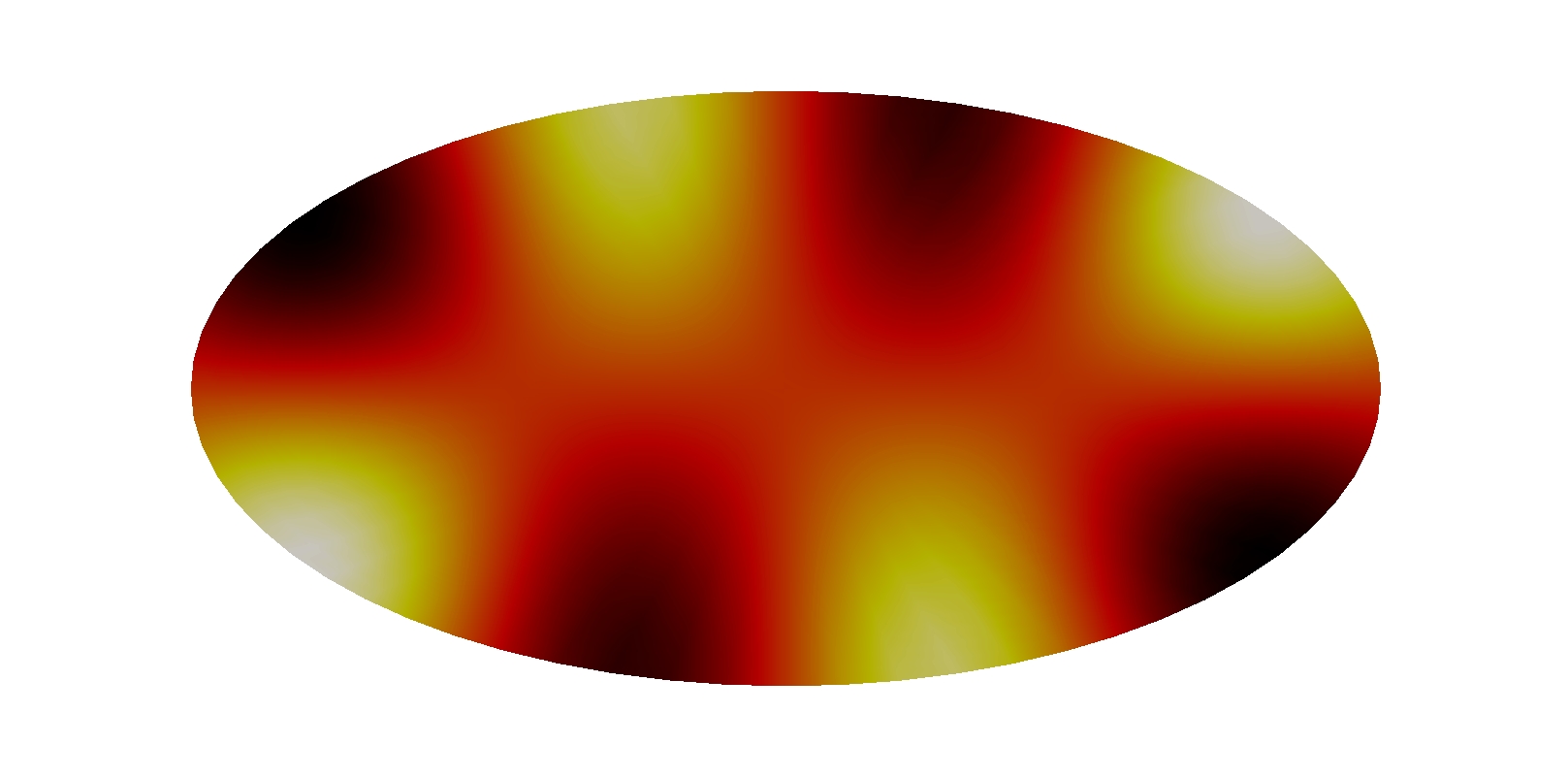}}
\caption{Eigenfunctions corresponding to the labelled eigenvalues on an ellipse. These are solutions of \eqref{eq:geneig} approximated using {\bf deal.II} (Colour version online)}\label{fig:ellipseeigen}
\centering
\subfloat[$d$=10, $\gamma$=10]{\includegraphics[trim = 100mm 0 0 0, clip, width=36mm]{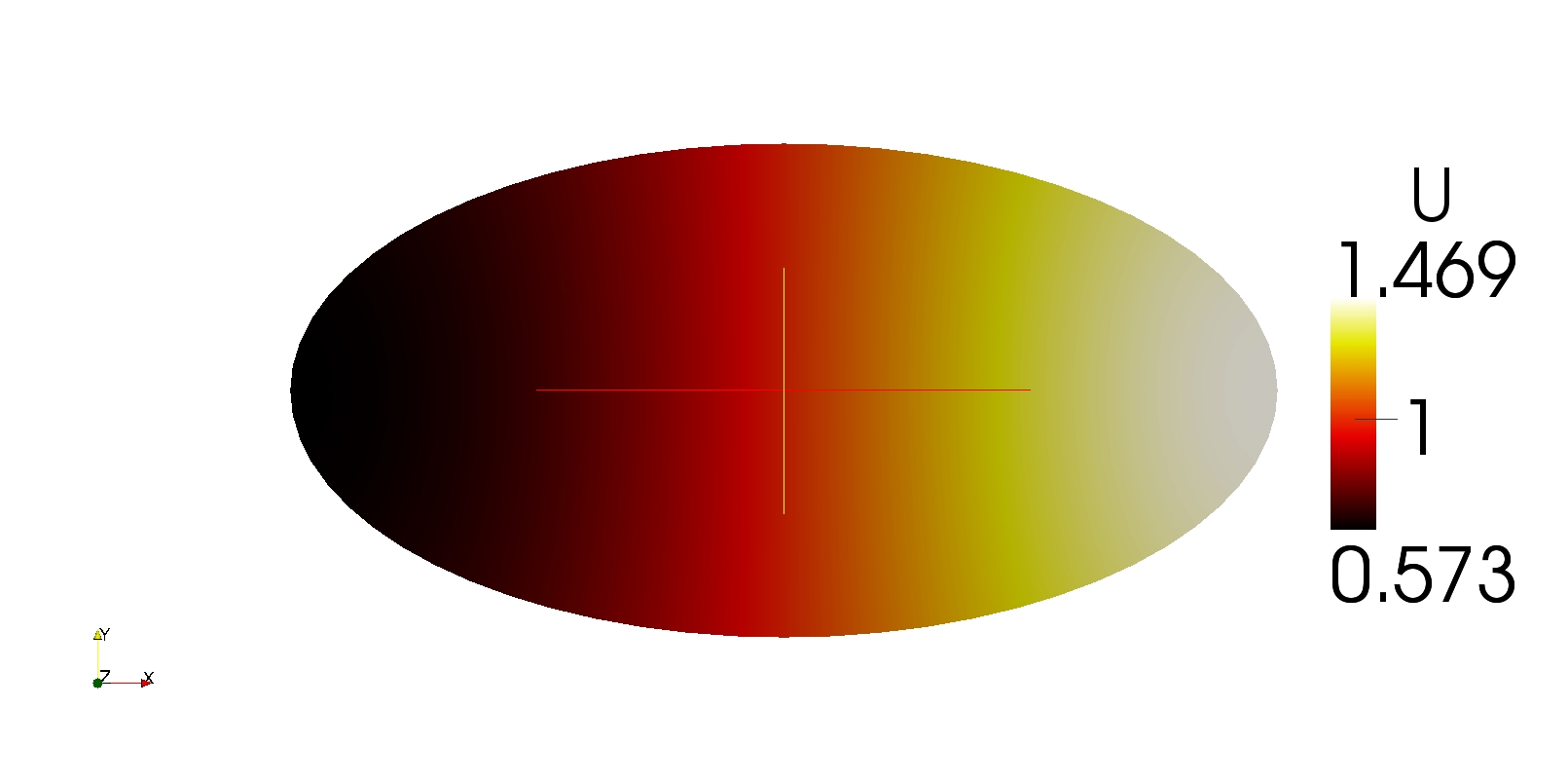}}\;
\subfloat[$d$=8.8, $\gamma$=30]{\includegraphics[trim = 100mm 0 0 0, clip, width=36mm]{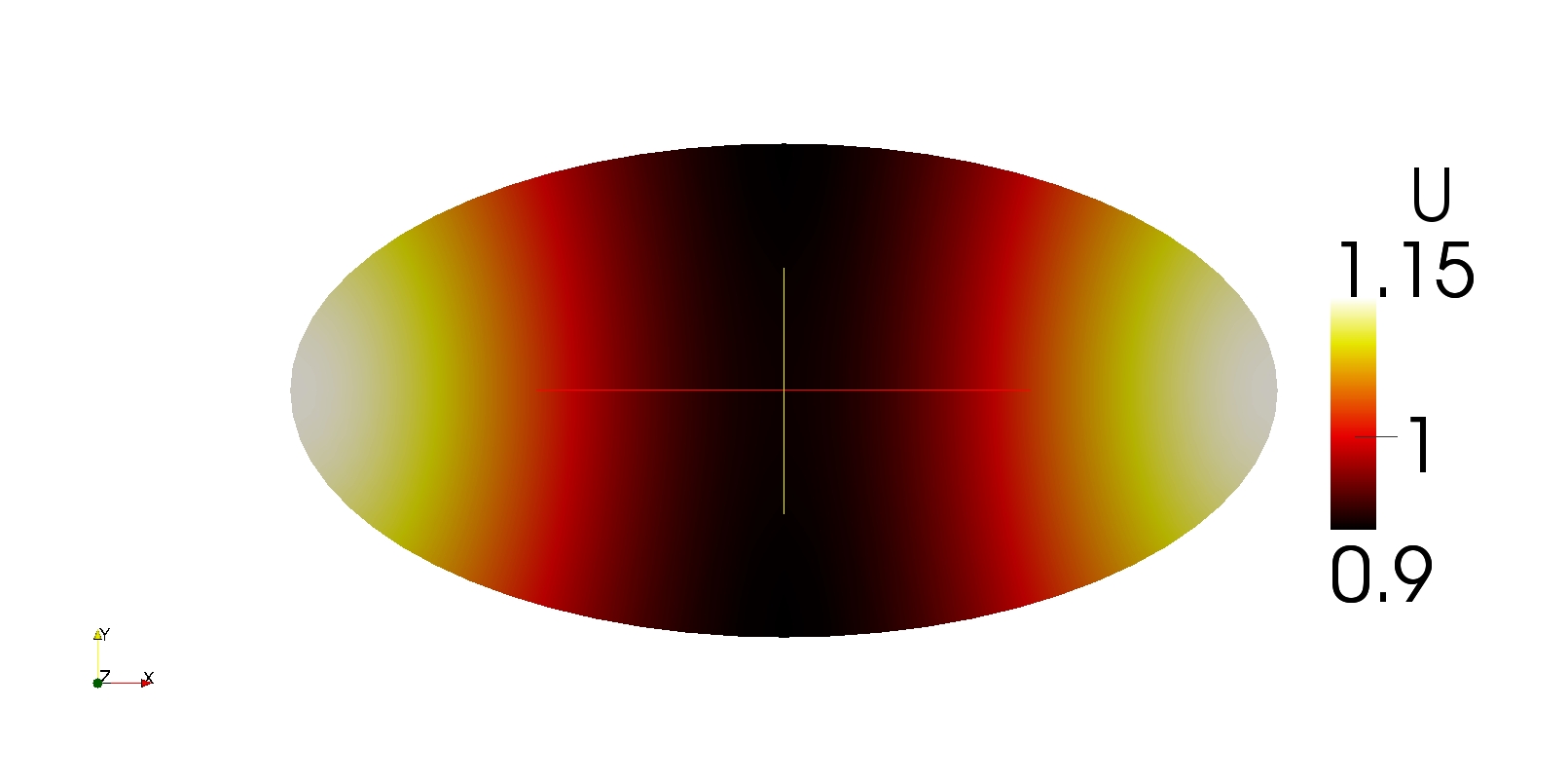}\label{fig:ell30}}\;
\subfloat[$d$=8.8, $\gamma$=44]{\includegraphics[trim = 100mm 0 0 0, clip, width=36mm]{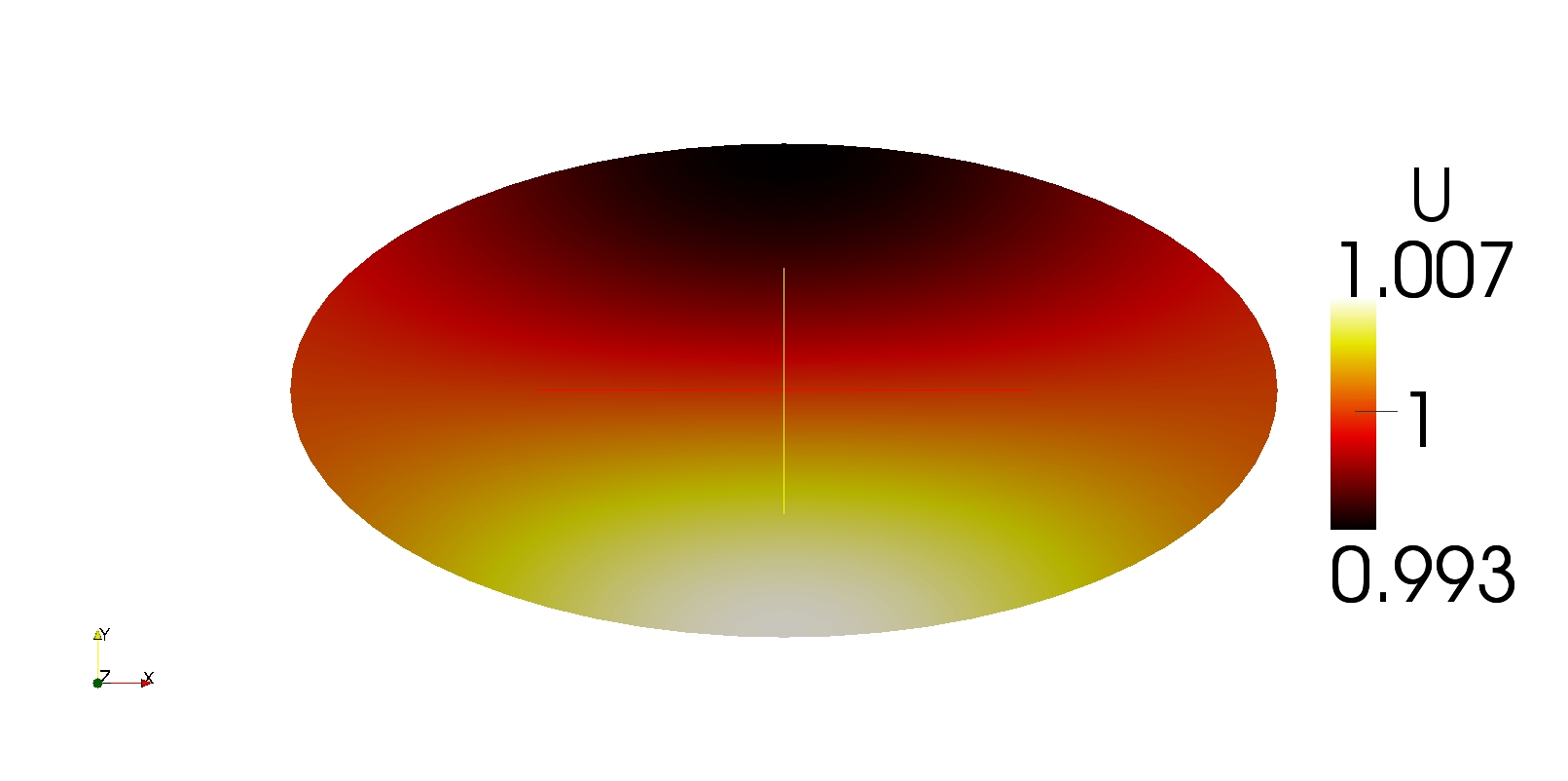}}\\
\subfloat[$d$=8.8, $\gamma$=57]{\includegraphics[trim = 100mm 0 0 0, clip, width=36mm]{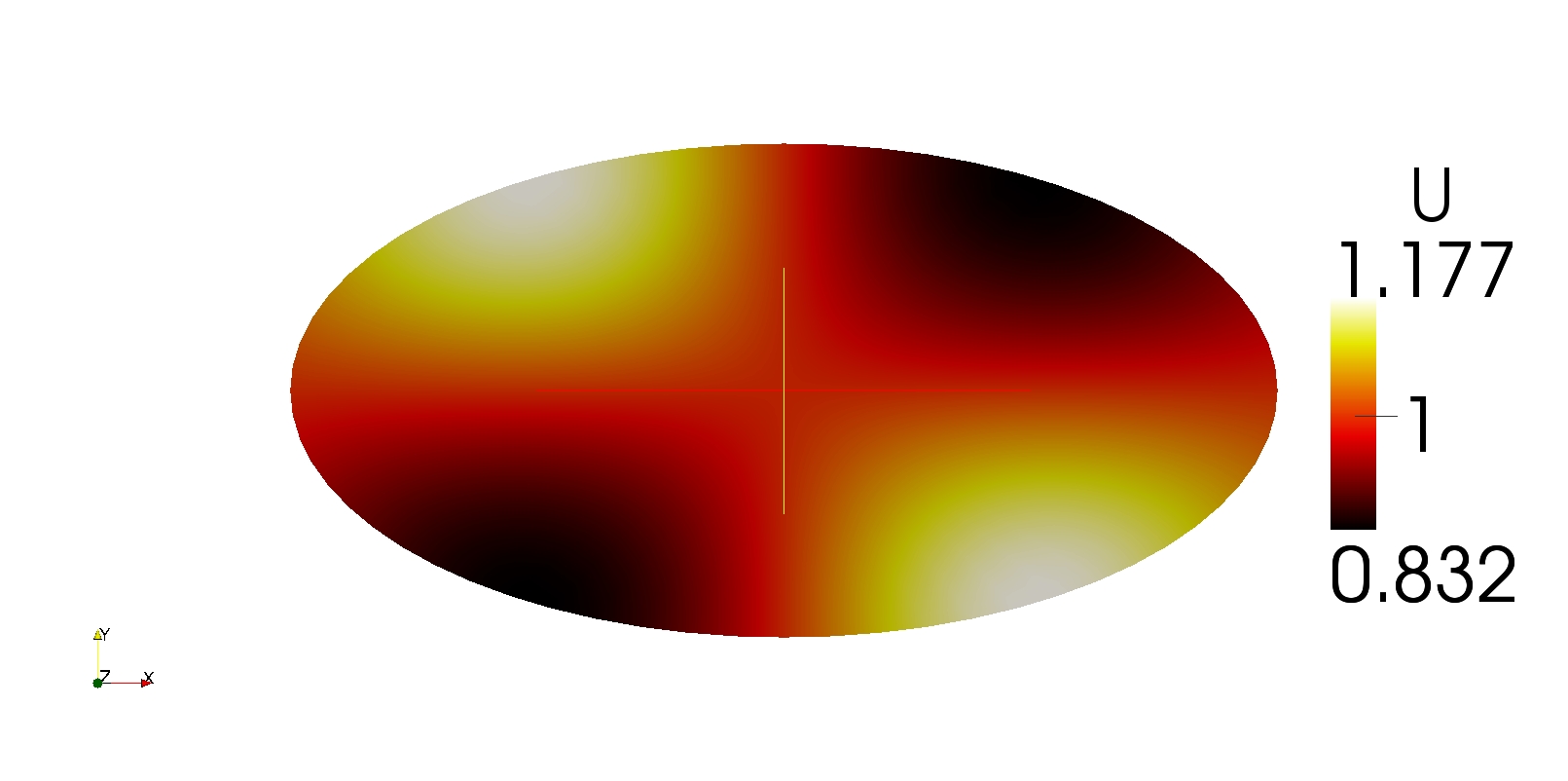}}\;
\subfloat[$d$=8.7, $\gamma$=77]{\includegraphics[trim = 100mm 0 0 0, clip, width=36mm]{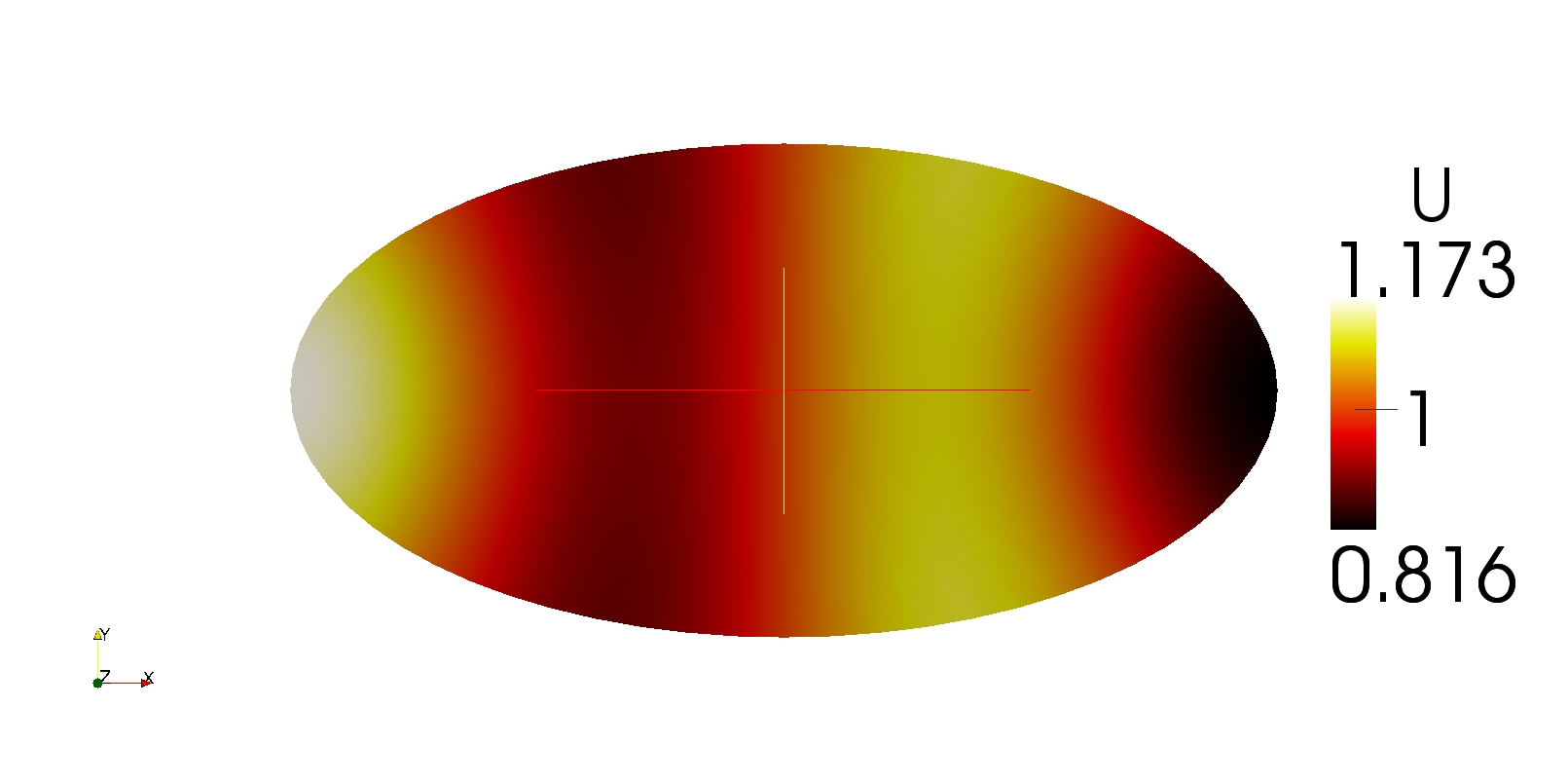}}\;
\subfloat[$d$=8.7, $\gamma$=95]{\includegraphics[trim = 100mm 0 0 0, clip, width=36mm]{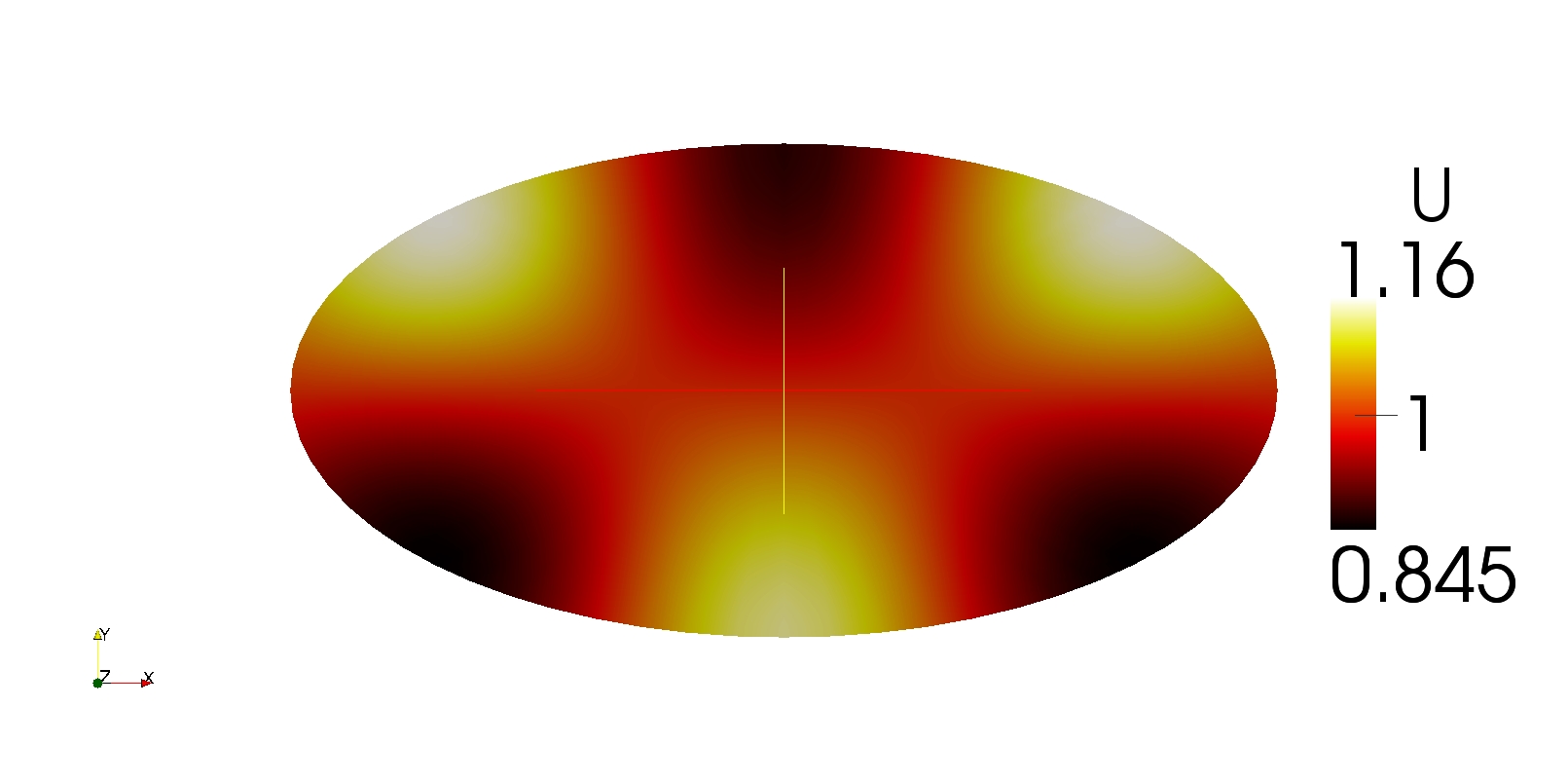}}\\
\subfloat[$d$=8.63, $\gamma$=115]{\includegraphics[trim = 100mm 0 0 0, clip, width=36mm]{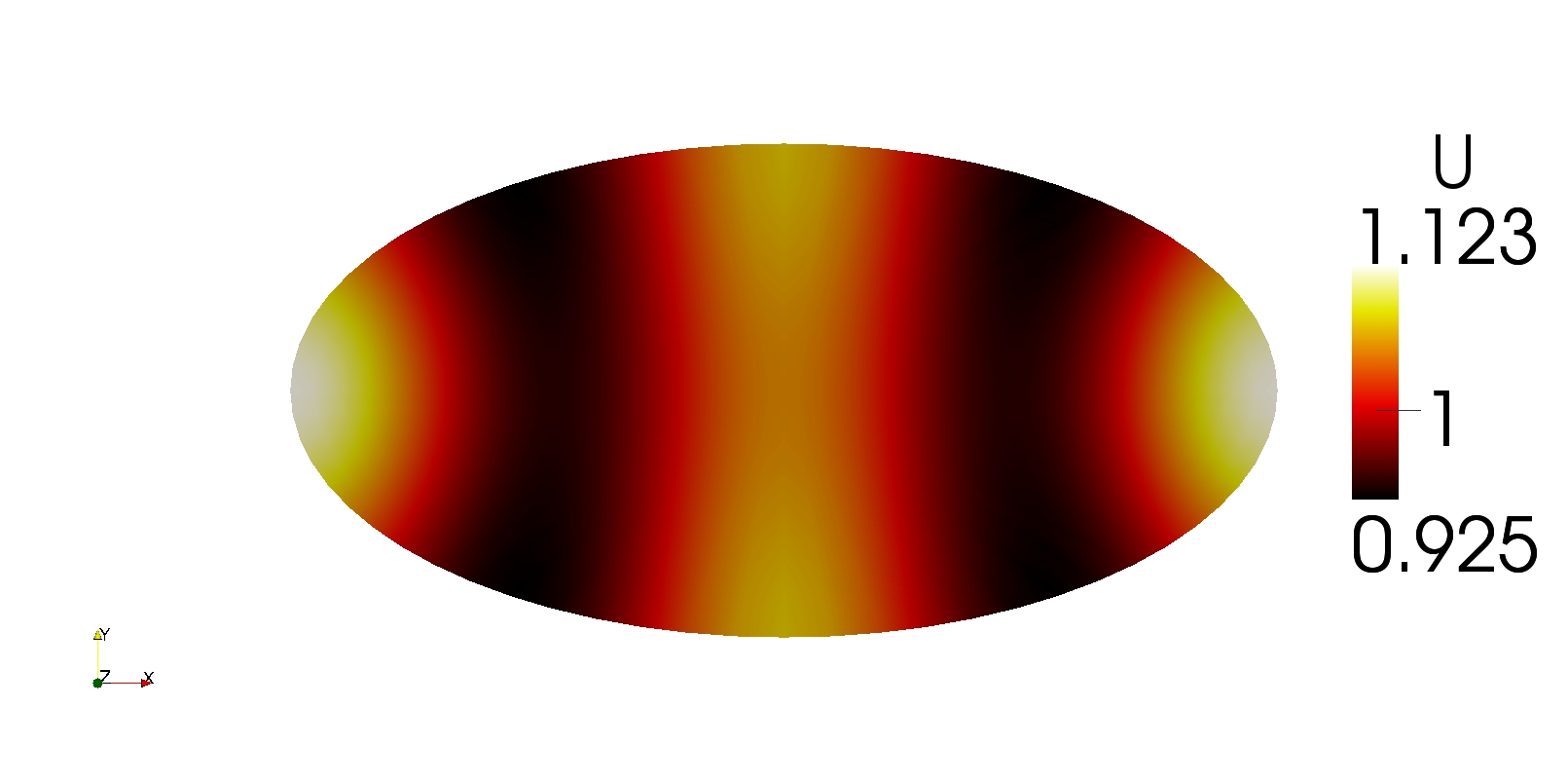}}\;
\subfloat[$d$=8.61, $\gamma$=135]{\includegraphics[trim = 100mm 0 0 0, clip, width=36mm]{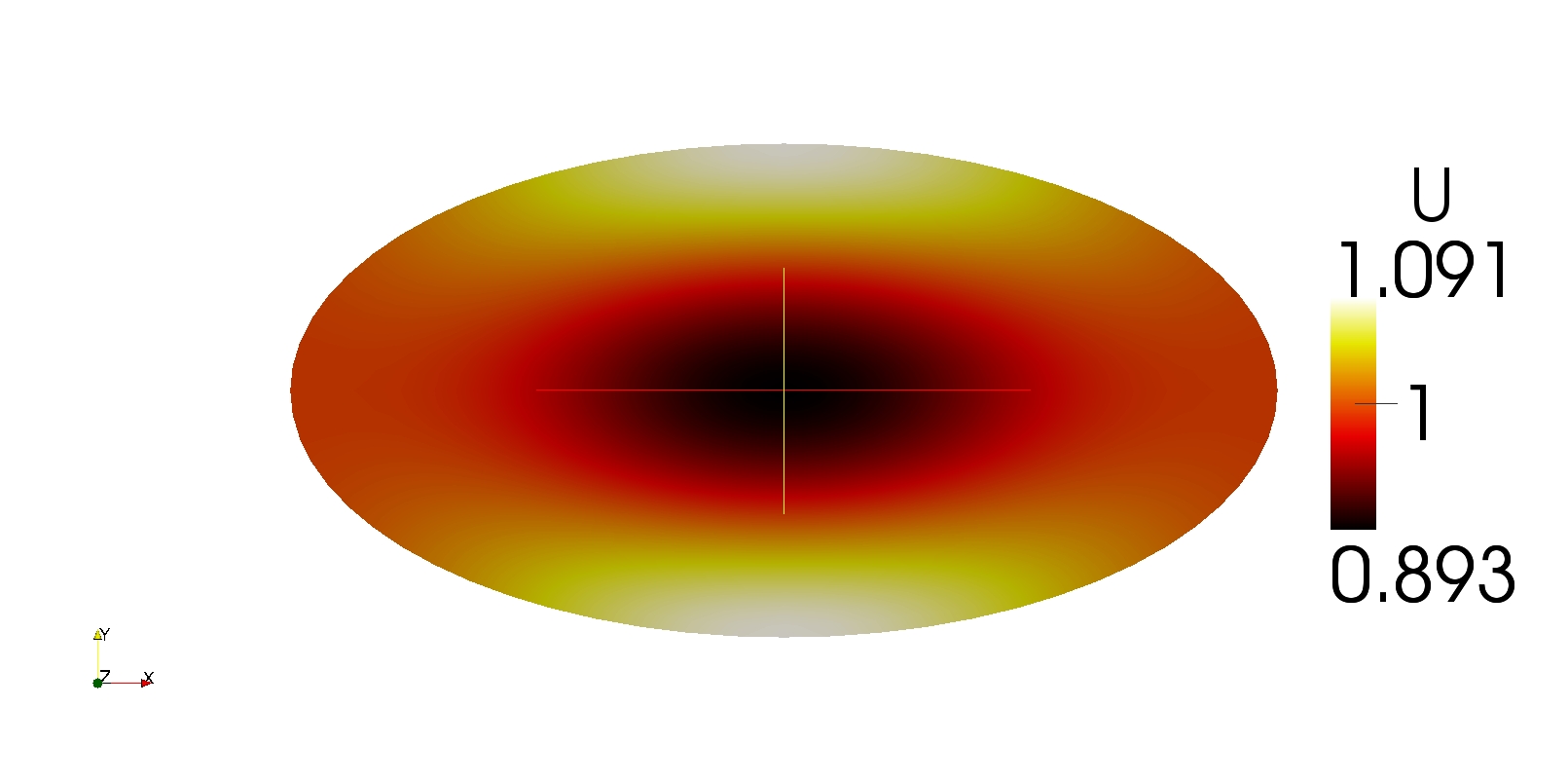}}\;
\subfloat[$d$=8.61, $\gamma$=150]{\includegraphics[trim = 100mm 0 0 0, clip, width=36mm]{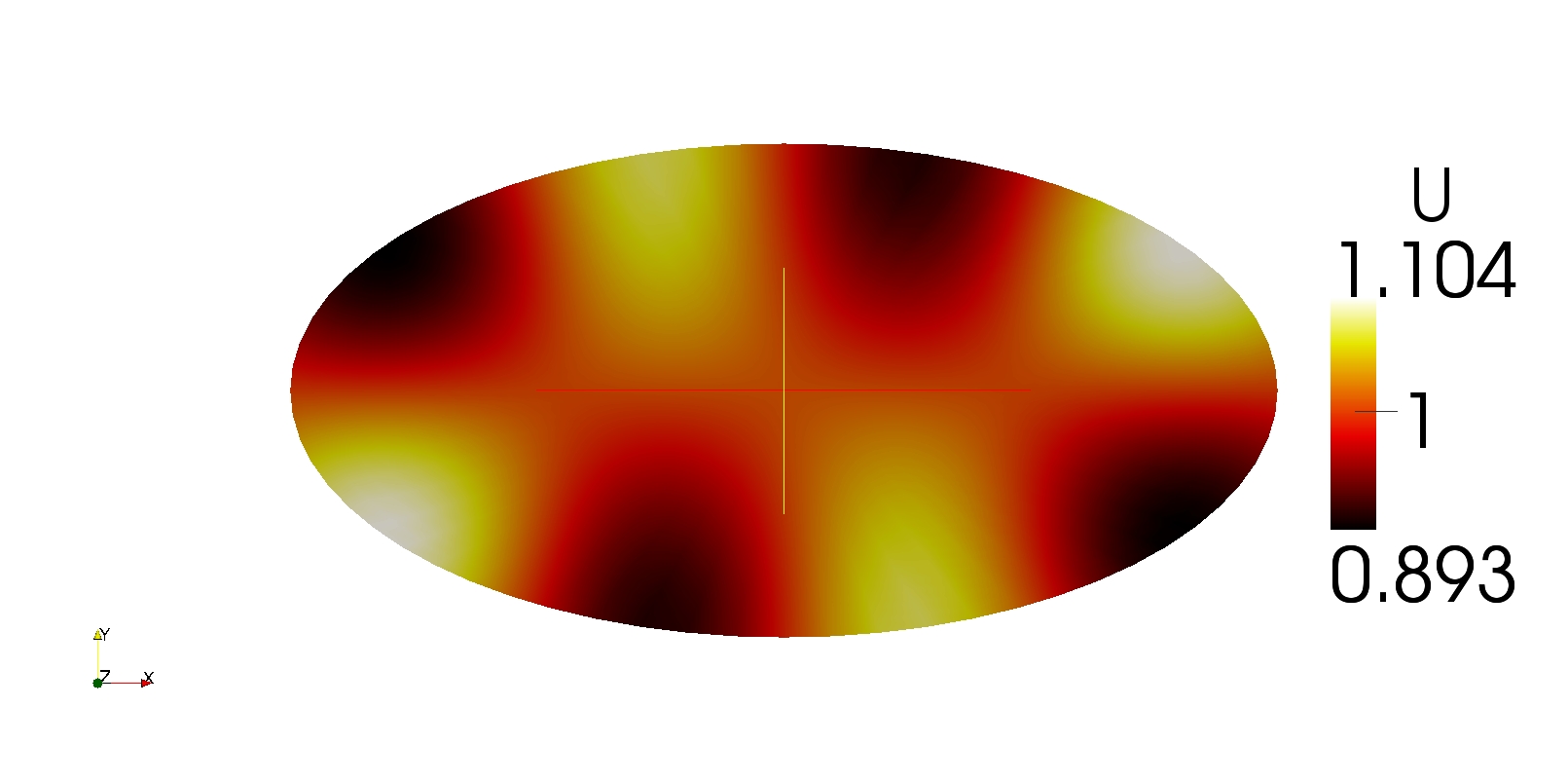}}
\caption{Converged solutions of system \eqref{eq:non}, with Schnakenberg kinetics \eqref{eq:schnakkin}, on an ellipse for the species $u$, they all match the associated eigenfunctions shown in Figure \ref{fig:ellipseeigen} (Colour version online)}\label{fig:ellipserd}
\end{figure}

\subsection{Example 2: Ellipse}
Eigenmodes on an ellipse have been investigated in various articles \citep{fox,greben,neves,wu}. Finding the solution involves numerically solving the Mathieu and modified Mathieu equations \citep{hbk}. In particular \cite{wu} analytically find the first eigenvalue of ellipses with Dirichlet boundary conditions, of various sizes of ellipse. Using the eigenvalue solver described in Section \ref{sec:genisolation}, with Dirichlet boundary conditions, we can reproduce their results (results not reported in the interests of brevity). In the following we consider Neumann conditions and choose the semimajor axis to be twice the semiminor axis. The eigenvalues and eigenfunctions are shown in Figure \ref{fig:ellipseeigen}. Figure \ref{fig:ellipserd} shows the converged solutions of the reaction diffusion system when the chosen values of $d$ and $\gamma$ isolate the corresponding wavenumbers $k^2_i=\lambda_i$.

\subsection{Example 3: Dumbbell}
\begin{figure}
\centering
\subfloat[$\lambda_1=1.49$]{\includegraphics[width=38mm]{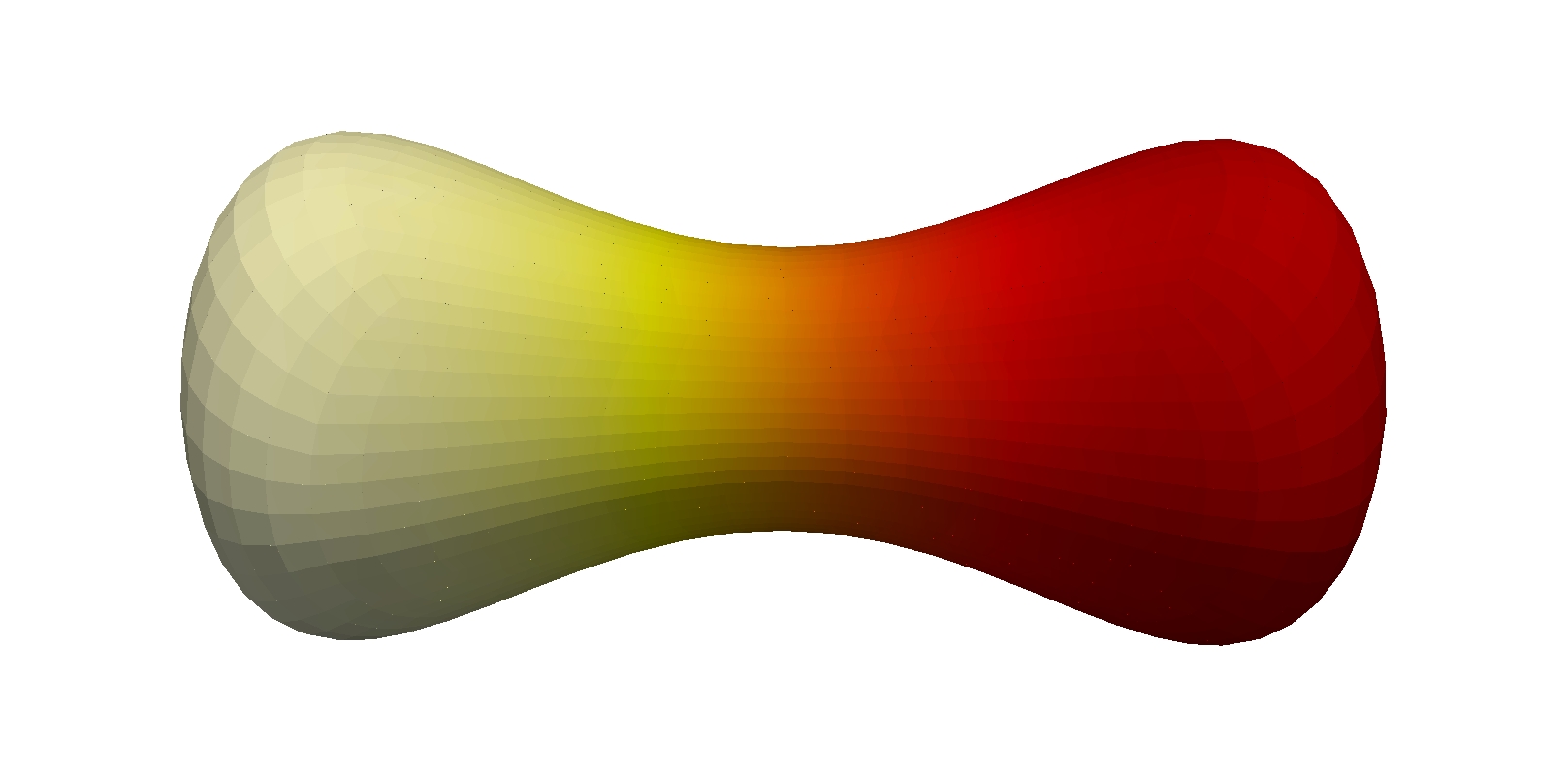}}\;
\subfloat[$\lambda_2=12.68$]{\includegraphics[width=38mm]{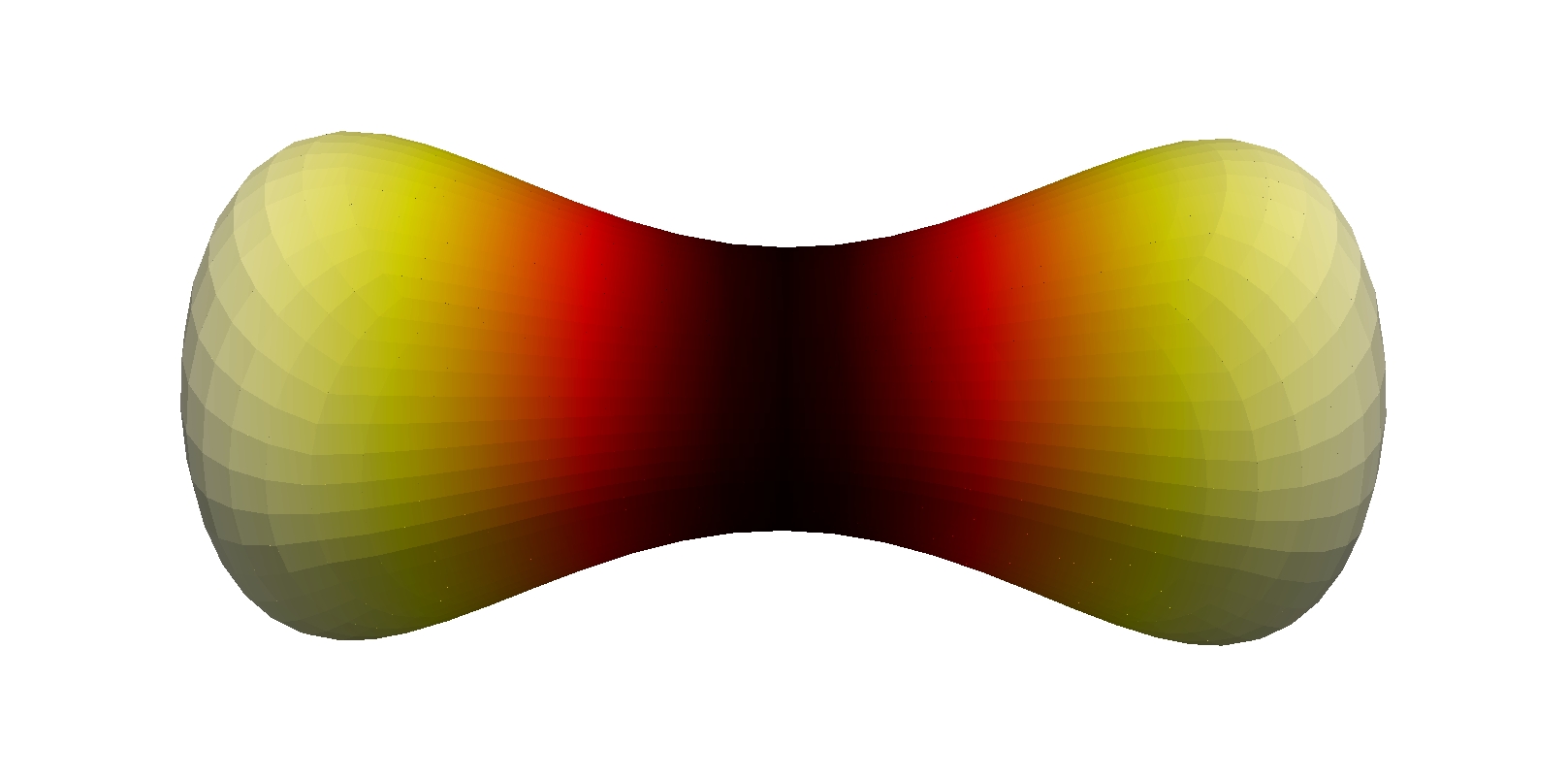}}\;
\subfloat[$\lambda_3=22.86$]{\includegraphics[trim = 80mm 15mm 40mm 0, clip, width=38mm]{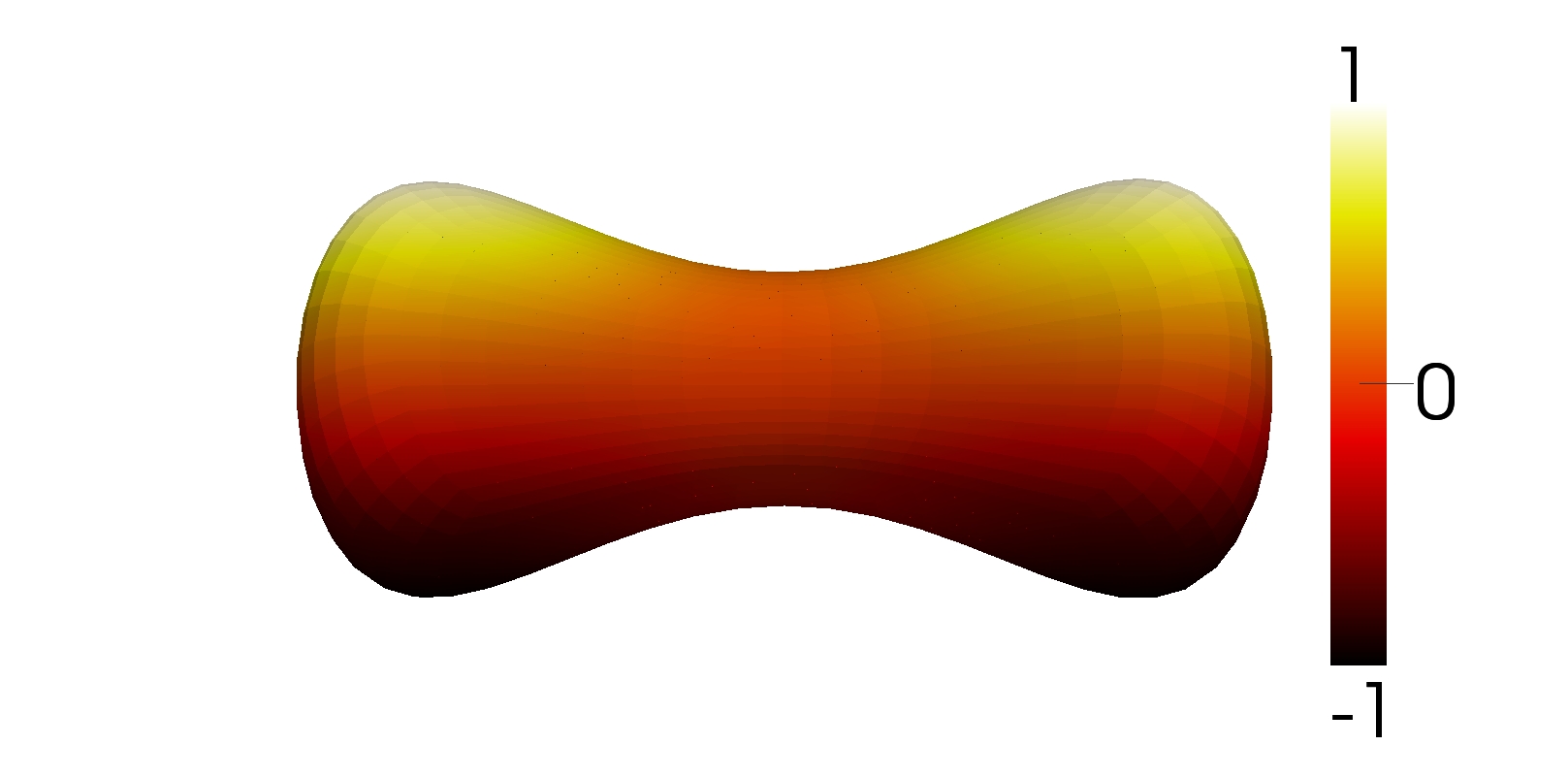}}\\
\subfloat[$\lambda_4=22.98$]{\includegraphics[width=38mm]{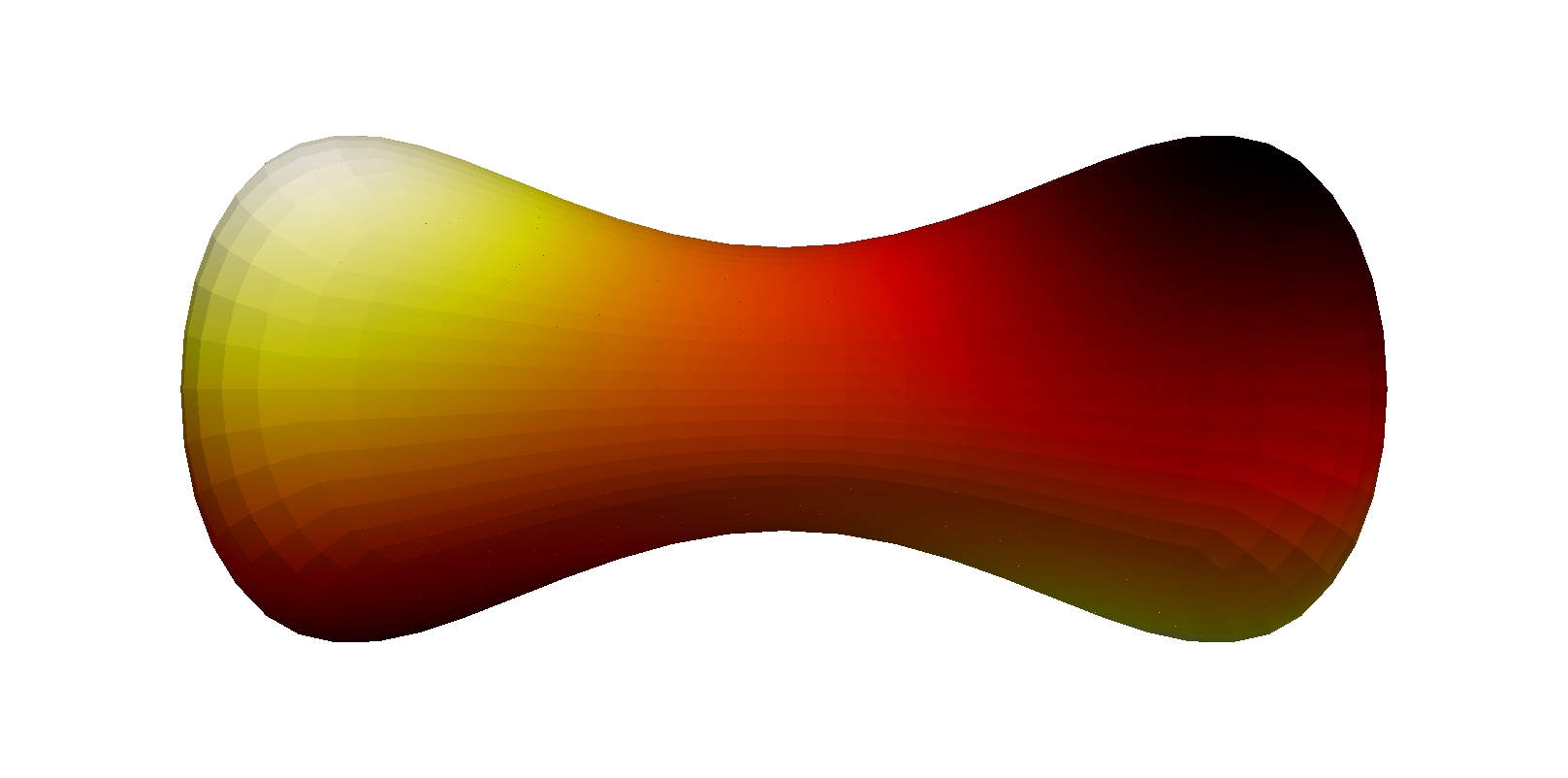}}\;
\subfloat[$\lambda_5=26.52$]{\includegraphics[width=38mm]{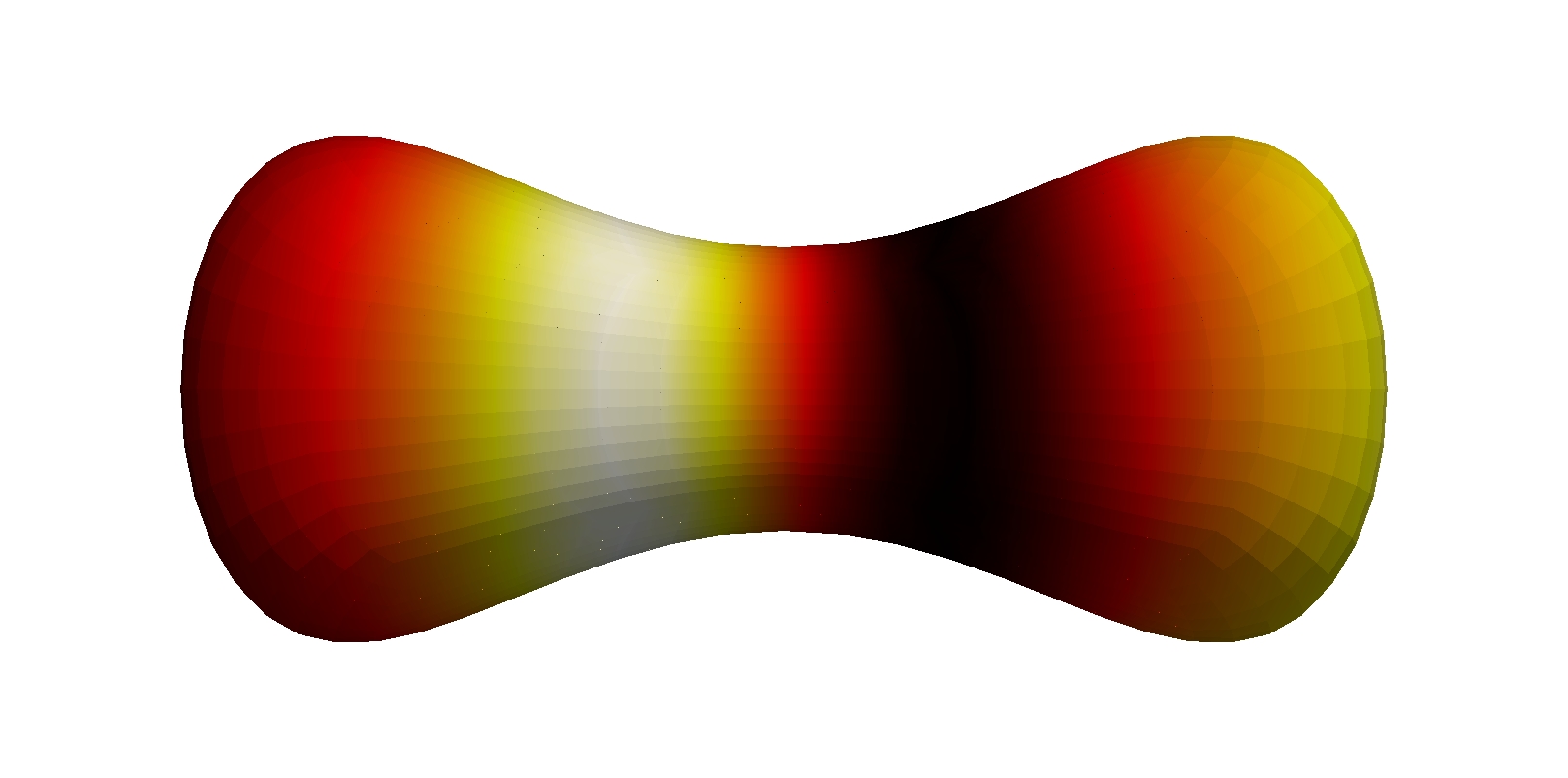}}\;
\subfloat[$\lambda_6=49.91$]{\includegraphics[width=38mm]{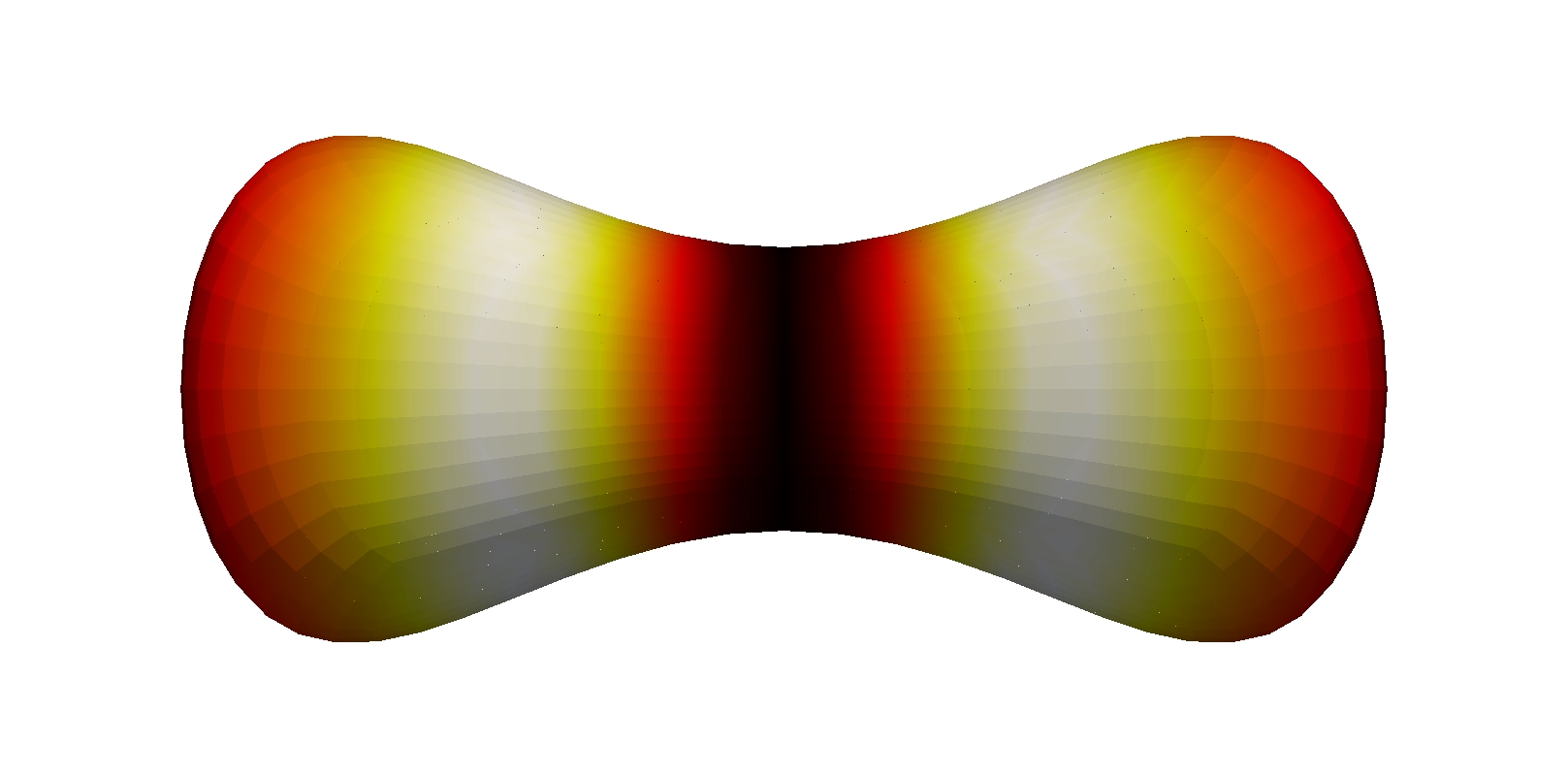}}\;
\caption{Eigenfunctions corresponding to the labelled eigenvalues on the dumbell. These are solutions of \eqref{eq:geneig} approximated using {\bf deal.II} (Colour version online)}\label{fig:dumeigen}
%
\centering
\subfloat[d=10, $\gamma=5$]{\includegraphics[trim = 90mm 15mm 20mm 0, clip, width=38mm]{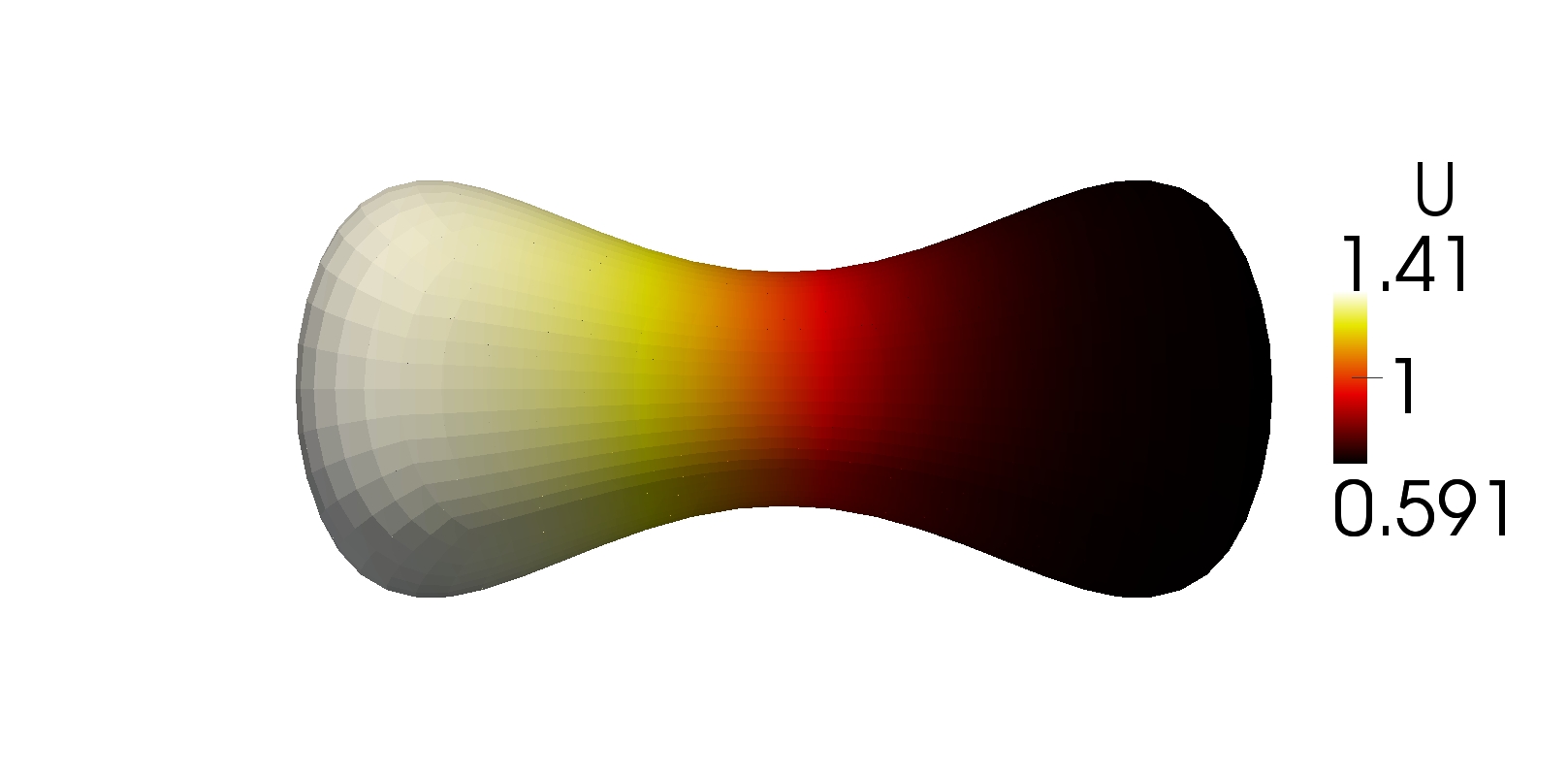}}\;
\subfloat[d=9, $\gamma=40$]{\includegraphics[trim = 90mm 15mm 10mm 0, clip, width=38mm]{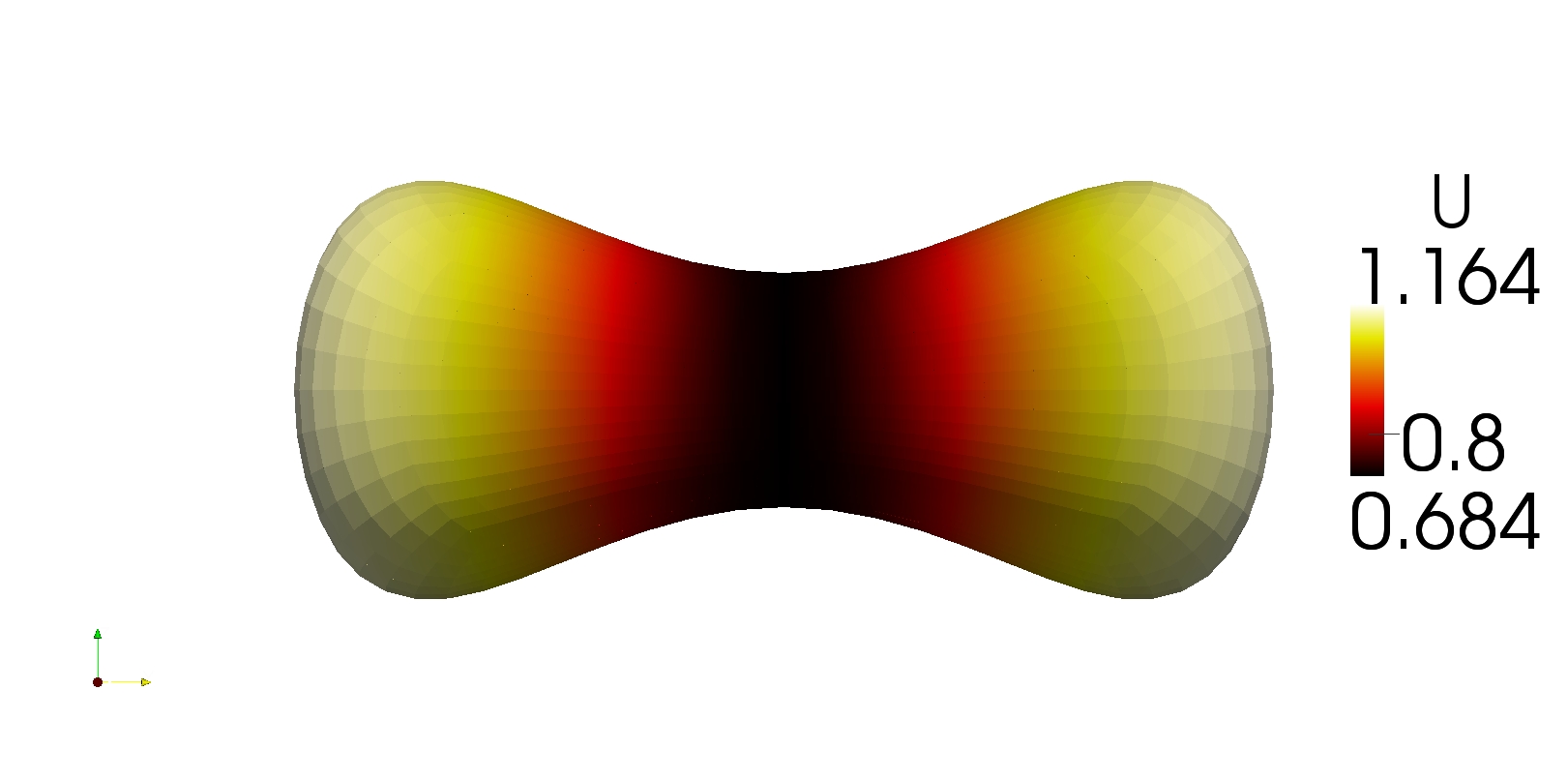}}\;
\subfloat[d=8.8, $\gamma=60$]{\includegraphics[trim = 90mm 15mm 20mm 0, clip, width=38mm]{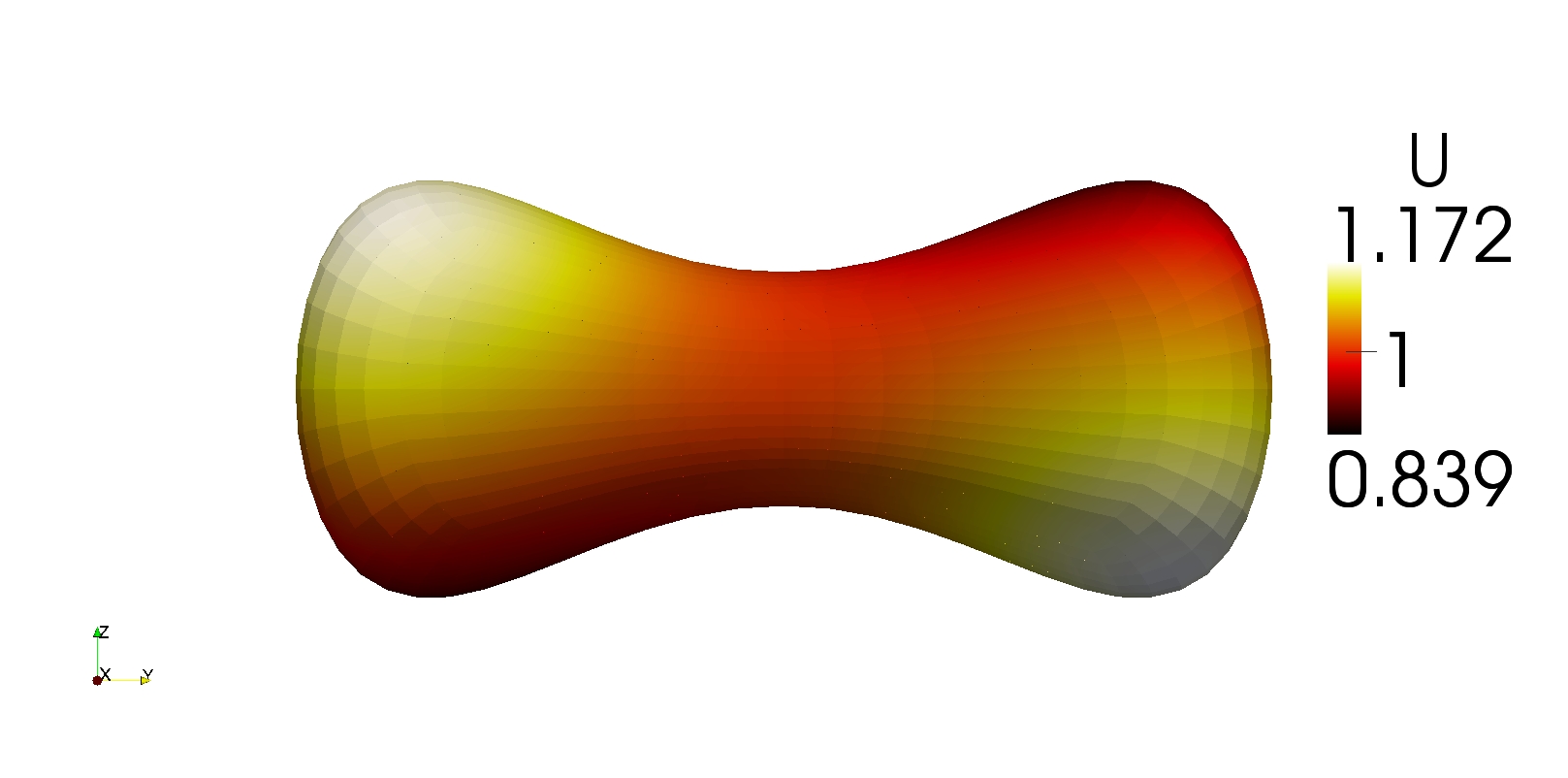}}\;
\subfloat[d=8.8, $\gamma=88$]{\includegraphics[trim = 90mm 15mm 10mm 0, clip, width=38mm]{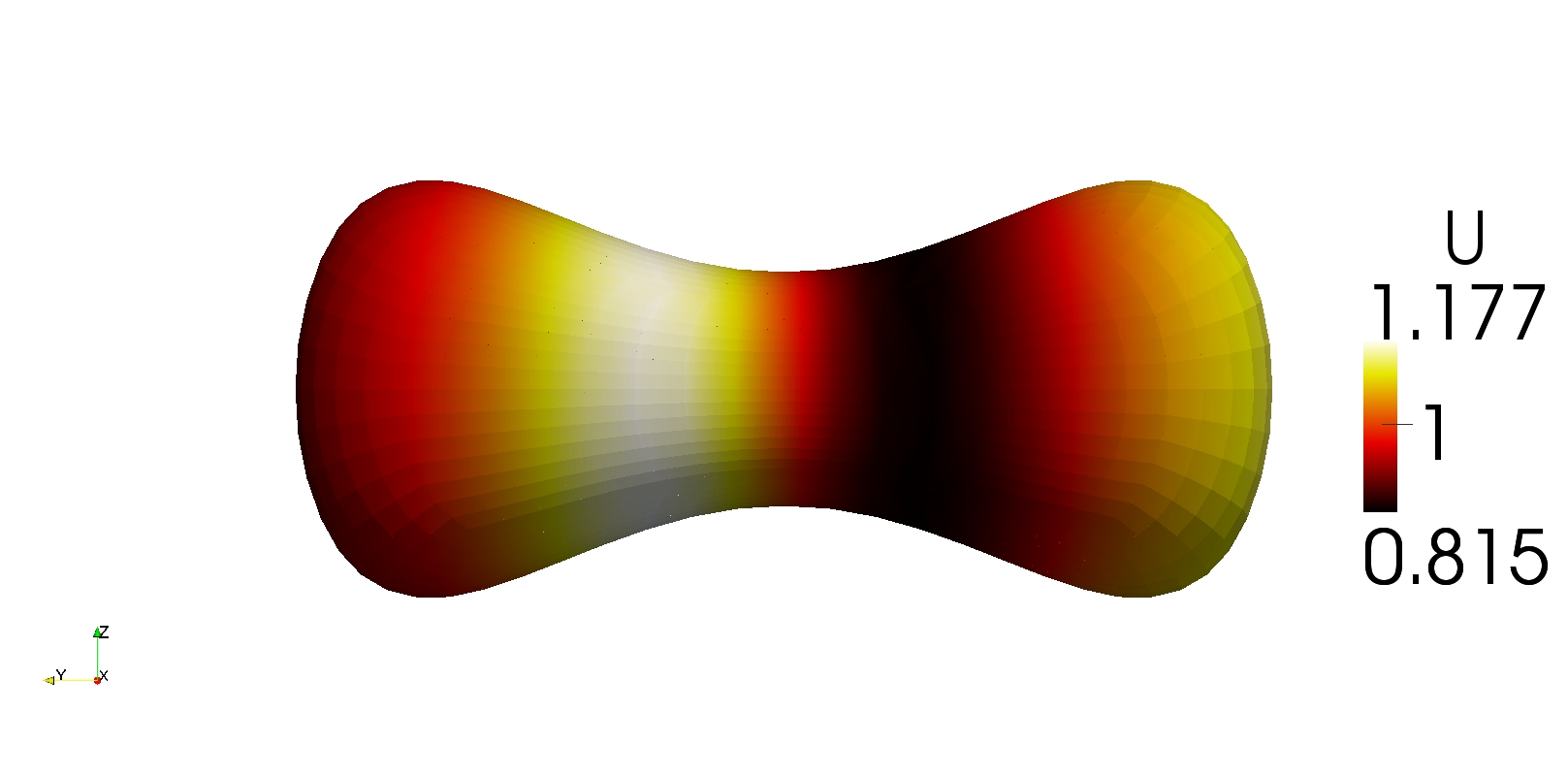}}\;
\subfloat[d=8.65, $\gamma=130$]{\includegraphics[trim = 90mm 15mm 10mm 0, clip, width=38mm]{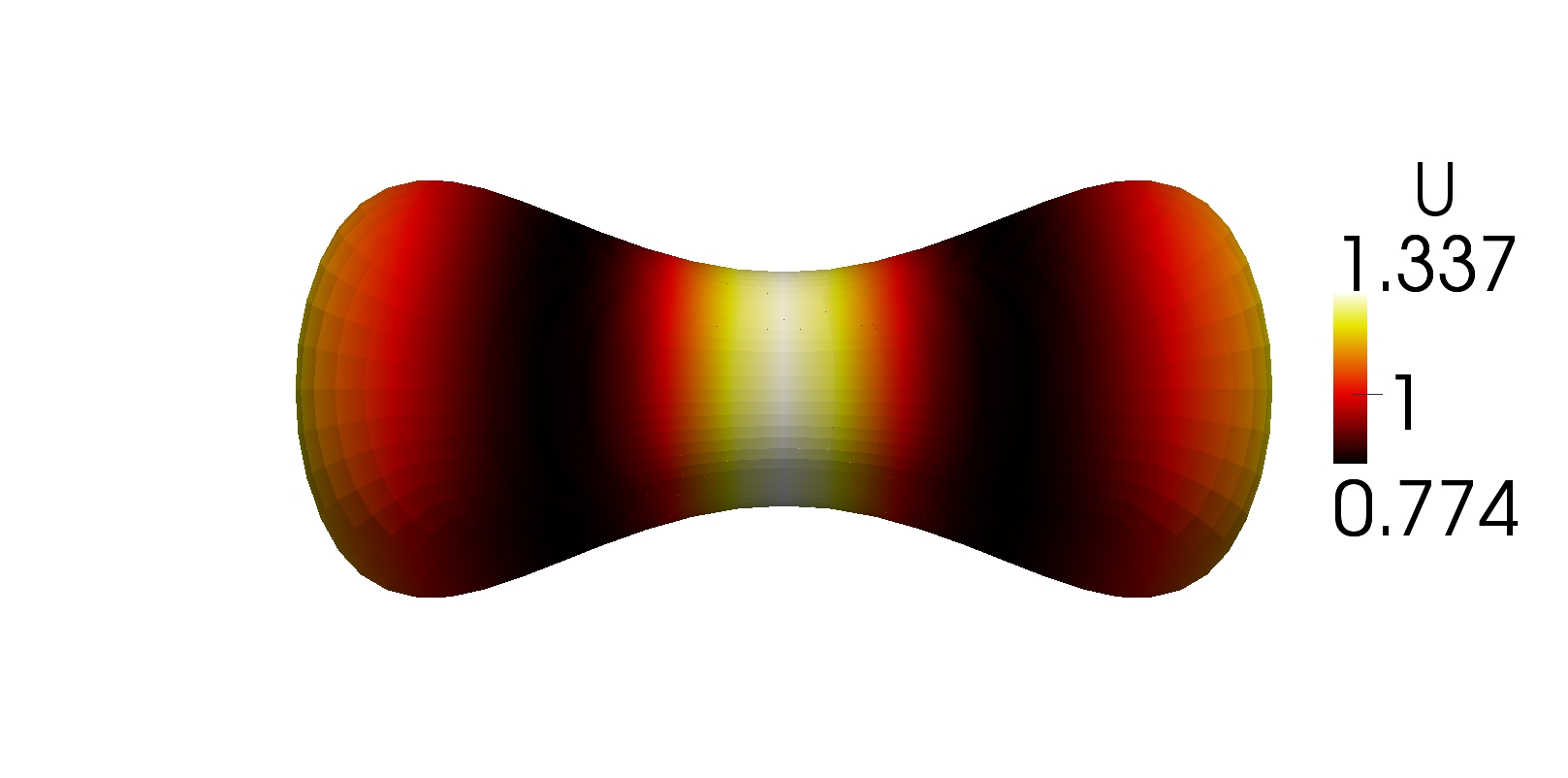}}
\caption{Converged $u$ solutions of system \eqref{eq:non} with Schnakenberg kinetics \eqref{eq:schnakkin} on a dumbell. Eigenvalues $\lambda_1,\lambda_2,\lambda_5,\lambda_6$ have been isolated, however since $\lambda_3\approx\lambda_4$ in (c) we see a linear combination of their eigenfunctions (Colour version online)}\label{fig:dumconverged}
\end{figure}

As a third example we consider the dumbbell shaped domain shown in Figure \ref{fig:meshes}\subref{fig:dumbbellgrid}. The solver for the eigenvalue problem on this mesh gives the output of eigenvalues and eigenfunctions shown in Figure \ref{fig:dumeigen}. The corresponding steady state solution with the parameters obtained by mode isolation are shown in Figure \ref{fig:dumconverged}.

\subsection{Example 4: Surface of a sphere} 
\begin{figure}
 \centering
\subfloat[The surface finite element solution with given  parameters $d=9$ and $\gamma=35$]{\includegraphics[width=55mm]{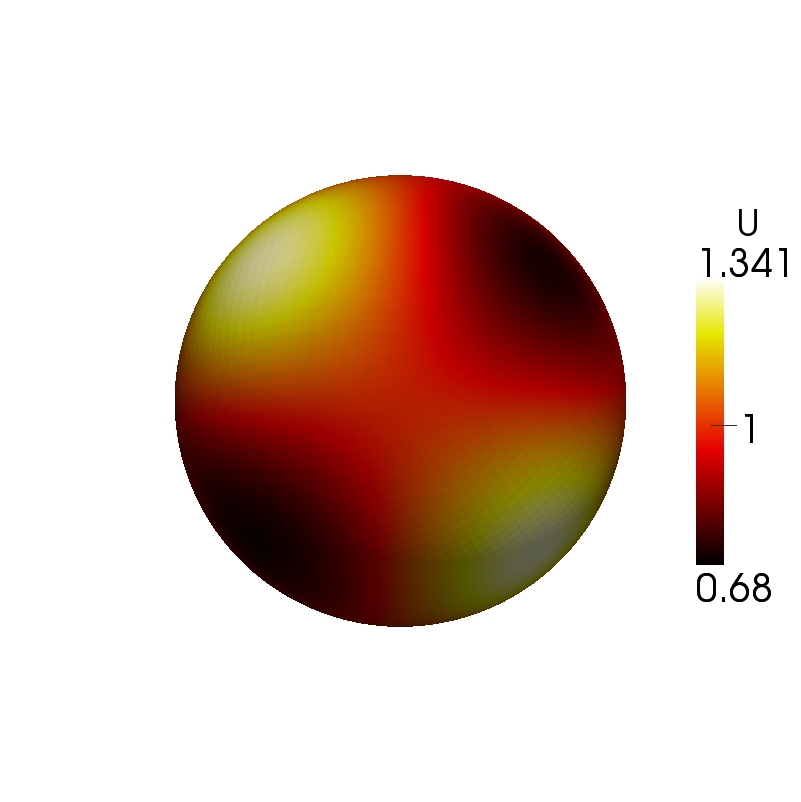}}\;\;\;
\subfloat[Numerically computed eigenfunction corresponding to eigenvalue $\lambda_9=12.0186$]{\includegraphics[width=55mm]{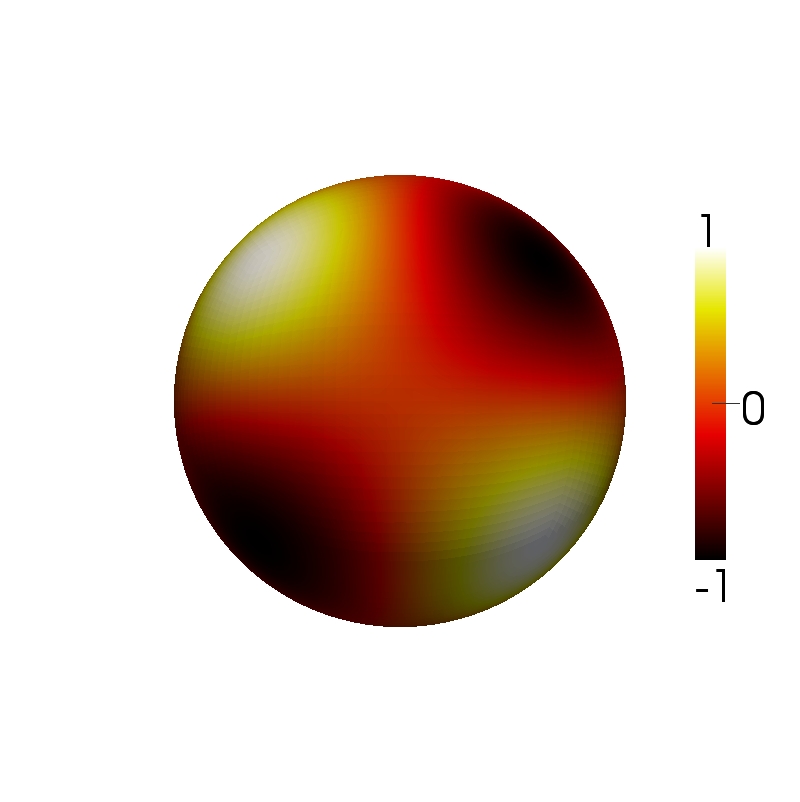}}
\caption{Mode isolation for the reaction-diffusion system with Schnakenberg kinetics on the surface of the sphere (Colour version online)}\label{fig:ss}
\end{figure}
In all the previous examples we considered bulk, volumetric domains. In this example we have a curved surface as the domain. This means using the Laplace Beltrami operator $\Delta_\Gamma$ instead of the Laplacian $\Delta$ in \eqref{eq:geneig} and \eqref{eq:non}.   To approximate solutions in this case, we employ the surface finite element method \citep{barreira,dziuk,dziukelliott,elliot2015,elliot,chung}. \\ 
The eigenpairs on the surface of the unit sphere can be found analytically and are well known and documented in \cite{chap} for example. The eigenfunctions are referred to as spherical harmonics. They are the restriction of the eigenfunctions \eqref{eq:efsph} to the surface. The eigenvalues are of the form $k^2=l(l+1)$, where $l$ is an integer, and the eigenfunctions are
\begin{equation}
w_l^m(\theta,\phi)=A_l^me^{im\phi}P_l^m(\cos\theta),
\end{equation}
where $m$ and $P_l^m$ are as in Section \ref{sec:sphere}. Therefore we can test the performance of the eigenvalue problem solver with this example. Using the eigenvalue solver on an approximated mesh of the surface of the sphere we obtain the following output of the first 30 eigenvalues computed to 4 decimal places
\begin{table}[H]
\begin{center}
 \begin{tabular}{ c c c c c c c c}
$k_h^2$    = & 2.0014, & 2.0014, & 2.0014, & & \\ & 6.00664, & 6.00664, & 6.00671, & 6.0085, & 6.00857, & & \\ & 12.0186, & 12.0224, & 12.0224, & 12.023, & 12.0279, & 12.0284, & 12.0284,\\ & 20.0484, & 20.0484, & 20.0622, & 20.0622, & 20.0717, & 20.0717, & 20.0749, \\ & 30.1043, & 30.1043, &  30.1102, & 30.1102, & 30.1523, & 30.1523, & 30.1591.
\end{tabular}
\end{center}
\end{table}
As expected these are the first 5 values of the form $k^2=l(l+1)$ with $l=1,2,3,4,5$. The values are not exact because the mesh is an approximation of the actual surface of the sphere. The eigenfunctions are analogous to those detailed in Section \ref{sec:sphere} restricted to the boundary. This shows that the eigenvalue solver gives the required output. Since the results are shown in Section \ref{sec:sphere} we only show one example of mode isolation in Figure \ref{fig:ss}.

\subsection{Example 5: "fish" surface}
\begin{figure}
 \centering
\subfloat[$d=8.9$, $\gamma=130$]{\includegraphics[trim = 0 35mm 0 50mm, clip,height=36mm]{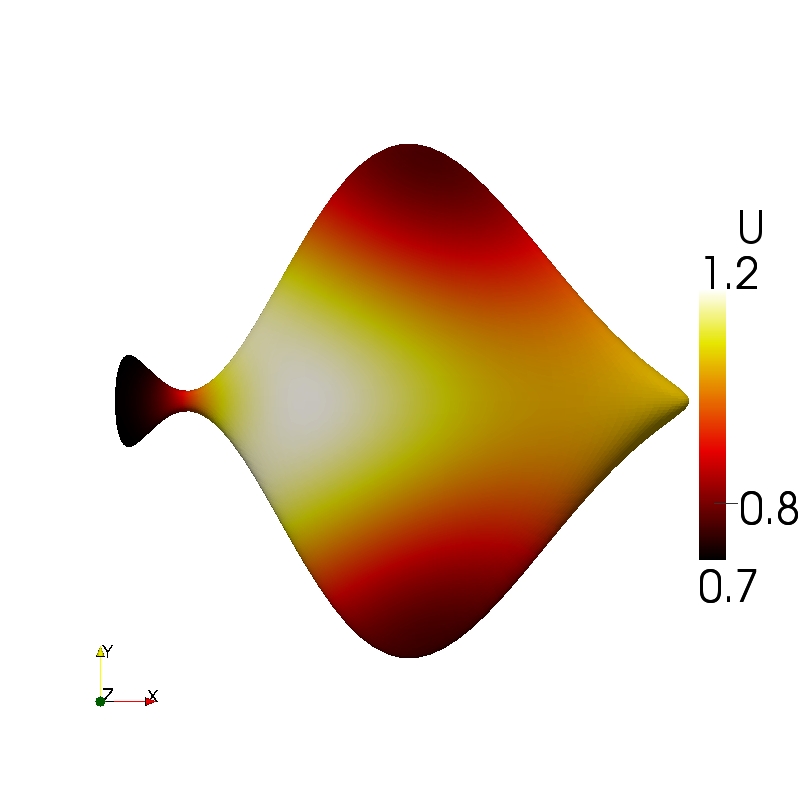}}\;
\subfloat[$\lambda_5=40.18$]{\includegraphics[trim = 0 35mm 0 50mm, clip,height=36mm]{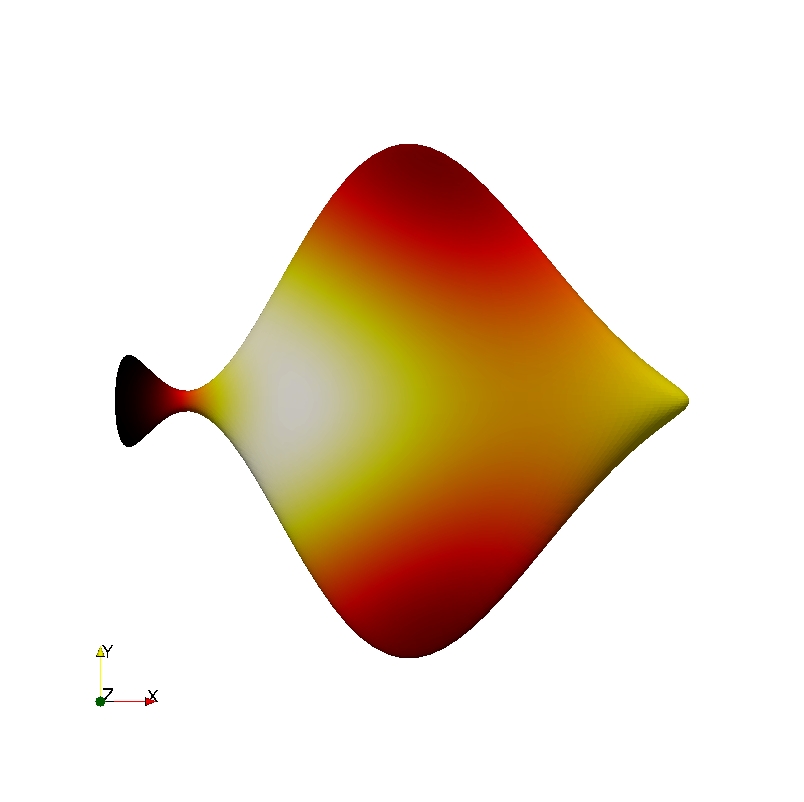}}\\
\subfloat[$d=8.58$, $\gamma=240$]{\includegraphics[trim = 0 35mm 0 50mm, clip,height=36mm]{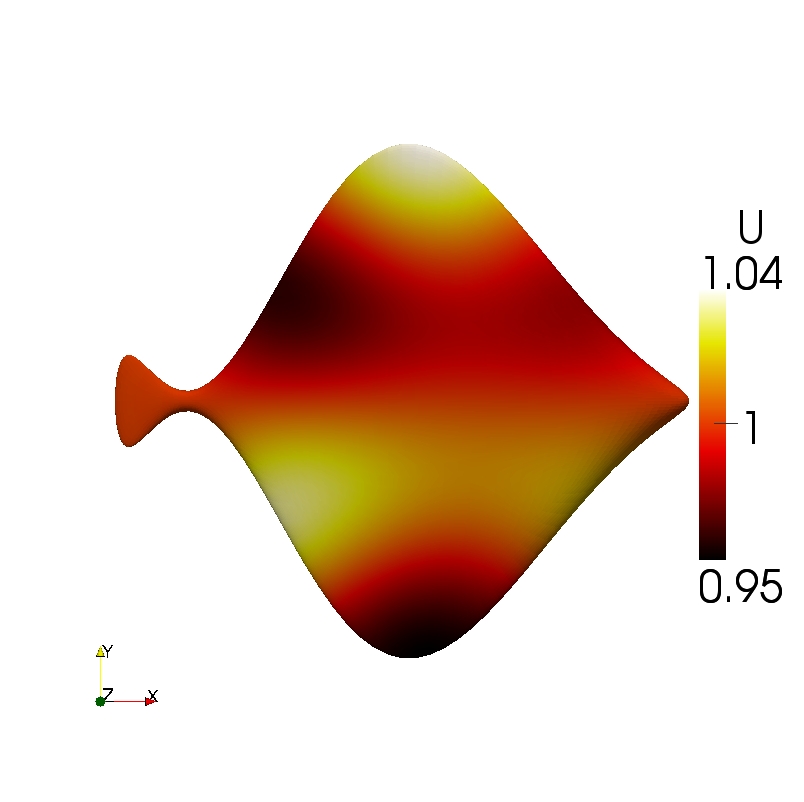}}\;
\subfloat[$\lambda_{10}=79.56$]{\includegraphics[trim = 0 35mm 0 50mm, clip,height=36mm]{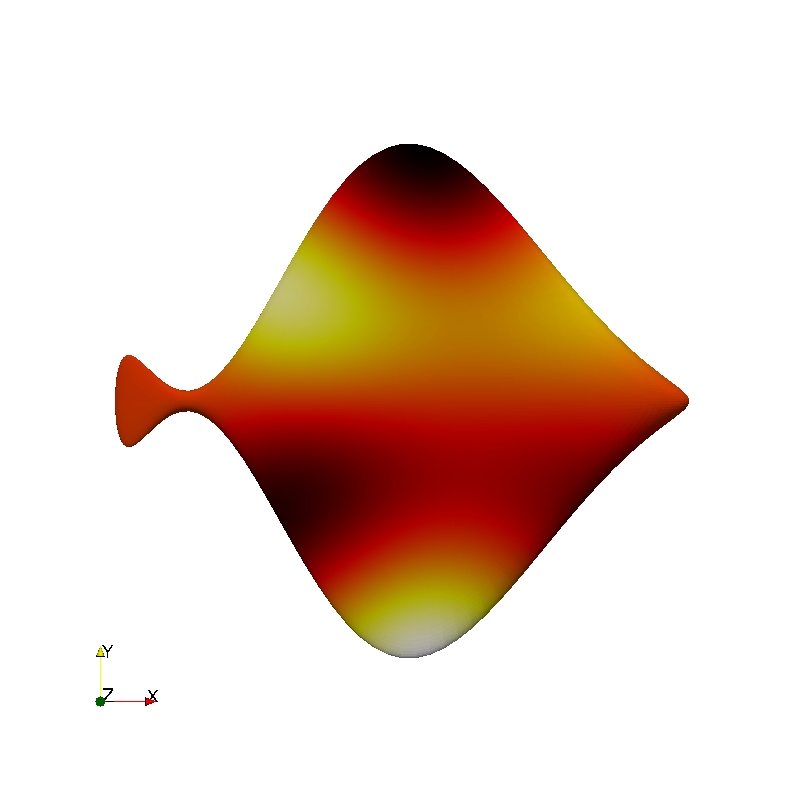}}\\
\subfloat[$d=8.58$, $\gamma=400$]{\includegraphics[trim = 0 35mm 0 50mm, clip,height=36mm]{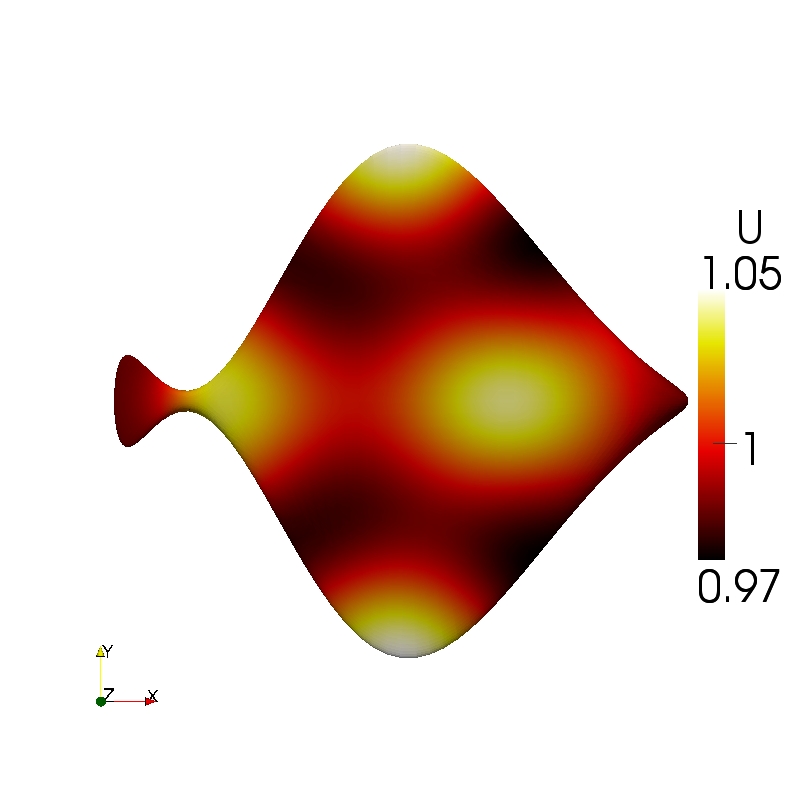}}\;
\subfloat[$\lambda_{15}=134.73$]{\includegraphics[trim = 0 35mm 0 50mm, clip,height=36mm]{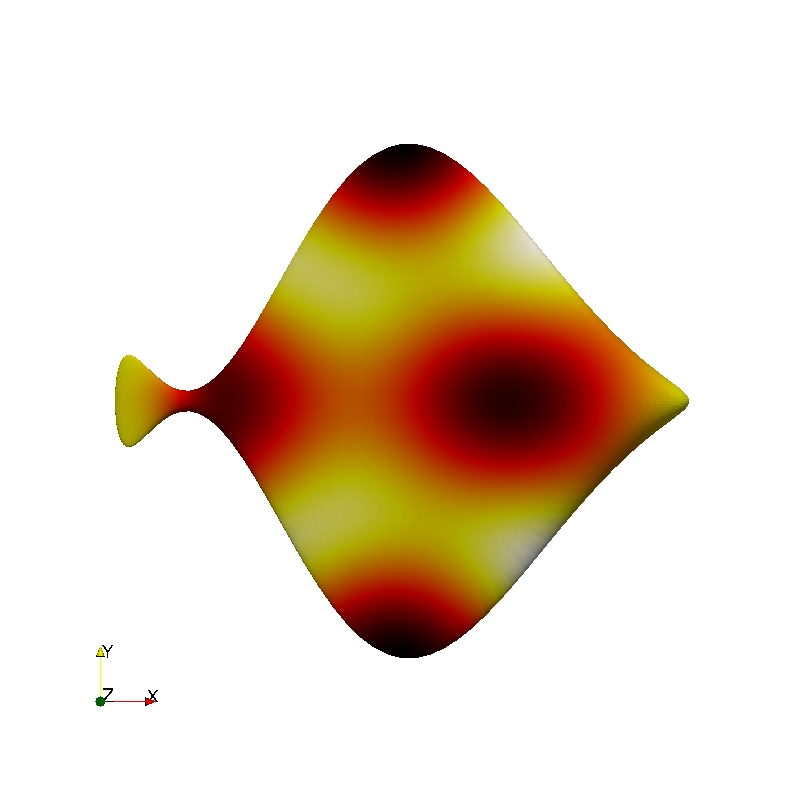}}\\
\subfloat[$d=8.58$, $\gamma=510$]{\includegraphics[trim = 0 35mm 0 50mm, clip,height=36mm]{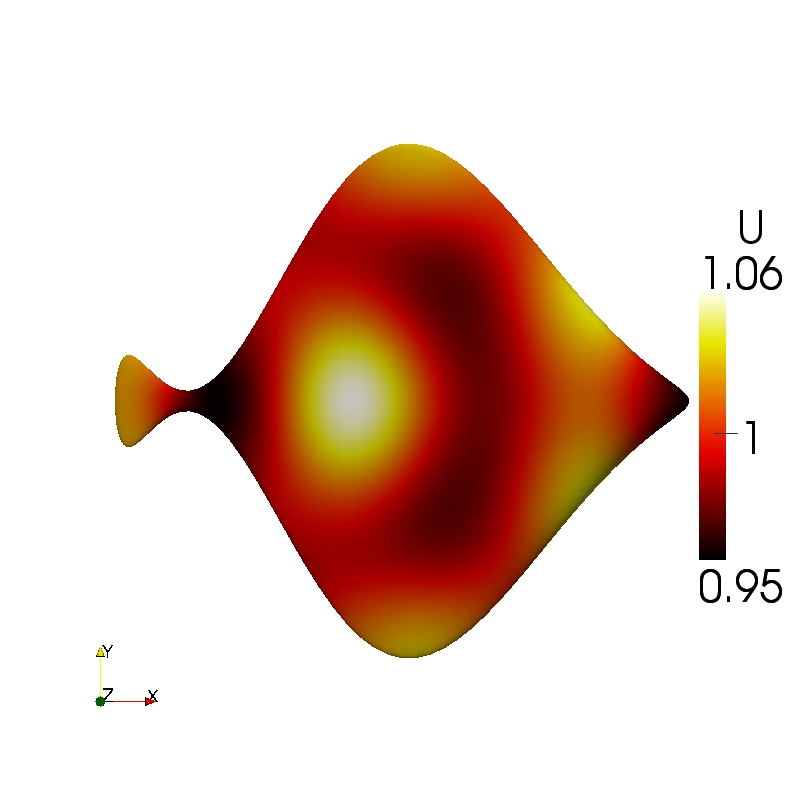}}\;
\subfloat[$\lambda_{19}=175.98$]{\includegraphics[trim = 0 35mm 0 50mm, clip,height=36mm]{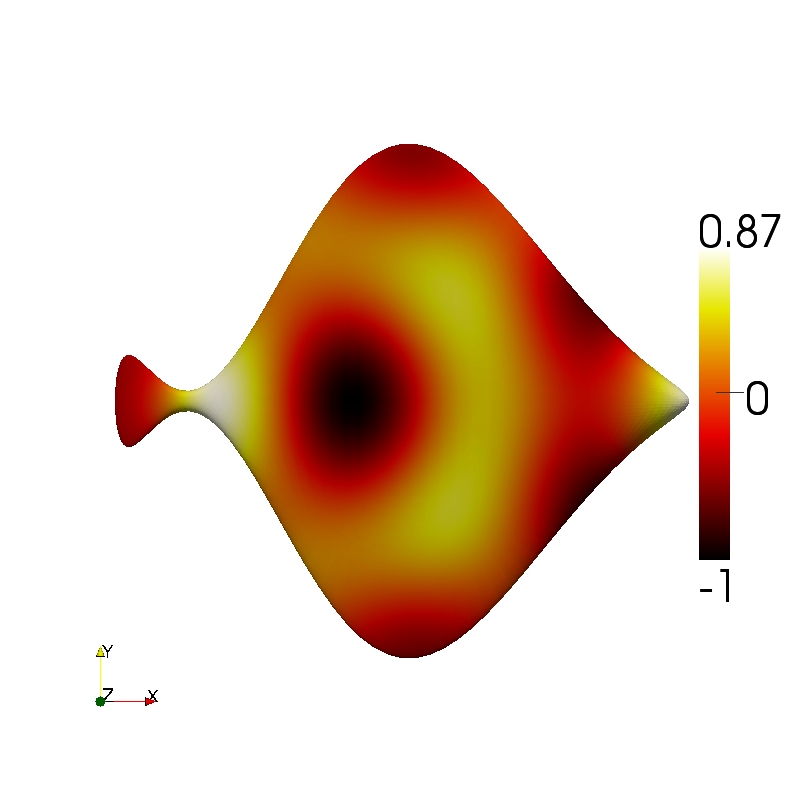}}
\caption{Surface finite element solutions corresponding to the $u$ species of the reaction-diffusion system with Schnakenberg kinetics with the given parameters on the left and numerically computed eigenfunctions corresponding to the given eigenvalue on the right (Colour version online)}\label{fig:smsols}
\end{figure}

We now consider a smooth surface on which no analytical expression for the eigenpairs is available, the surface is taken to be diffeomorphic to the sphere and is shown in Figure \ref{fig:meshes}\subref{fig:smmesh}, it is meant to (very loosely) mimic the shape of a fish. We found the first 100 eigenpairs then chose several to isolate. These are shown in Figure \ref{fig:smsols}. Various patterns are observed including stripes, spots and concentric rings.

\begin{figure}
 \centering
\subfloat[$\lambda_{4(open)}=54.43$]{\includegraphics[height=24mm]{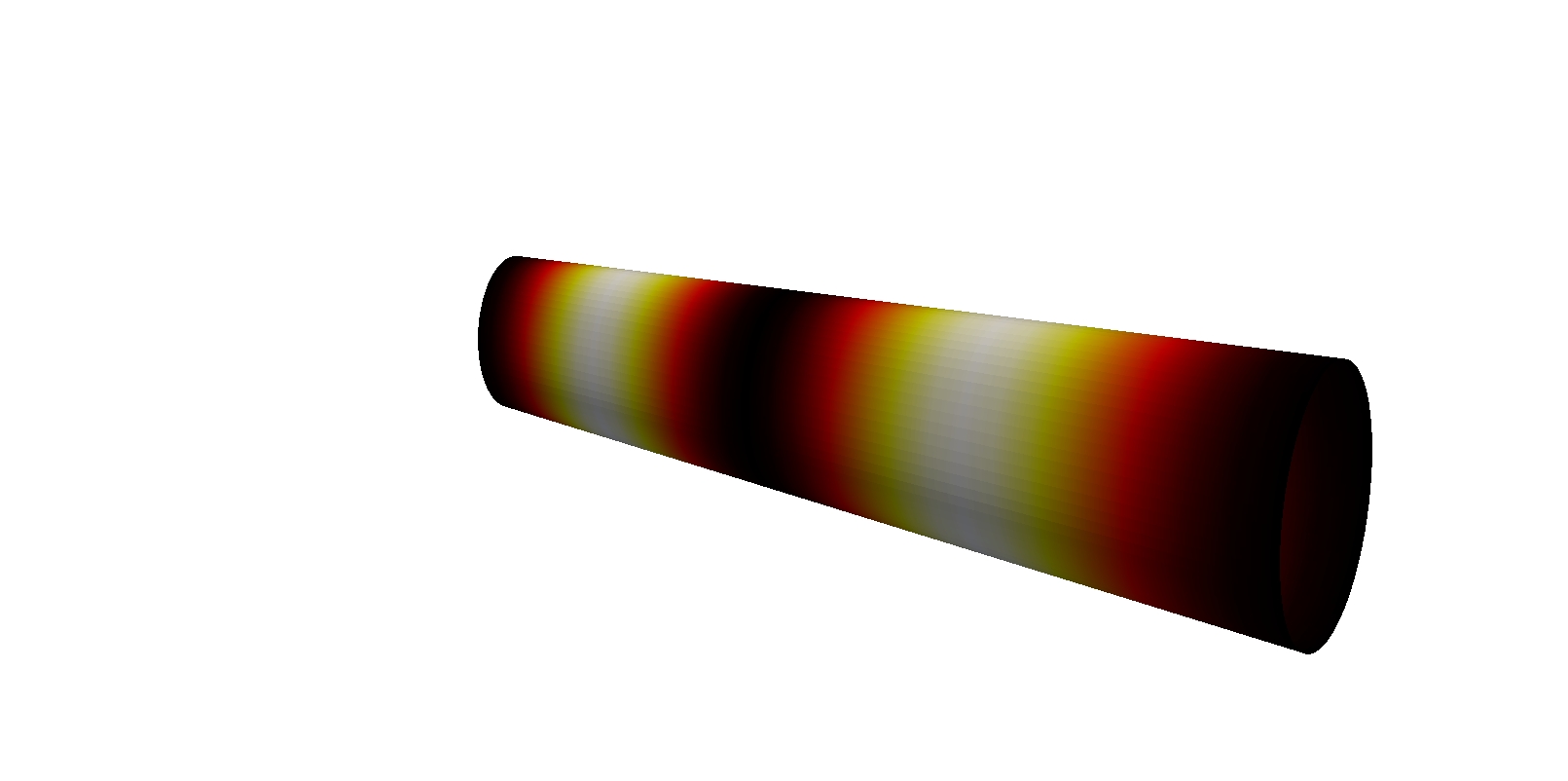}}\;
\subfloat[$\lambda_{4(closed)}=44.94$]{\includegraphics[height=24mm]{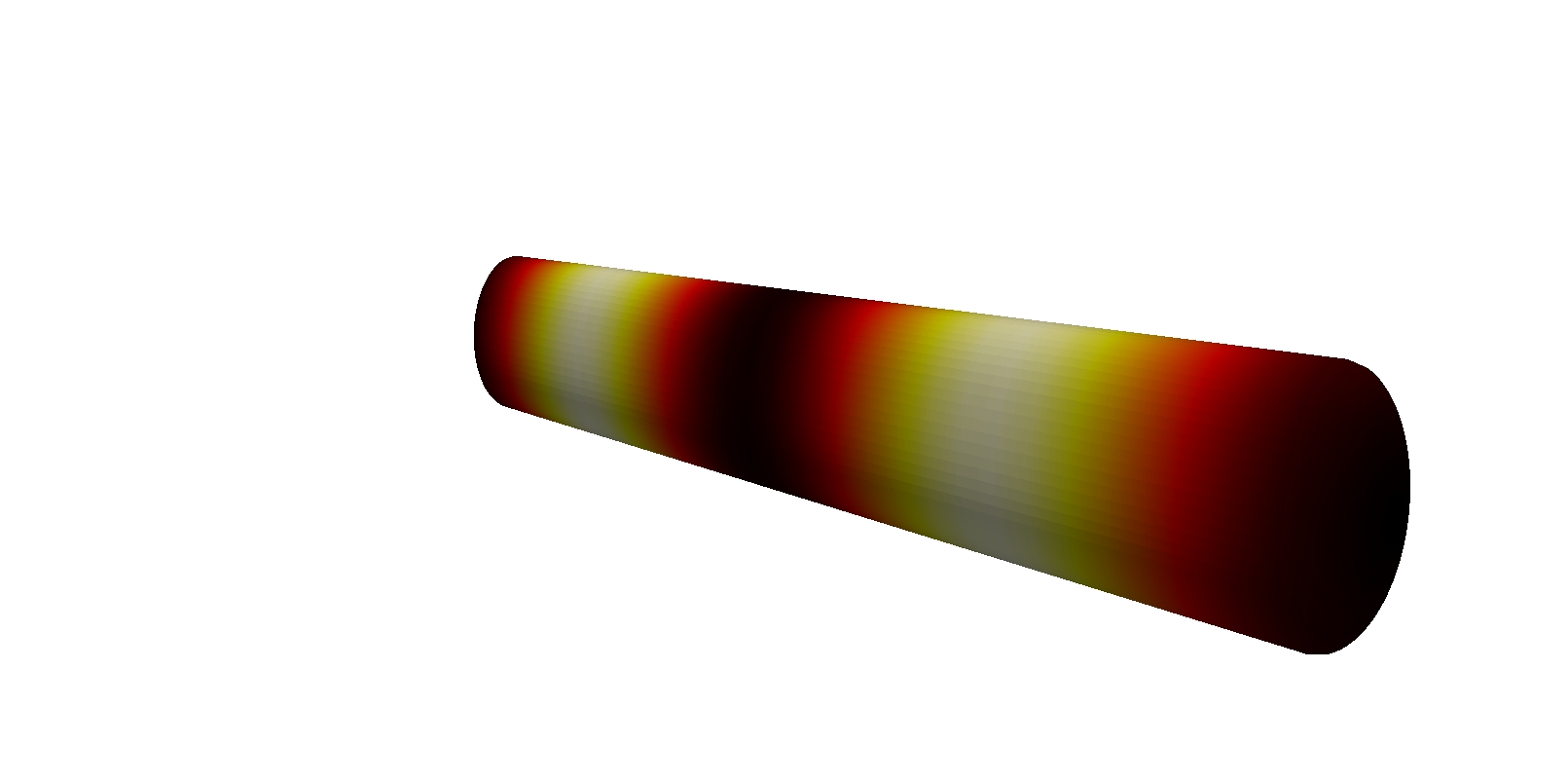}}\\
\subfloat[$\lambda_{23(open)}=253.69$]{\includegraphics[height=23mm]{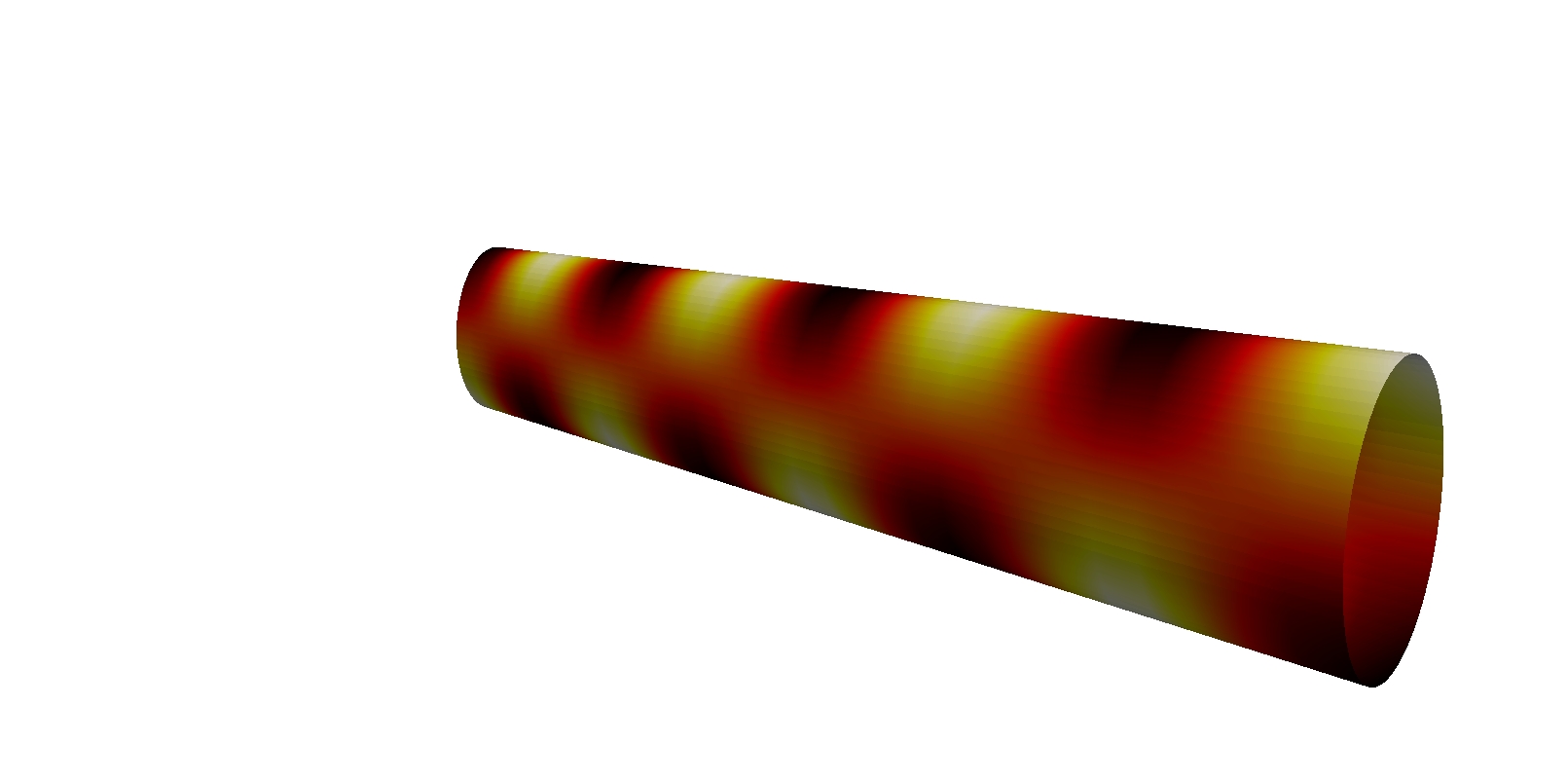}}\;
\subfloat[$\lambda_{25(closed)}=257.54$]{\includegraphics[height=24mm]{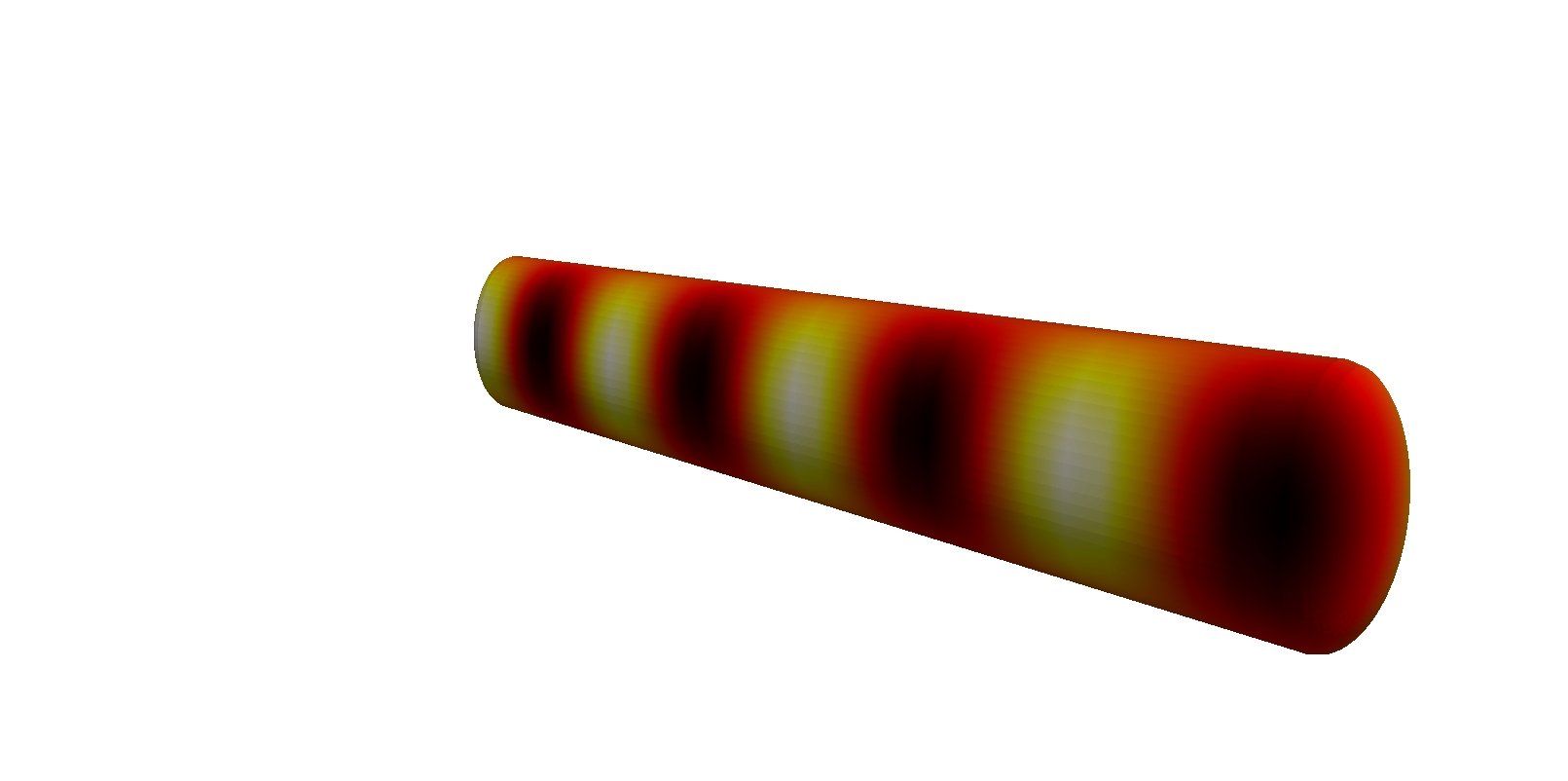}}\\
\subfloat[$\lambda_{24(open)}=253.73$]{\includegraphics[height=23mm]{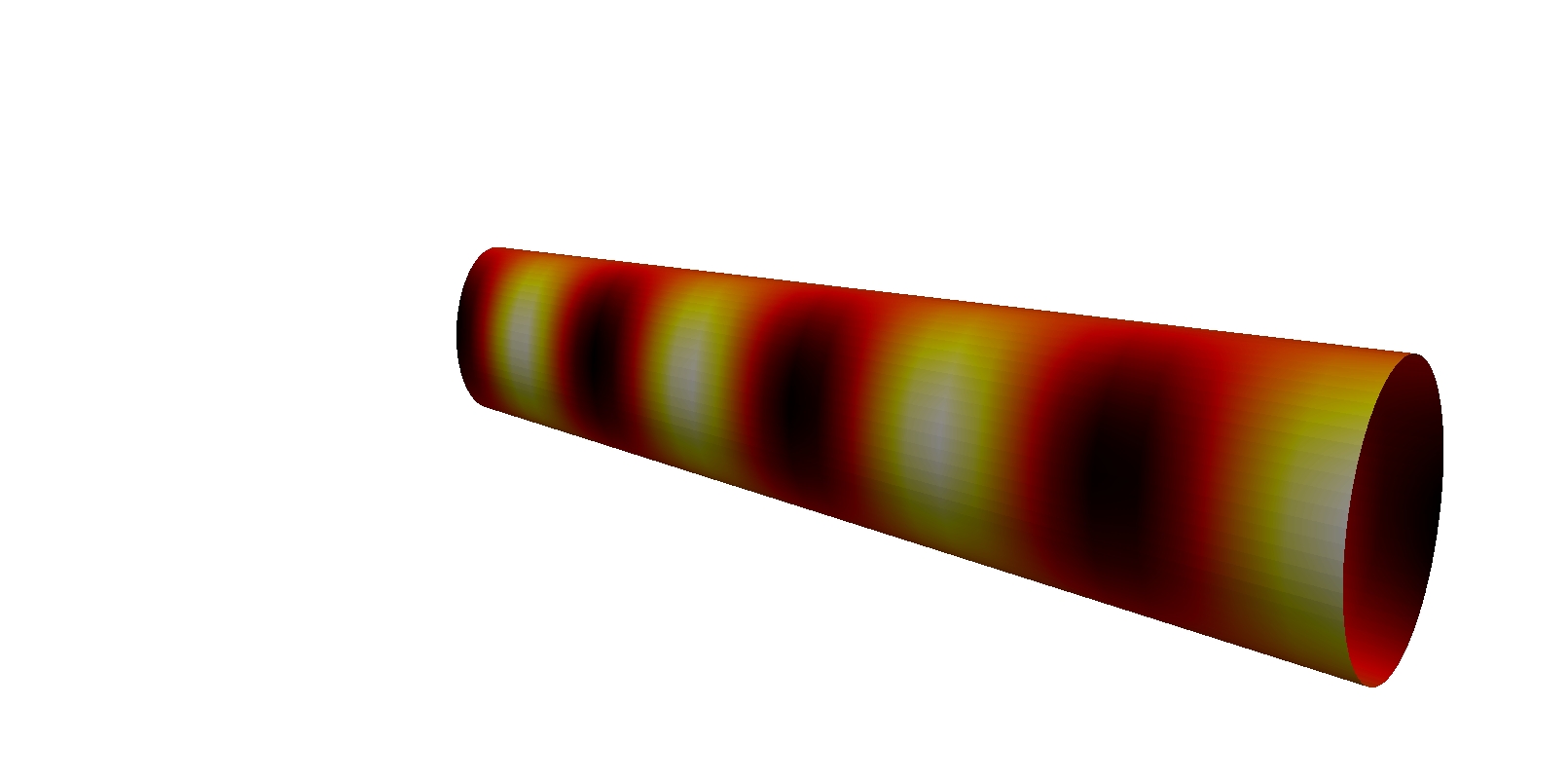}}\hspace{50mm}
\caption{Eigenfunctions of the Laplace-Beltrami operator on the "eel" shape with the corresponding eigenvalue. The left column shows the surface without a boundary and the right has a boundary. Note that, although the eigenfunctions are different, $\lambda_{23}\approx\lambda_{24}$ (Colour version online)}\label{fig:eelefs}
\subfloat[$d=8.8$, $\gamma=140$]{\includegraphics[height=24mm]{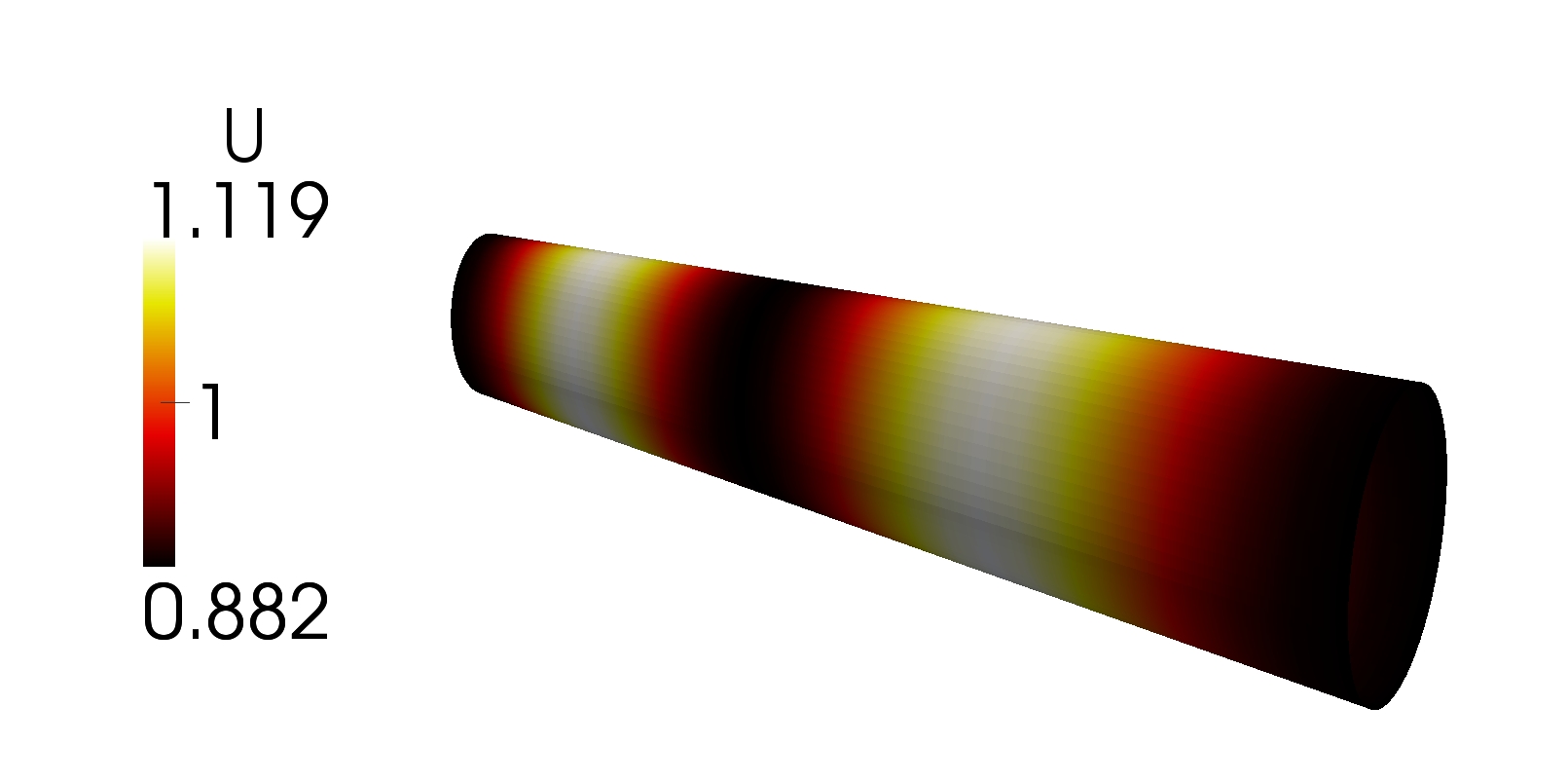}}\;
\subfloat[$d=8.8$, $\gamma=140$]{\includegraphics[height=24mm]{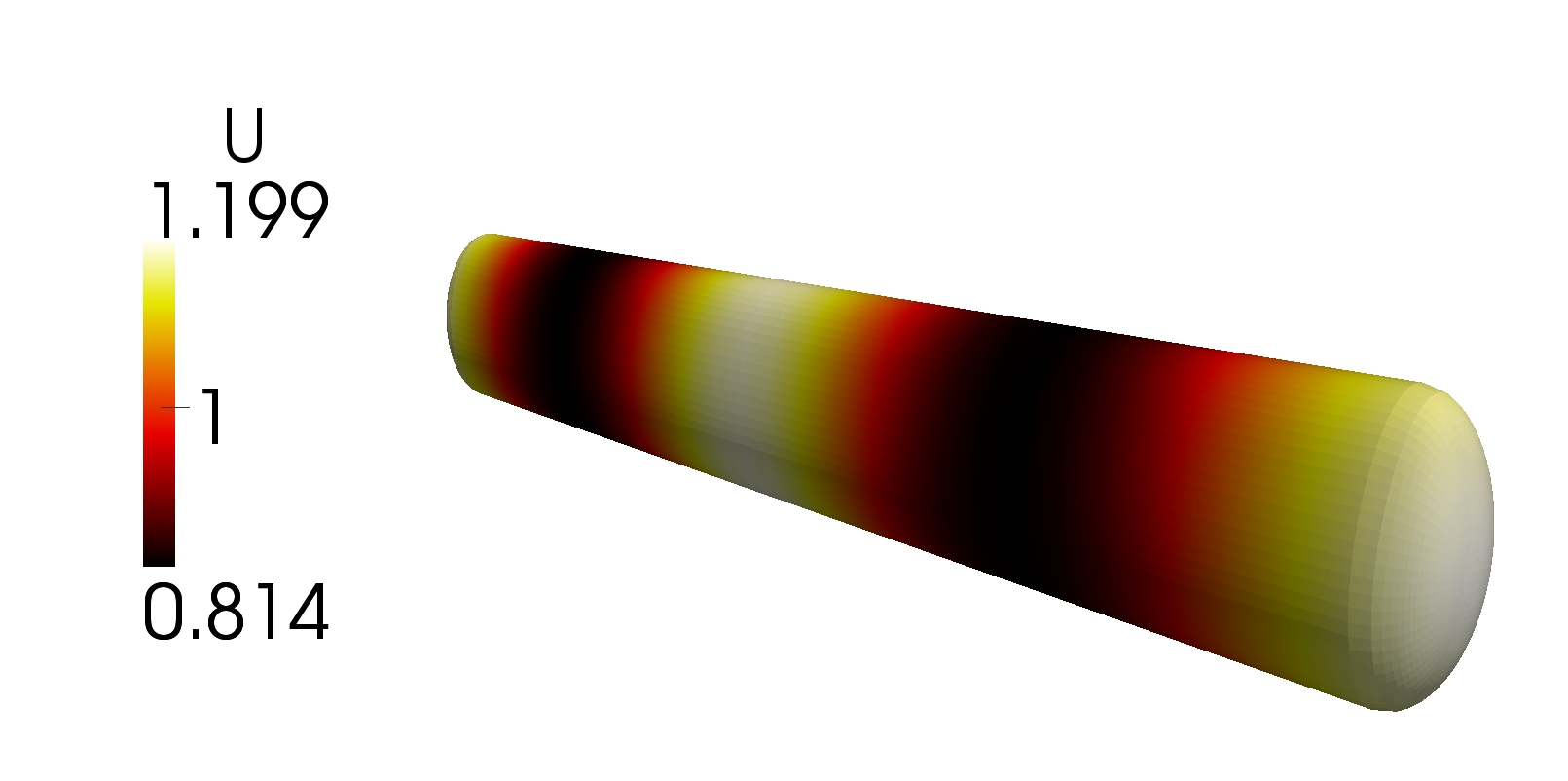}}\\
\subfloat[$d=8.6$, $\gamma=750$]{\includegraphics[height=24mm]{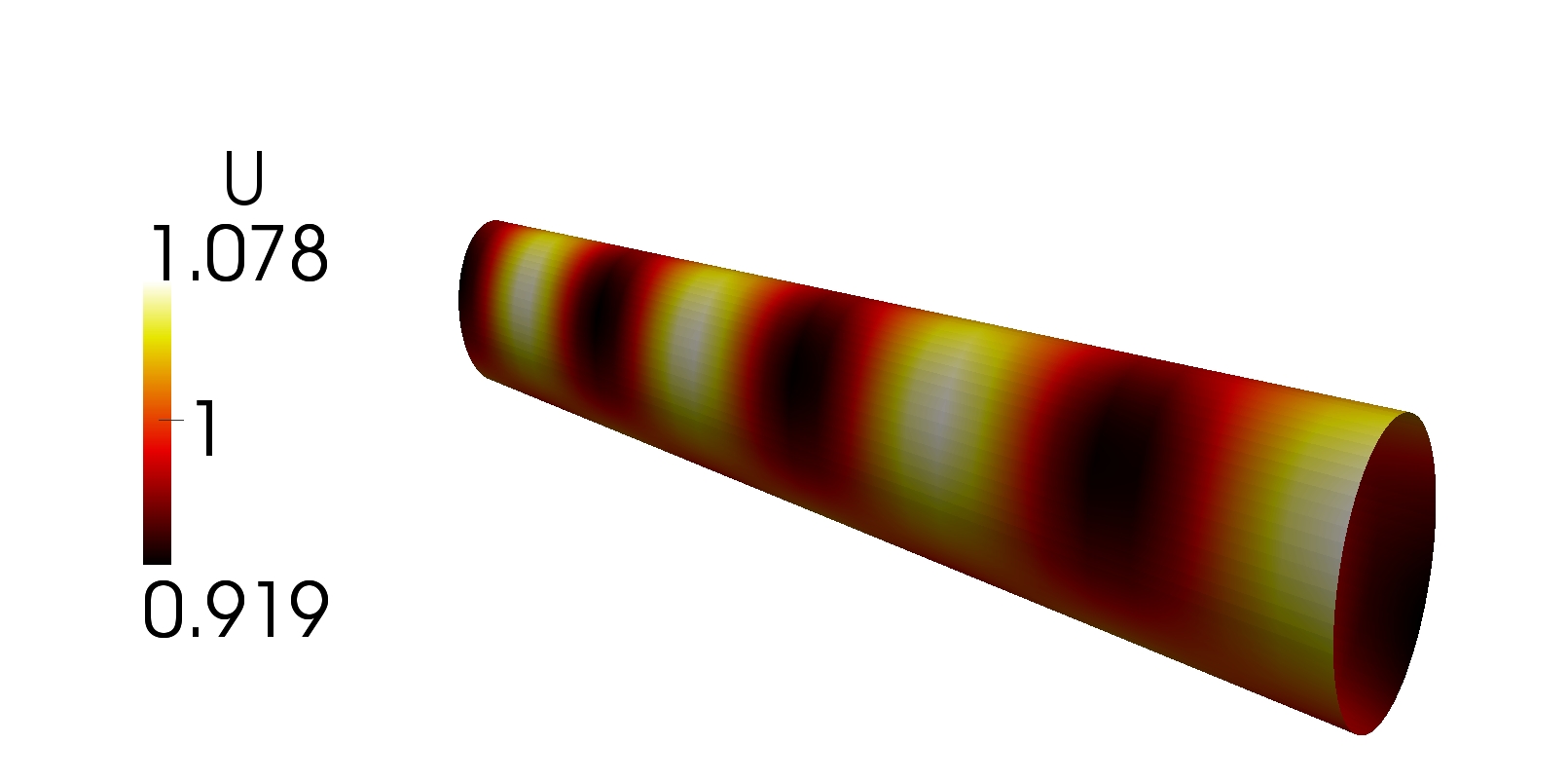}\label{fig:o750}}\;
\subfloat[$d=8.6$, $\gamma=750$]{\includegraphics[height=24mm]{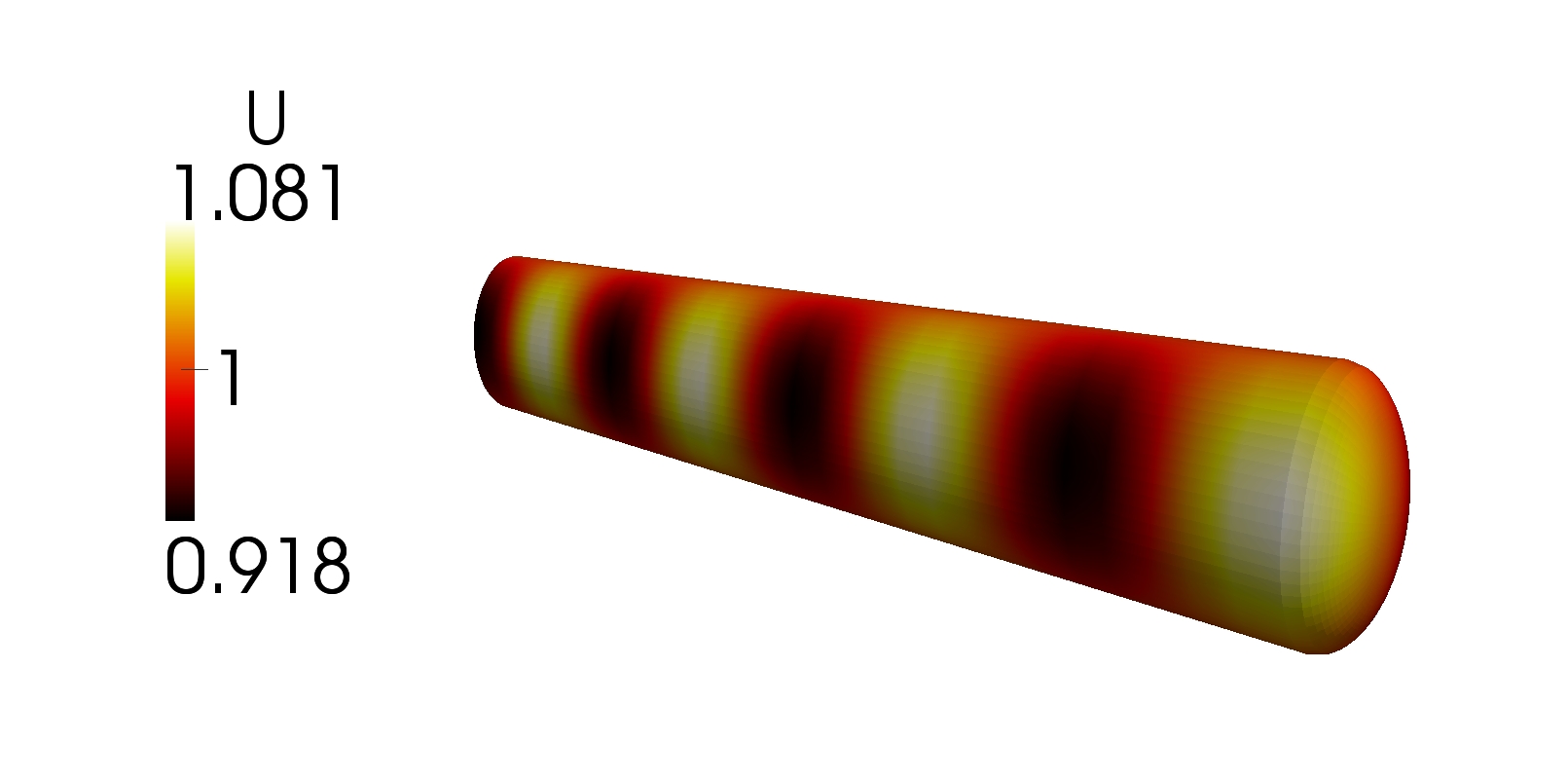}}\\
\caption{Converged solutions corresponding to the $u$ species of the reaction-diffusion system with Schnakenberg kinetics on the surface of an eel. The surfaces on the right have a boundary whereas those on the left do not. We find that using the same parameter values on both surfaces gives very similar results (Colour version online)}\label{fig:eelsols}
\end{figure}

\subsection{Example 6 and 7 "eel" shapes}
When computing on surfaces, one has to consider whether or not the surface has a boundary. In papers modelling fish or eel patterns (see for example \cite{venkataraman}), a surface with a boundary is often used. To investigate whether having a boundary is significant in this example we consider a surface  with  and  without  boundary. We see that the eigenvalues and eigenfunctions are very similar and it is possible to isolate similar patterns using the same parameter values.

\section{Conclusion and further challenges}
In this paper we have considered reaction-diffusion systems and have presented a framework for isolating particular spatially inhomogeneous patterns. The method involves finding eigenpairs  of the Laplacian and computing parameters  such that when the reaction-diffusion system is solved numerically, only patterns analogous to a particular eigenfunction will grow. In previous works the eigenvalue problem is solved analytically whereas in this paper both the eigenvalue problem and the reaction-diffusion system are solved using the finite element method. Advances in numerical software mean that we can find 100 eigenpairs in a few minutes and we have demonstrated that these eigenpairs match analytical results. The approach is shown to work for 3 different examples of nonlinear reaction kinetics and on a variety of domains and surfaces. 
In summary, the main observations are:
\begin{itemize}
 \item Mode isolation is straightforward for low values of $k^2$ but can become slightly more difficult for higher values of $k^2$. This is due to the approximation of the nonlinear terms and clustering of the eigenvalues of a linear problem.
 \item When two or more eigenvalues are clustered close to each other it becomes difficult to isolate them computationally. If two or more eigenvalues are in the permissible range then the inhomogeneous steady state could be a linear combination of the corresponding eigenfunctions.
 \item We display an example of two surfaces where pattern formation appears to be robust despite the fact one has a boundary while the other does not. An interesting investigation would be to see if this can be true for other geometries. Note that this is only the case for zero-flux boundary conditions. Imposing Dirichlet or Robin-type boundary conditions would result in substantially different patterns.
\end{itemize}
In this paper we have only considered stationary domains/volumes and surfaces. However the domains of biological processes generally evolve with time \citep{barreira,elliot,lakkis2013,madzvamuse2006,venkataraman}. This adds more complexity to solving the reaction-diffusion systems. An interesting and natural extension of this work would be to introduce domain growth and surface evolution. For this extension, studies on the effects of initial conditions would also be worthwhile.

\section*{Data management}
All the computational data output is included in the present manuscript.

\section*{Acknowledgements}
This work (LM) was supported by an EPSRC Doctoral Training Centre Studentship through the University of Sussex. CV and AM acknowledge support from the Leverhulme Trust
Research Project Grant (RPG-2014-149) and the EPSRC grant (EP/J016780/1). This research was partly undertaken
whilst LM, CV and AM were participants in the Isaac Newton Institute Program,
Coupling Geometric PDEs with Physics for Cell Morphology, Motility and Pattern
Formation. This work (AM) has received funding from the European Union Horizon
2020 research and innovation programme under the Marie Sklodowska-Curie grant
agreement (No 642866). AM was partially supported by a grant from the Simons
Foundation. LM  acknowledges the support from the University of Sussex ITS for
computational purposes. 

\newpage 
\bibliographystyle{plainnat}
 \bibliography{mybib}{}

\begin{appendix}
\section{The Krylov-Schur algorithm}\label{krsch}
The Krylov-Schur algorithm was introduced by \cite{stewart} and is an alternative to the method of \cite{arnoldi}. The aim of the algorithm is to compute a number of eigenpairs of a given square matrix A.\\
The basic {\bf Arnoldi algorithm} has input matrix $A$ and initial vector $\bv_1$ of norm 1 ($\bv_j$ will make up the columns of an $n\times m$ matrix $\bV_m$) and output $\bV_m,\bH_m,\bfn,\beta$ such that 
\begin{equation}
 \bA\bV_m=\bV_m\bH_m+fe_m^*,\; \beta=||\bfn||_2.                                                                                                                                                                                                                                                                                                                                                                                                                \end{equation}
A Krylov decompostion is a generalised version of this and is given by
\begin{equation}
\bA\bV_m= \bV_m\bB_m+\bv_{m+1}\bb^*_{m+1},
\end{equation}
where $\bB_m$ is not necessarily upper Hessenberg and $\bb^*_{m+1}$ is an arbitrary vector.
The Krylov-Schur method is described in the SLEPc Technical Report \citep{str7} as follows
\blockquote{
Input: Matrix $\bA$, initial vector $\bv_1$, and dimension of the subspace $m$\\
Output: A partial Schur decompostion $\bA\bV_{1:k}=\bV_{1:k}\bH_{1:k,1:k}$
\begin{itemize}
\item Normalize $\bv_1$
\item Initialize $\bV_m=[v_1],\quad k=0,\quad p=0$
\item Restart loop
      \begin{itemize}
\item[$\diamond$] Perform $m-p$ steps of Arnoldi with deflation
\item[$\diamond$] Reduce $\bH_m$ (part of the output of the Arnoldi algorithm) to (quasi-)triangular form, $\bH_m \leftarrow \bU_1^*\bH_m\bU_1$
\item[$\diamond$] Sort the $1\times 1$ or $2\times 2$ diagonal blocks: $\bH_m \leftarrow \bU_2^*\bH_m\bU_2$
\item[$\diamond$] $\bU=\bU_1\bU_2$
\item[$\diamond$] Compute eigenpairs of $\bH_m$, $\bH_m\by_i=\by_i\theta_i$
\item[$\diamond$] Compute residual norm estimates, $\tau_i=\beta |\be^*_m\by_i|$
\item[$\diamond$] $\bV_m\leftarrow \bV_m\bU$
\item[$\diamond$] Exit if enough converged eigenpairs, otherwise lock newly converged vectors
\item[$\diamond$] Choose $p$ ($k$ (the number of currently converged eigenpairs) $<p<m$) and set $\tilde{\bv}_{p+1}=\bv_{m+1}$
\item[$\diamond$] Compute $\bb_w$ (the leading subvector of $\bb^*_{m+1}\bU$) and insert in the appropriate positions of $\bH_p$
      \end{itemize}
\item end
\end{itemize}
}
If the eigenpairs of $\bH_{1:k,1:k}$ (ie solutions of $\bH\by=\theta \by$) are $(\theta_i, \by_i)$ then the approximate eigenvalues of $\bA$ are $\lambda_i=\theta_i$ and eigenvectors are $\bx_i=\bV_{1:k}\by_i$.
In our problem we have \eqref{eq:ep} (the generalised eigenvalue problem) instead of $\bA \bx=\lambda \bx$, and here one works with a spectral transformation $\bT_S=\bM^{-1}\bA$ or $\bT_{SI}=(\bA-\sigma \bI)\bM$ instead of $\bA$.
\end{appendix}

\newpage
\appendix

\end{document}